\documentclass[a4paper,reqno]{amsart}
\usepackage{amssymb,amsthm,amscd,textcomp,url}
\usepackage[all,ps,dvips]{xy}

\newtheoremstyle{statement}%
{}{}{\slshape}{}{\upshape\bfseries}{.}{ }{}
\theoremstyle{statement}
\newtheoremstyle{example}%
{}{}{}{}{\upshape\bfseries}{.}{.5em}{}
\theoremstyle{statement}
\newtheorem{theorem}{Theorem}[section]
\newtheorem{proposition}[theorem]{Proposition}
\theoremstyle{example}
\newtheorem{example}{Example}
\theoremstyle{remark}
\newtheorem{remark}{Remark}
\newtheorem*{coordinates}{Coordinates}

\newcommand{\fl}{\hspace*{-6pc}}

\newcommand{\J}{J}
\newcommand{\JJ}{J_h}
\newcommand{\Se}{\Gamma}
\newcommand{\sd}[2]{\{\,{#1}\,\mid\,{#2}\,\}}
\newcommand{\cprime}{\/{\mathsurround=0pt$'$}}

\DeclareMathAlphabet{\bi}{OML}{cmm}{b}{it}

\newcommand{\F}{\mathcal{F}}
\newcommand{\C}{\mathcal{C}}
\newcommand{\Ev}{\bi{E}}
\newcommand{\R}{\mathbb{R}}
\newcommand{\X}{\mathcal{X}}
\newcommand{\Xc}{\X_{\scriptscriptstyle{\mathcal{C}}}}
\newcommand{\La}{\Lambda}
\newcommand{\Lah}{\La_h}
\newcommand{\Lac}{\La_{\scriptscriptstyle{\mathcal{C}}}}
\newcommand{\abs}[1]{\left\vert#1\right\vert}

\let\CDiff\C
\newcommand{\CDiffsk}{\CDiff^{\mathrm{sk}}}
\newcommand{\CDiffskalt}{\CDiff^{\mathrm{sk\,*}}}

\DeclareMathOperator{\Hom}{Hom}
\DeclareMathOperator{\im}{im}
\DeclareMathOperator{\sym}{sym}
\DeclareMathOperator{\cl}{cl}
\DeclareMathOperator{\T}{\mathcal{T}}
\DeclareMathOperator{\Ts}{\mathcal{T}^*}

\let\phi=\varphi
\let\kappa=\varkappa
\let\epsilon=\varepsilon

\newcommand*{\rmd}{\mathinner{\!}\mathrm{d}}
\newcommand{\rmdh}{\rmd_h}
\newcommand{\rmdc}{\rmd_{\scriptscriptstyle{\mathcal{C}}}}

\newcommand{\lshad}{[\![}
\newcommand{\rshad}{]\!]}

\newcommand{\rme}{\mathrm{e}}
\newcommand{\id}{\mathrm{id}}

\newcommand{\eval}[2][\right]{\relax
\ifx#1\right\relax \left.\fi#2#1\vert}

\begin{document}
\title{Geometry of jet spaces and integrable systems}
\author{Joseph Krasil{\cprime}shchik}
\address{Independent University of Moscow,
B.~Vlasevsky~11,
119002 Moscow,
Russia}
\email{josephk@diffiety.ac.ru}

\author{Alexander Verbovetsky}
\address{Independent University of Moscow,
B.~Vlasevsky~11,
119002 Moscow,
Russia}
\email{verbovet@mccme.ru}

\begin{abstract}
  An overview of some recent results on the geometry of partial differential
  equations in application to integrable systems is given.  Lagrangian and
  Hamiltonian formalism both in the free case (on the space of infinite jets)
  and with constraints (on a PDE) are discussed. Analogs of tangent and
  cotangent bundles to a differential equation are introduced and the
  variational Schouten bracket is defined. General theoretical constructions
  are illustrated by a series of examples.
\end{abstract}

\maketitle

\section*{Introduction}
\label{sec:gjsis_intro:introduction}

The main task of this paper is to overview a series of our results achieved
recently in understanding integrability properties of partial differential
equations (PDEs) arising in mathematical physics and
geometry~\cite{KerstenKrasilshchikVerbovetsky:HOpC,KerstenKrasilshchikVerbovetsky:InCSSREvDEq,KerstenKrasilshchikVerbovetsky:MAmEqHSSRH,KerstenKrasilshchikVerbovetsky:NHSSRHNApApNSKEq,KerstenKrasilshchikVerbovetsky:GSDBTEq,GolovkoKrasilshchikVerbovetsky:InCHEqIn,GolovkoKrasilshchikVerbovetsky:VPNSPDEq,IgoninKerstenKrasilshchikVerbovetskyVitolo:VBGPDE,KerstenKrasilshchikVerbovetskyVitolo:HSGP,KerstenKrasilshchikVerbovetskyVitolo:InKD}. These
results are essentially based on the geometrical approach to PDEs developed
since the 1970s by A.~Vinogradov and his school
(see~\cite{KrasilshchikLychaginVinogradov:GJSNPDEq,KrasilshchikVinogradov:SCLDEqMP,KrasilshchikVerbovetsky:HMEqMP,Vinogradov:CAnPDEqSC,KrasilshchikKersten:SROpCSDE}
and references therein). The approach treats a PDE as an
(infinite-dimensional) submanifold in the space~$\J^\infty(\pi)$ of infinite
jets for a bundle~$\pi\colon E\to M$ whose sections play the r\^{o}le of
unknown functions (fields). This attitude allowed to apply to PDEs powerful
techniques of differential geometry and homological algebra. The latter, in
particular, made it possible to give an invariant and efficient formulation of
higher-order Lagrangian formalism with constraints and calculus of variations
(see~\cite{Vinogradov:SSLFCLLTNT,Vinogradov:AlGFLFT,Vinogradov:CAnPDEqSC,Vinogradov:SSAsNDEqAlGFLFTC}),
see
also~\cite{Takens:GVInPCV,Anderson:InVB,Tsujishita:VBAsDEq,Tsujishita:HMCInSDEq,Tulczyjew:LC,Kupershmidt:GJBSLHF,Manin:AlAsNDEq}
and, as it became clear later, was a bridge to BRST cohomology in gauge
theories, anti-field formalism and related topics,~\cite{Henneaux:CInGFCAp};
see
also~\cite{BarnichBrandtHenneaux:LBCAnFIGT,HenneauxTeitelboim:QGS,GomisParisSamuel:AnAnGTQ}.

Geometrical treatment of differential equations has a long history and
originates in the works by Sophus Lie~\cite{Lie:VDT,Lie:GB,Lie:GAb}, as well
as in research by A.V.~B\"{a}cklund~\cite{Baecklund:ZTF},
G.~Monge~\cite{Monge:ApG}, G.~Darboux~\cite{Darboux:LTGSApGCIn},
L.~Bianchi~\cite{Bianchi:LGD} and, later, by \'{E}lie
Cartan~\cite{Cartan:SDExApG}. Note incidentally that Cartan's theory of
involutivity for external differential systems was an inspiration for another
cohomological theory associated to PDEs and developed in papers by D.~Spencer
and his school,~\cite{Spencer:OvSLPDEq,Goldschmidt:ExTAnLPDEq}. Spencer's work
(the so-called formal theory) closely relates to earlier and unfairly
forgotten results by M.~Janet~\cite{Janet:LSDP} and
Ch.~Riquier~\cite{Riquier:LSDP}; see also~\cite{Pommaret:SPDEqLP} as well
as~\cite{Rashevskii:GTPDEq,Seiler:FTDEqApCAl,BryantChernGardnerGoldschmidtGriffiths:ExDS,KruglikovLychagin:GDEq}.

A milestone in the geometry of differential equations was introduction by
Charles Ehresmann the notion of jet
bundles~\cite{Ehresmann:InTSInPL1,Ehresmann:InTSInPL2} that became the most
adequate language for Lie's theory, but a real revival of the latter came with
the works by L.V.~Ovsyannikov (see his book~\cite{Ovsiannikov:GAnDEq} on group
analysis of PDEs; see
also~\cite{Ibragimov:TGApMP,Stephani:DEqTSUsS,BlumanKumei:SDEq,BlumanCheviakovAnco:ApSMPDEq}).

A new impulse for the reappraisal of the Sophus Lie heritage was given by the
discovery of integrability phenomena in nonlinear
systems~\cite{GardnerGreeneKruskalMiura:MSKVEq,Miura:KVEqGIRExNT,MiuraGardnerKruskal:KVEqGIIExCLCM,SuGardner:KVEqGIIIDKVEqBEq,Gardner:KVEqGIVKVEqHS,KruskalMiuraGardnerZabusky:KVEqGVUnNPCL,GardnerGreeneKruskalMiura:KEqGVMExS,BulloughCaudrey:S,NovikovManakovPitaevskiiZakharov:TS}
in the fall of the 1960s\footnote{Though discussions on ``what is
  integrability'' continued~\cite{Zakharov:WIsIn} and are still held now (see,
  e.g. quite recent papers~\cite{FerapontovKhusnutdinovaKlein:LDInQSHD}
  or~\cite{OdesskiiSokolov:InPRGHF}).} and Hamiltonian interpretation of this
integrability~\cite{ZakharovFaddeev:KdVEqisFInHS,Gardner:KVEqGIVKVEqHS}. In
particular, it became clear that integrable equations possess infinite series
of higher, or generalised, symmetries
(see~\cite{Olver:ApLGDEq,KrasilshchikLychaginVinogradov:GJSNPDEq}), and
classification of evolution equations with respect to this property allowed to
discover new, at that time, integrable
equations,~\cite{MikhailovShabatYamilov:SApCNEqCLInS,SvinolupovSokolovYamilov:BTInEvEq,Mikhailov:In}. Later
the notion of a higher symmetry was generalised further to that of a nonlocal
one~\cite{VinogradovKrasilshchik:MCHSNEvEqNS} and the search for a geometrical
background of nonlocality led to the concept of a differential
covering~\cite{KrasilshchikVinogradov:NTGDEqSCLBT}. The latter proved to play
an important r\^{o}le in the geometry of PDEs and we discuss it in our review.

It also became clear that the majority of integrable evolutionary systems
possess a bi-Hamiltonian structure~\cite{Magri:SInHP,Magri:SMInHEq}, i.e., can
be represented as Hamiltonian flows on the space of infinite jets in at least
two different ways and the corresponding Hamiltonian structures are
compatible. The bi-hamiltonian property, by Magri's
scheme~\cite{Magri:SMInHEq}, leads to the existence of infinite series of
commuting symmetries and conservation laws. In addition, it gives rise to a
recursion operator for higher symmetries that is an efficient tool for
practical construction of symmetry hierarchies. Nevertheless, recursion
operators exist for equations possessing no Hamiltonian structure at all
(e.g., for the Burgers equation). A self-contained cohomological approach to
recursion operators based on Nijenhuis brackets and related to the theory of
deformations for PDE structures is exposed
in~\cite{KrasilshchikKersten:SROpCSDE}.

The literature on the Hamiltonian theory of PDEs is vast and we confine
ourselves here to the key
references~\cite{FaddeevTakhtadzhyan:HMTS,Dickey:SEqHS,Dorfman:DSInNEvEq,DubrovinKricheverNovikov:InSI},
but one feature is common to all research: theories and techniques are
applicable to evolution equations only. Then a natural question arises: what
to do if the equation at hand is not represented in the evolutionary form? We
believe that (at least, a partial) answer to this question can be found in
this paper.

Of course, one of possible solutions is to transform the equation to the
evolutionary form. But:
\begin{itemize}
\item Not all equations can be rendered to this form\footnote{For example,
    gauge invariant equations, such the Yang-Mills, Maxwell, Einstein
    equation, etc., can not be presented in the evolutionary form.}.
\item How to check independence of Hamiltonian (and other) structures on a
  particular representation of our equation? In other words, if we found a
  Hamiltonian operator in one representation what guarantees that it survives
  when the representation is changed?
\item Even if the answer to the previous question is positive, how to
  transform the results when passing to the initial form of the equation?
\end{itemize}

In what follows, we treat any concrete equation ``as is'' and try to uncover
those objects and constructions that are naturally associated to this
equation. In particular, we do not assume existence of any additional
structures that enrich the equation. Such structures are by all means
extremely interesting and lead to very nontrivial classes of equations (e.g.,
equations of hydrodynamical type~\cite{DubrovinNovikov:PBHT}, Monge-Amp\'{e}re
equations~\cite{KushnerLychaginRubtsov:CGNLDEq} or equations associated to Lie
groups~\cite{DrinfeldSokolov:LAlEqKVT}), but here we look for internal
properties of an arbitrary PDE.

As it was said in the very beginning, an equation (or, to be more precise, its
infinite prolongation, i.e., equation itself together with all its
differential consequences) is a submanifold in a jet
space~$\J^\infty(\pi)$. To escape technical difficulties, we consider the
simplest case, when~$\pi\colon E\to M$ is a locally trivial vector bundle,
though all the results remain valid in a more complicated situation (e.g., for
jets of submanifolds). The reader who is interested in local results only may
keep in mind the trivial bundle~$\R^m\times\R^n\to\R^n$ instead of~$\pi$.

Understood in such a way, any equation~$\mathcal{E}\subset\J^\infty(\pi)$ is
naturally endowed with a $(\dim M)$-dimensional integrable distribution~$\C$
(the Cartan distribution) which consists, informally speaking, of planes
tangent to formal solutions of~$\mathcal{E}$. This is the main and essentially
the only geometric structure that we use.

In the research of the PDE differential geometry we use somewhat informal but
quite productive guidelines which were originally introduced
in~\cite{Vinogradov:CNDEq} (see also~\cite{Vinogradov:InInGJS,Vinogradov:GNDEq}) and may be
formulated as
\begin{description}
\item[The structural principle] Any construction and concept must take into
  account the Cartan distribution on~$\mathcal{E}$.
\item[The correspondence principle] ``Physical dimension'' of~$\mathcal{E}$
  is~$n=\dim\C=\dim M$ and differential geometry of~$\mathcal{E}$ reduces to
  the finite-dimensional one when passing to the ``classical limit~$n\to 0$''.
\item[The invariance principle] All constructions must be independent of the
  embedding~$\mathcal{E}\to\J^\infty(\pi)$ and defined by the
  equation~$\mathcal{E}$ itself.
\end{description}
Below we accompany our exposition by toy dictionaries that illustrate the
correspondence between two languages, those of the geometry of PDEs and
classical differential geometry.

The paper consists of three sections. In
Section~\ref{sec:gjsis_jets:jet-spaces} we describe the geometry of the
``empty equation'', i.e., of the jet space~$\J^\infty(\pi)$. In particular, we
define the tangent and cotangent bundles to~$\J^\infty(\pi)$, introduce
variational differential forms and multivectors and define the variational
Schouten bracket. We discuss geometry of Hamiltonian flows on the space of
infinite jets (i.e., Hamiltonian evolutionary equations) and Lagrangian
formalism without constraints. Section~\ref{sec:gjsis_eqs:diff-equat} deals
with the same matters, but in the context of a differential
equation~$\mathcal{E}\subset\J^\infty(\pi)$. Although the exposition in this
part is quite general, the result on the Hamiltonian theory (the definition of
the cotangent bundle, in particular) are valid for the so-called $2$-line
equations only (we call such equations normal in
Section~\ref{sec:gjsis_eqs:diff-equat}). This notion is related to the
cohomological length of the compatibility complex for the linearization
operator of~$\mathcal{E}$ and manifests itself, for example, in the number of
nontrivial lines in Vinogradov's $\C$-spectral
sequence,~\cite{Vinogradov:SSLFCLLTNT} (see
also~\cite{Vinogradov:SSAsNDEqAlGFLFTC,Vinogradov:CAnPDEqSC,BryantGriffiths:CCDSGT}
and~\cite{KrasilshchikVerbovetsky:HMEqMP,Verbovetsky:NHC}). From this point of
view, jet spaces are $1$-line equations and this is the reason why they have
to be treated separately. Finally, in
Section~\ref{sec:gjsis_nonloc:nonlocal-theory} we briefly overview the theory
of differential coverings for PDEs and some of its applications: nonlocal
symmetries and shadows, B\"{a}cklund transformations, etc.

Our exposition of general facts related to the geometry of jet spaces and
infinitely prolonged equations, including the nonlocal theory, is essentially
based on the
books~\cite{KrasilshchikLychaginVinogradov:GJSNPDEq,KrasilshchikVinogradov:SCLDEqMP,KrasilshchikVerbovetsky:HMEqMP}. Lagrangian
formalism, both in the free case (on jets) and with constraints (on
equations), is exposed using the material of~\cite{Vinogradov:CAnPDEqSC} (see
also~\cite{KrasilshchikVinogradov:SCLDEqMP,KrasilshchikVerbovetsky:HMEqMP}). The
geometrical approach to Hamiltonian formalism (including the theory of the
Schouten bracket) is based
on~\cite{KerstenKrasilshchikVerbovetsky:HOpC,KerstenKrasilshchikVerbovetskyVitolo:HSGP}.

A practical implementation of the general theory, in the majority of cases,
needs the use of an appropriate computer algebra software. To avoid technical
details that obscure the essentials, we chose for a ``tutorial example'' the
well known Korteweg-de Vries equation for which all computations are
transparent and can be done by hand. We did our best to illustrate the theory
by a reasonable number of less trivial examples and really do hope that the
result will be interesting to the readers. 

\section{Jet spaces}
\label{sec:gjsis_jets:jet-spaces}

Jet spaces constitute a natural geometric environment for differential
equations and for equations of mathematical physics, in particular. But these
spaces are themselves an interesting geometric object that contains
information on Lagrangian and Hamiltonian formalisms without
constraints. Thus, we begin our exposition with a description of these spaces
and structures related to them.

\subsection{Definition of jet spaces}
\label{sec:gjsis_jets:defin-jet-spac}

Let $\pi\colon E\to M$ be a locally trivial smooth vector bundle\footnote{For
  a definition of jets in a more general setting see,
  e.g.,~\cite{KrasilshchikLychaginVinogradov:GJSNPDEq,Vinogradov:CAnPDEqSC,Saunders:GJB,Saunders:JMNB}.}
over a smooth manifold~$M$, $\dim M=n$, $\dim E=m+n$. In what follows, $M$
will be the manifold of independent variables while sections of~$\pi$ will
play the r\^{o}le of unknown functions (fields). The set of all
sections~$s\colon M\to E$ will be denoted by~$\Se(\pi)$ and it forms a module
over the algebra~$C^\infty(M)$. Two sections~$s$, $s'\in\Se(\pi)$ are said to
be \emph{$k$-equivalent} at a point~$x\in M$ if their graphs are tangent to
each other with order~$k$ at the point~$s(x)=s'(x)\in E$. The equivalence
class of~$s$ with respect to this relation is denoted by~$[s]_x^k$ and is
called the \emph{$k$-jet} of~$s$ at~$x$. The set
\begin{equation*}
  \J^k(\pi)=\sd{[s]_x^k}{x\in M,s\in\Se(\pi)}
\end{equation*}
is endowed with a natural structure of a smooth manifold; the latter is called
the \emph{manifold of $k$-jets} of sections of~$\pi$. Moreover, the maps
\begin{equation}\label{eq:gjsis_jets:1}
  \pi_k\colon\J^k(\pi)\to M,\quad [s]_x^k\mapsto x,
\end{equation}
and
\begin{equation}\label{eq:gjsis_jets:2}
  \pi_{k,l}\colon\J^k(\pi)\to\J^l(\pi),\quad [s]_x^k\mapsto[s]_x^l,\qquad
  k\ge l,
\end{equation}
are smooth fibre bundles, $\pi_k$ being vector bundles. For any
section~$s\in\Se(\pi)$ the map
\begin{equation}\label{eq:gjsis_jets:3}
  j_k(s)\colon M\to\J^k(\pi),\quad x\mapsto[s]_x^k,
\end{equation}
is a smooth section of~$\pi_k$ that is called the \emph{$k$-jet} of~$s$.

Here we are mostly interested in the case~$k=\infty$, i.e., in the
space~$\J^\infty(\pi)$. It can be understood as the inverse limit of the chain
\begin{equation}\label{eq:gjsis_jets:4}
    \qquad\cdots\xrightarrow{}\J^{k+1}(\pi)\xrightarrow{\pi_{k+1,k}}\J^k(\pi)\xrightarrow{}\cdots\xrightarrow{}
    \J^1(\pi)\xrightarrow{\pi_{1,0}}\J^0(\pi)=E\xrightarrow{\pi}M.
\end{equation}
Due to projections~\eqref{eq:gjsis_jets:4} there exist monomorphisms of
function algebras
\begin{equation}\label{eq:gjsis_jets:5}
  C^\infty(M)\subset\F_0(\pi)\subset\dots\subset\F_k(\pi)\subset\F_{k+1}(\pi)
  \subset\dots,
\end{equation}
where~$\F_k(\pi)=C^\infty(\J^k(\pi))$, and we define the \emph{algebra of
  smooth functions} on~$\J^\infty(\pi)$ as the filtered
algebra~$\F(\pi)=\cup_k\F_k(\pi)$. Elements of~$\F(\pi)$ are identified with
nonlinear scalar differential operators acting on sections of~$\pi$ by the
following rule:
\begin{equation}\label{eq:gjsis_jets:6}
  \Delta_f(s)=j_\infty(s)^*(f),\qquad s\in\Se(\pi),f\in\F(\pi).
\end{equation}

More general, let~$\pi'\colon E'\to M$ be another vector bundle
and~$\pi^*(\pi')$ be its pull-back to~$\J^\infty(\pi)$. Introduce the
notation~$\F(\pi,\pi')=\Se(\pi^*(\pi'))$. Then any section~$f\in\F(\pi,\pi')$
is identified, by a formula similar to~\eqref{eq:gjsis_jets:6}, with a
nonlinear differential operator that acts from~$\Se(\pi)$ to~$\Gamma(\pi')$.

\subsection{Vector fields and differential forms}
\label{sec:gjsis_jets:vect-fields-diff}

A \emph{vector field} on~$\J^\infty(\pi)$ is a derivation of the function
algebra~$\F(\pi)$, i.e., an $\R$-linear map~$X\colon\F(\pi)\to\F(\pi)$ such
that
\begin{equation*}
  X(fg)=fX(g)+gX(f)
\end{equation*}
for all $f$, $g\in\F(\pi)$. The set of all vector fields is denoted
by~$\X(\pi)$ and it is a Lie algebra with respect to the commutator (the
\emph{Lie bracket}).

The definition of a \emph{differential form} of degree~$r$ on~$\J^\infty(\pi)$
is similar to that of smooth functions. Using
projections~\eqref{eq:gjsis_jets:4} we consider the
embeddings~$\La^r(\J^k(\pi))\subset\La^r(\J^{k+1}(\pi))$ and
set~$\La^r(\pi)=\cup_k\La^r(\J^k(\pi))$. We shall also consider the Grassmann
algebra of all forms~$\La^*(\pi)=\oplus_{r\ge 0}\La^r(\pi)$ with respect to
the wedge product.

\begin{coordinates}
  \label{sec:gjsis_jets:vect-fields-diff-1}
  Let~$\mathcal{U}\subset M$ be a coordinate neighbourhood such that the
  bundle~$\pi$ becomes trivial over~$\mathcal{U}$. Choose local
  coordinates~$x^1,\dots,x^n$ in~$\mathcal{U}$ and~$u^1,\dots,u^m$ along the
  fibres of~$\pi$ over~$\mathcal{U}$. Then the \emph{adapted} coordinates
  in~$\pi^{-1}(\mathcal{U})\subset\J^\infty(\pi)$ naturally arise. These
  coordinates are denoted by~$u_I^j$, $I$ being a multi-index, and are defined
  by
  \begin{equation*}
    j_\infty(s)^*(u_I^j)=\frac{\partial^{\abs{I}}s^j}{\partial x^I},
  \end{equation*}
  where~$s=(s^1,\dots,s^m)$ is a local section of~$\pi$ over~$\mathcal{U}$. In
  other words, the coordinate functions~$u_I^j$ correspond to partial
  derivatives of local sections.

  In these coordinates, smooth function on~$\J^\infty(M)$ are of the form
  \begin{equation*}
    f=f(x^i,u_I^j),
  \end{equation*}
  where the number of arguments is \emph{finite}. Vector fields are
  represented as \emph{infinite} sums
  \begin{equation*}
    X=\sum_ia_i\frac{\partial}{\partial x^i}+
    \sum_{I,j}a_I^j\frac{\partial}{\partial u_I^j},\qquad
    a_i,a_I^j\in\F(\pi),
  \end{equation*}
  while differential forms of degree~$r$ are \emph{finite} sums
  \begin{equation*}
    \omega=\sum b_{i_1,\dots,i_c,j_{c+1},\dots,j_r}^{I_{c+1},\dots,I_r}
    \rmd x^{i_1}\wedge\dots\wedge\rmd x^{i_c}\wedge
    \rmd u_{I_{c+1}}^{j_{c+1}}\wedge\dots\wedge\rmd u_{I_r}^{j_r}.
  \end{equation*}
\end{coordinates}

\subsection{Main structure: the Cartan distribution}
\label{sec:gjsis_jets:main-struct-cart}

Let~$\theta\in\J^\infty(\pi)$. Then the graphs of all sections~$j_\infty(s)$,
$s\in\Se(\pi)$, passing through the point~$\theta$ have a common
$n$-dimensional tangent plane~$\C_\theta$ (the \emph{Cartan plane}). The
correspondence~$\C\colon\theta\mapsto\C_\theta$ is an
integrable\footnote{\label{fn:1}Integrability is understood formally here and
  means that if~$X$ and~$Y$ are two vector fields lying in~$\C$ then their
  bracket~$[X,Y]$ lies in~$\C$ as well. Since~$\J^\infty(\pi)$ is
  infinite-dimensional, this does not mean that the Frobenius theorem holds
  for~$\C$: for any point~$\theta\in\J^\infty(\pi)$ there exist infinitely
  many maximal integral manifolds that contain~$\theta$. On the other hand,
  if~$\C$ is restricted to an equation (see below
  Section~\ref{sec:gjsis_eqs:diff-equat}), there may exist no maximal integral
  manifold at all.}  $n$-dimensional distribution on~$\J^\infty(\pi)$ that is
called the \emph{Cartan distribution}. This distribution is the basic
geometric structure on the manifold~$\J^\infty(\pi)$. In particular, the
following result is valid:
\begin{proposition}
  \label{prop:gjsis_jets:1}
  A submanifold in~$\J^\infty(M)$ is a maximal integral manifold of~$\C$ if
  and only if it is the graph of~$j_\infty(s)$, where~$s$ is a local section
  of~$\pi$.
\end{proposition}

Moreover, since the planes~$\C_\theta$ project to~$M$ non-degenerately, any
vector field~$X$ on~$M$ can be uniquely lifted up to a field~$\C X$
on~$\J^\infty(\pi)$. In such a way, one obtains a connection in the
bundle~$\pi_\infty$ called the \emph{Cartan connection}. This connection is
flat, i.e.,
\begin{equation}\label{eq:gjsis_jets:7}
  \C[X,Y]=[\C X,\C Y]
\end{equation}
for all vector fields~$X$, $Y$ on~$M$. Due to \eqref{eq:gjsis_jets:7}, the
space~$\C\X(\pi)$ of all vector fields lying in the Cartan distribution is a
Lie subalgebra in~$\X(\pi)$. Vector fields belonging to~$\X(\pi)$ will be
called \emph{Cartan fields}.

Any vector field~$Z\in\X(\pi)$ can be uniquely decomposed to its
\emph{vertical} and \emph{horizontal} components,
\begin{equation}\label{eq:gjsis_jets:8}
  Z=Z^v+Z^h,
\end{equation}
where~$Z^v$ is the projection of~$X$ to the fibre of the bundle~$\pi_\infty$
along Cartan planes, while~$Z^h$ lies in the Cartan distribution. Thus, one
has
\begin{equation}\label{eq:gjsis_jets:9}
  \X(\pi)=\X^v(\pi)\oplus\C\X(\pi),
\end{equation}
where~$\X^v(\pi)$ is the Lie algebra of vertical vector fields.

Dually to~\eqref{eq:gjsis_jets:9}, the module of differential
forms~$\La^1(\pi)$ splits into the direct sum
\begin{equation}\label{eq:gjsis_jets:10}
  \La^1(\pi)=\Lac^1(\pi)\oplus\Lah^1(\pi),
\end{equation}
where~$\Lac^1(\pi)$ consists of $1$-forms that annihilate the Cartan
distribution (they will be called \emph{Cartan forms}, or \emph{higher contact
  forms}), while elements of~$\Lah^1(\pi)$ are \emph{horizontal forms}.

\begin{coordinates}
  \label{sec:gjsis_jets:main-struct-cart-1}
  Choose an adapted coordinate system~$(x^i,u_I^j)$ in~$\J^\infty(\pi)$. Then
  one has
  \begin{equation}\label{eq:gjsis_jets:11}
    \C\colon\frac{\partial}{\partial x^i}\mapsto D_i=
    \frac{\partial}{\partial x^i}+
    \sum_{I,j}u_{Ii}^j\frac{\partial}{\partial u_I^j}.
  \end{equation}
  The fields~$D_i$ are called \emph{total derivatives} and they span the
  Cartan distribution.

  For a basis in the module~$\Lac^1(\pi)$ one can choose the forms
  \begin{equation}\label{eq:gjsis_jets:12}
    \omega_I^j=\rmd u_I^j-\sum_i u_{Ii}^j\rmd x^i,
  \end{equation}
  while horizontal forms are
  \begin{equation}\label{eq:gjsis_jets:13}
    \omega=\sum_i a_i\rmd x^i,\qquad a_i\in\F(\pi).
  \end{equation}
\end{coordinates}

\begin{remark}
  \label{rem:gjsis_jets:6}
  It should be noted that all results and constructions below are valid not
  for the entire jet space only, but for an arbitrary open domain
  in~$\J^\infty(\pi)$. Everywhere below, when speaking about~$\J^\infty(\pi)$,
  we actually mean an open domain.
\end{remark}

\subsection{Evolutionary vector fields and linearizations}
\label{sec:gjsis_jets:evol-vect-fields}

We shall now describe infinitesimal symmetries of the Cartan distribution
on~$\J^\infty(\pi)$. A vector field~$X\in\X(\pi)$ is a \emph{symmetry}
if~$[X,Z]\in\C\X(\pi)$ as soon as~$Z\in\C\X(\pi)$. The space~$\C\X(\pi)$ is an
ideal in the Lie algebra~$\Xc(\pi)$ of symmetries. Due to integrability of the
Cartan distribution, any~$Z\in\C\X(\pi)$ is a symmetry, and we call such
symmetries \emph{trivial}. Thus, we introduce the Lie algebra of
\emph{nontrivial symmetries} as
\begin{equation*}
  \sym\pi=\Xc(\pi)/\C\X(\pi).
\end{equation*}
By \eqref{eq:gjsis_jets:9},~$\sym\pi$ is identified with the vertical part
of~$\Xc(\pi)$.

Take a vector field~$X\in\sym\pi$ and restrict it to the
subalgebra~$\F_0(\pi)\subset\F(\pi)$. Then this restriction can be identified
with an element~$\phi_X\in\F(\pi,\pi)$. For shortness, we shall use the
notation~$\F(\pi,\pi)=\kappa(\pi)$.
\begin{theorem}
  \label{thm:gjsis_jets:1}
  The correspondence~$X\mapsto\varphi_X$ defines a bijection between~$\sym\pi$
  and~$\kappa(\pi)$.
\end{theorem}
The element~$\phi_X$ is called the \emph{generating section} of a
symmetry~$X$, while the symmetry corresponding to a
section~$\phi\in\kappa(\pi)$ is called an \emph{evolutionary vector field} and
is denoted by~$\Ev_\phi$.

Theorem~\ref{thm:gjsis_jets:1} allows to introduce an $\F(\pi)$-module
structure into the Lie algebra~$\sym\pi$ by setting
\begin{equation*}
  f\cdot\Ev_\phi=\Ev_{f\phi}.
\end{equation*}
This multiplication differs from the usual multiplication of a vector field by
functions and does not survive when passing from the space of jets to
equations (see Section~\ref{sec:gjsis_eqs:diff-equat} below).

On the other hand, the same theorem defines a Lie algebra structure
in~$\kappa(\pi)$: the \emph{Jacobi bracket}~$\{\phi,\psi\}$ is uniquely given
by the equality
\begin{equation}\label{eq:gjsis_jets:14}
  \Ev_{\{\phi,\psi\}}=[\Ev_\phi,\Ev_\psi].
\end{equation}
The Jacobi bracket can also be computed using the formula
\begin{equation}\label{eq:gjsis_jets:15}
  \{\phi,\psi\}=\Ev_\phi(\psi)-\Ev_\psi(\phi).
\end{equation}

\begin{remark}\label{sec:gjsis_jets:Lie-Back}
  We pointed out in Footnote~\ref{fn:1} that integrability of distributions
  on~$\J^\infty(\pi)$ differs from the one on finite-dimensional
  manifolds. The same holds for integrability (i.e., existence of the
  corresponding one-parameter group of transformations) of vector
  fields. Generally speaking, a vector field on~$\J^\infty(\pi)$ is not
  integrable in this sense, but there exists an important class of vector
  fields that are integrable. Namely, if~$X$ is a field on~$\J^k(\pi)$,
  $k<\infty$, that preserves the Cartan distribution then it can be
  \emph{lifted} in a natural way to~$\J^{k+1}(\pi)$,
  see~\cite{KrasilshchikVinogradov:SCLDEqMP}. The entire collection of such
  fields determines a vector field on~$\J^\infty(\pi)$, and any such a field
  possesses a one-parameter group of transformations. The infinitesimal
  version of the Lie-B\"acklund theorem states that all such fields are the
  lifts of arbitrary fields on~$\J^0(\pi)$ (when~$\dim\pi>1$) or of a contact
  vector field on~$\J^1(\pi)$ (when~$\dim\pi=1$). A complete description of
  integrable vector fields on~$\J^\infty(\pi)$ can be found
  in~\cite{Chetverikov:SInFInPEq,Chetverikov:SInF}.
\end{remark}

\begin{coordinates}
  \label{sec:gjsis_jets:evol-vect-fields-1}
  Let, in an adapted coordinate system, a section~$\phi\in\kappa(\pi)$ be of
  the form~$\phi=(\phi^1,\dots,\phi^m)$. Then the corresponding evolutionary
  vector field is
  \begin{equation}\label{eq:gjsis_jets:16}
    \Ev_\phi=\sum_{I,j}D_I(\phi^j)\frac{\partial}{\partial u_I^j},
  \end{equation}
  where~$D_I=D_{i_1}\circ\dots\circ D_{i_l}$ is the composition of total
  derivatives corresponding to the multi-index~$I=i_1\dots i_l$.
  
  If~$\psi=(\psi^1,\dots,\psi^m)$ is another element of~$\kappa(\pi)$ then the
  components of the Jacobi bracket are
  \begin{equation}\label{eq:gjsis_jets:17}
    \{\phi,\psi\}^j=
    \sum_{I,\alpha}\left(D_I(\phi^\alpha)\frac{\partial\psi^j}{\partial u_I^\alpha}
      -D_I(\psi^\alpha)\frac{\partial\phi^j}{\partial u_I^\alpha}\right),\qquad
    j=1,\dots,m.
  \end{equation}
\end{coordinates}

Fix a section~$\psi\in\kappa(\pi)$ and consider the map
\begin{equation}\label{eq:gjsis_jets:18}
  \ell_\psi\colon\kappa(\pi)\to\kappa(\pi),\qquad
  \ell_\psi(\phi)=\Ev_\phi(\psi).
\end{equation}
This map is called the \emph{linearization} of the element~$\psi$ (recall
that~$\psi$ may be identified with a nonlinear differential operator acting
from~$\pi$ to~$\pi$).

More generally, let~$\pi'\colon E'\to M$ be a vector bundle. Then the action
of an evolutionary vector field~$\Ev_\phi$ can be extended to
\begin{equation*}
  \Ev_\phi\colon\F(\pi,\pi')\to\F(\pi,\pi')
\end{equation*}
in a well defined way. Consider a section~$\psi\in\F(\pi,\pi')$, i.e., a
nonlinear differential operator from~$\Se(\pi)$ to~$\Se(\pi')$. Its
linearization is the map
\begin{equation}\label{eq:gjsis_jets:19}
  \ell_\psi\colon\kappa(\pi)\to\F(\pi,\pi')
\end{equation}
is defined similar to~\eqref{eq:gjsis_jets:18}.

\begin{coordinates}\label{sec:gjsis_jets:evol-vect-fields-2}
  Let, in adapted coordinates,~$\phi=(\phi^1,\dots,\phi^m)$
  and~$\psi=(\psi^1,\dots,\psi^{m'})$. Then the $j$th component
  of~$\ell_\psi(\phi)$ is
  \begin{equation*}
    \sum_{I,\alpha}\frac{\partial\psi^j}{\partial u_I^\alpha}D_I(\phi^\alpha),
  \end{equation*}
  i.e., the linearization is a matrix operator of the form
  \begin{equation}\label{eq:gjsis_jets:20}
    \ell_\psi=
    \left\Vert\sum_I\frac{\partial\psi^j}{\partial u_I^\alpha}D_I
    \right\Vert_{\alpha=1,\dots,m.}^{j=1,\dots,m'}
  \end{equation}

  Using linearizations, formula~\eqref{eq:gjsis_jets:15} can be rewritten as
  \begin{equation}\label{eq:gjsis_jets:21}
    \{\phi,\psi\}=\ell_\psi(\phi)-\ell_\phi(\psi).
  \end{equation}
\end{coordinates}

\subsection{$\C$-differential operators}
\label{sec:gjsis_jets:c-diff-oper}

From \eqref{eq:gjsis_jets:20} we see that linearizations are differential
operators in total derivatives. We call such operators \emph{$\C$-differential
  operators}. More precisely, let~$\xi$ and~$\xi'$ be two vector bundles
over~$\J^\infty(\pi)$ and~$P$, $P'$ be the $\F(\pi)$-modules of their
sections. An $\R$-linear map~$\Delta\colon P\to P'$ is a $\C$-differential
operator of order~$k$ if for any point~$\theta\in\J^\infty(\pi)$ and a
section~$p\in P$ the value of~$\Delta(p)$ at~$\theta$ is completely determined
by the values of~$D_I(p)$, $\abs{I}\le k$, at this point. The space of all
such operators is denoted by~$\CDiff_k(P,P')$ and we also
set~$\CDiff(P,P')=\cup_k\CDiff_k(P,P')$.

A closely related notion to that of a $\C$-differential operator is
\emph{horizontal jets}. Let~$P$ be as above. We say that two sections~$p$,
$p'\in P$ are \emph{horizontally $k$-equivalent} (the case~$k=\infty$ is
included) at a point~$\theta\in\J^\infty(\pi)$ if~$D_I(p)=D_I(p')$ at~$\theta$
for all~$I$ such that~$\abs{I}\le k$. Denote the equivalence class
by~$\{p\}_\theta^k$. The set
\begin{equation*}
  \JJ^k(P)=\sd{\{p\}_\theta^k}{\theta\in\J^\infty(\pi),p\in P}
\end{equation*}
forms a smooth manifold which is fibred over~$\J^\infty(\pi)$,
\begin{equation*}
  \xi_k\colon\JJ^k(P)\to\J^\infty(\pi)\quad\{\rho\}_\theta^k\mapsto\theta.
\end{equation*}
The section~$j_k^h(p)\colon\J^\infty(\pi)\to\JJ^k(P)$,
$\theta\mapsto\{p\}_\theta^k$, is called the \emph{horizontal jet} of~$p\in
P$.
\begin{proposition}
  \label{prop:gjsis_jets:2}
  For any $\C$-differential operator~$\Delta\in\CDiff_k(P,P')$ there exists a
  unique morphism~$\Phi_\Delta$ of vector bundles~$\xi_k$ and~$\xi'$ such
  that~$\Delta(p)=\Phi_\Delta(j_k^h(p))$ for any~$p\in P$.
\end{proposition}

Two natural identifications will be useful below.
\begin{proposition}
  \label{prop:gjsis_jets:3}
  For any vector bundle~$\pi$ one has:
  \begin{enumerate}
  \item The module~$\Lac^1(\pi)$ is isomorphic
    to~$\CDiff(\kappa(\pi),\F(\pi))$.
  \item The module~$\JJ^\infty(\kappa(\pi))$ is isomorphic to~$\X^v(\pi)$. The
    vector fields corresponding to sections of the form~$j_\infty^h(p)$ are
    evolutionary fields.
  \end{enumerate}
\end{proposition}

\begin{coordinates}
  \label{sec:gjsis_jets:c-diff-oper-1}
  Choose an adapted coordinate system in the manifold~$\J^\infty(\pi)$ and
  let~$r$, $r'$ be dimensions of the bundles~$\xi$, $\xi'$, respectively. Then
  any operator~$\Delta\in\CDiff(P,P')$ is of the form
  \begin{equation}\label{eq:gjsis_jets:22}
    \Delta=
    \left\Vert\sum_Ia_{\alpha\beta}^ID_I\right\Vert_{\alpha=1,\dots,r',}^{\beta=1,\dots,r}
    \qquad a_{\alpha\beta}^I\in\F(\pi).
  \end{equation}

  If~$v^1,\dots,v^r$ are fibre-wise coordinates in the bundle~$\xi$ then the
  adapted coordinates~$v_K^l$ in~$\JJ^\infty(\xi)$, $K$ being a multi-index,
  $l=1,\dots,r$, are determined by the equalities
  \begin{equation}\label{eq:gjsis_jets:23}
    j_\infty^h(s)^*(v_K^l)=D_K(s^l),
  \end{equation}
  where~$s=(s^1,\dots,s^r)$ is a local section of the bundle~$\xi$.
\end{coordinates}

\begin{remark}
  \label{rem:gjsis_jets:1}
  The space of horizontal jets~$\JJ^\infty(P)$ is also endowed with an
  integrable distribution similar to the Cartan one:
  if~$\theta\in\JJ^\infty(P)$ then the corresponding plane~$\C_\theta$ is
  tangent to the graphs of horizontal jets passing through this point. The
  differential of the map~$\xi_\infty\colon\JJ^\infty(P)\to\J^\infty(\pi)$
  isomorphically projects~$\C_\theta$ to~$\C_{\xi_\infty(\theta)}$.

  Moreover, if~$P$ is of the form~$P=\Se(\pi_\infty^*(\xi))$, where~$\xi$ is a
  vector bundle over~$M$ and~$\pi_\infty^*(\xi)$ is its pull-back, then one
  has a diffeomorphism
  \begin{equation*}
    \JJ^\infty(P)=\J^\infty(\pi\times_M\xi),
  \end{equation*}
  where~$\pi\times_M\xi$ is the Whitney product, and this isomorphism takes
  the Cartan distribution on~$\JJ^\infty(P)$ to the one
  on~$\J^\infty(\pi\times_M\xi)$.
\end{remark}

\begin{remark}
  \label{rem:gjsis_jets:2}
  In the case when the modules~$P=\Se(\pi_\infty^*(\xi))$
  and~$P'=\Se(\pi_\infty^*(\xi'))$ are of the form considered in
  Remark~\ref{rem:gjsis_jets:1} $\C$-differential operators from~$P$ to~$P'$
  may be understood as non-linear differential operators that take sections
  of~$\pi$ to linear differential operators acting from~$\xi$ to~$\xi'$.
\end{remark}

\subsection{The variational complex and Lagrangian formalism}
\label{sec:gjsis_jets:vari-compl-lagr}

Consider now decomposition~\eqref{eq:gjsis_jets:10} and note that it implies a
more general splitting
\begin{equation}\label{eq:gjsis_jets:24}
  \La^k(\pi)=\bigoplus_{p+q=k}\Lac^p(\pi)\otimes\Lah^q(\pi),
\end{equation}
where
\begin{equation*}
  \Lac^p(\pi)=
  \underbrace{\Lac^1(\pi)\wedge\dots\wedge\Lac^1(\pi)}_{p\text{
      times}},\qquad
  \Lah^q(\pi)=
  \underbrace{\Lah^1(\pi)\wedge\dots\wedge\Lah^1(\pi)}_{q\text{
      times}}.
\end{equation*}
Introduce the notation~$\Lac^p(\pi)\otimes\Lah^q(\pi)=E_0^{p,q}(\pi)$. By
Proposition~\ref{prop:gjsis_jets:3}~(1), this space is identified with the
module~$\CDiffsk_p(\kappa(\pi),\Lah^q(\pi))$ of $p$-linear skew-symmetric
$\C$-differential operators acting from~$\kappa(\pi)$ to~$\Lah^q(\pi)$.

The de~Rham differential~$\rmd\colon\La^k(\pi)\to\La^{k+1}(\pi)$,
by~\eqref{eq:gjsis_jets:24}, splits into two parts~$\rmd=\rmdc+\rmdh$, where
\begin{equation}\label{eq:gjsis_jets:25}
  \rmdc=\rmdc^{p,q}\colon E_0^{p,q}(\pi)\to E_0^{p+1,q}(\pi),\qquad
  \rmdh=\rmdh^{p,q}\colon E_0^{p,q}(\pi)\to E_0^{p,q+1}(\pi)
\end{equation}
are the \emph{vertical} (or \emph{Cartan}) and \emph{horizontal}
differentials, respectively These differentials anti-commute, i.e.,
\begin{equation}\label{eq:gjsis_jets:26}
  \rmdc\circ\rmdh+\rmdh\circ\rmdc=0.
\end{equation}

\begin{coordinates}
  For a function on~$\J^\infty(\pi)$ (i.e., a $0$-form) the action of the
  Cartan differential is given by
  \begin{equation}\label{eq:gjsis_jets:27}
    \rmdc f=\sum_{I,j}\frac{\partial f}{\partial u_I^j}\rmdc u_I^j,
  \end{equation}
  where~$\rmdc u_I^j=\omega_I^j$ are the Cartan forms presented in
  \eqref{eq:gjsis_jets:12}, while the horizontal differential acts as follows
  \begin{equation}\label{eq:gjsis_jets:28}
    \rmdh f=\sum_i D_i(f)\rmd x^i.
  \end{equation}
  To compute the action on arbitrary forms it suffices to
  use~\eqref{eq:gjsis_jets:27} and~\eqref{eq:gjsis_jets:28} and the fact
  that~$\rmdc$ and~$\rmdh$ differentiate the wedge product and anti-commute
  with the de~Rham differential.
\end{coordinates}

Let~$E_1^{p,q}(\pi)$ be the cohomology of~$\rmdh$ at the term~$E_0^{p,q}(\pi)$,
i.e.,
\begin{equation}\label{eq:gjsis_jets:29}
  E_1^{p,q}(\pi)=\ker\rmdh^{p,q}/\im\rmdh^{p,q-1}.
\end{equation}
Due to~\eqref{eq:gjsis_jets:26}, the vertical differential~$\rmdc$ induces the
differentials~$E_1^{p,q}(\pi)\to E_1^{p+1,q}(\pi)$ which will be denoted
by~$\delta^{p,q}$. The groups~$E_1^{p,q}(\pi)$, together with the
differentials~$\delta^{p,q}$, play one of the most important r\^{o}les in the
geometry of jets providing the background for Lagrangian formalism without
constraints. To describe them, we shall need new notions.

Let~$\xi$ be a vector bundle over~$\J^\infty(\pi)$ and~$P=\Se(\xi)$. Introduce
the \emph{adjoint} module $\hat{P}=\Hom_{\F(\pi)}(P,\Lah^n(\pi))$. Consider
another module of sections~$Q$ and a $\C$-differential operator~$\Delta\colon
P\to Q$. Then the \emph{adjoint} operator~$\Delta^*\colon\hat{Q}\to\hat{P}$ is
defined and it enjoys the \emph{Green formula}
\begin{equation}\label{eq:gjsis_jets:30}
  \langle\Delta(p),\hat{q}\rangle-\langle p,\Delta^*(\hat{q})\rangle=\rmdh\omega
\end{equation}
for all~$p\in P$, $\hat{q}\in\hat{Q}$ and
some~$\omega=\omega(p,\hat{q})\in\Lah^{n-1}(\pi)$,
where~$\langle\cdot\,,\cdot\rangle$ is the natural pairing between the module
and its adjoint. It is useful to keep in mind that the
correspondence~$(p,\hat{q})\mapsto\omega(p,\hat{q})$ is a $\C$-differential
operator with respect to both arguments.

An operator~$\Delta$ is called \emph{self-adjoint} if~$\Delta^*=\Delta$ and
\emph{skew-adjoint} if~$\Delta^*=-\Delta$.

\begin{coordinates}
  \label{sec:gjsis_jets:vari-compl-lagr-1}
  If~$\Delta=\sum_I a_ID_i$ is a scalar operator then
  \begin{equation}\label{eq:gjsis_jets:31}
    \Delta^*=\sum_I(-1)^{\abs{I}}D_I\circ a_I.
  \end{equation}
  For a matrix $\C$-differential
  operator~$\Delta=\left\Vert\Delta_{ij}\right\Vert$ one has
  \begin{equation}\label{eq:gjsis_jets:32}
    \Delta^*=\left\Vert\Delta_{ji}^*\right\Vert,
  \end{equation}
  where~$\Delta_{ji}^*$ is given by \eqref{eq:gjsis_jets:31}.
\end{coordinates}

We can now describe the groups~$E_1^{p,q}(\pi)$.
\begin{theorem}[One-line Theorem]
  \label{thm:gjsis_jets:2}
  Let~$\pi\colon E\to M$, $\dim M=n$, be a locally trivial
  vector bundle. Then:
  \begin{enumerate}
  \item the groups~$E_1^{0,q}(\pi)$, $q=0,\dots,n-1$, are isomorphic to the
    de~Rham cohomology groups~$H^q(M)$ of the manifold~$M$;
  \item the group~$E_1^{0,n}(\pi)$ consists of the Lagrangians depending on
    the fields that are sections of~$\pi$;
  \item the groups~$E_1^{p,n}(\pi)$, $p>0$, are identified with the
    modules~$\CDiffskalt_{p-1}(\kappa(\pi),\hat{\kappa}(\pi))$ of
    $(p-1)$-linear skew-symmetric $\C$-differential operators that are
    skew-adjoint in each argument (in particular,
    $E_1^{1,n}(\pi)=\hat{\kappa}(\pi)$);
  \item all other terms are trivial.
  \end{enumerate}
\end{theorem}

\begin{remark}
  The group~$E_1^{0,n}(\pi)$ is also called the $n$th horizontal cohomology
  group of~$\pi$ and denoted by~$H_h^n(\pi)$.
\end{remark}

\begin{remark}
  \label{rem:gjsis_jets:3}
  The construction above is a particular case of Vinogradov's
  \emph{$\C$-spectral
    sequence},~\cite{Vinogradov:AlGFLFT,Vinogradov:SSAsNDEqAlGFLFTC,Vinogradov:SSLFCLLTNT,Vinogradov:CAnPDEqSC,KrasilshchikVerbovetsky:HMEqMP}.
\end{remark}

In what follows, we shall assume the manifold~$M$ to be cohomologically
trivial, i.e., its de~Rham cohomology is isomorphic to~$\R$.

Define the operator~$\delta\colon\Lah^n(\pi)\to\hat{\kappa}(\pi)$ as the
composition of the projection~$\Lah^n(\pi)\to E_1^{0,n}$ and the
differential~$\delta_1^{0,n}\colon E_1^{0,n}\to E_1^{1,n}$.

\begin{remark}
  In what follows, we shall also use the notation~$\delta$ for the
  differential~$\delta_1^{0,n}$ itself.
\end{remark}

\begin{coordinates}
  If~$\omega\in\Lah^n(\pi)$ is of the form~$L\rmd x^1\wedge\dots\wedge\rmd x^n$
  then
  \begin{equation}\label{eq:gjsis_jets:33}
    \delta(\omega)=
    \left(\frac{\delta L}{\delta u^1},\dots,\frac{\delta L}{\delta u^m}\right),
  \end{equation}
  where
  \begin{equation}\label{eq:gjsis_jets:34}
    \frac{\delta L}{\delta u^j}=\sum_I(-1)^{\abs{I}}D_I\frac{\partial
      L}{\partial u_I^j}
  \end{equation}
  are \emph{variational derivatives}. Thus,~$\delta$ is the \emph{Euler}
  operator and it takes a Lagrangian density~$\omega$ to the corresponding
  \emph{Euler operator}.
\end{coordinates}

\begin{proposition}
  \label{prop:gjsis_jets:4}
  Let~$\pi$ be a locally trivial vector bundle over a cohomologically trivial
  manifold~$M$. Then the complex
  \begin{multline}\label{eq:gjsis_jets:35}
    \F(\pi)\xrightarrow{\rmdh^{0,0}}\Lah^1(\pi)\xrightarrow{\rmdh^{0,1}}
      \cdots\xrightarrow{}\Lah^{n-1}(\pi)\xrightarrow{\rmdh^{0,n-1}}\Lah^n(\pi) \\
      \xrightarrow{\delta}\hat{\kappa}\xrightarrow{\delta_1^{1,n}}
      \CDiffskalt_1(\kappa,\hat{\kappa})\xrightarrow{\delta_1^{2,n}}
      \CDiffskalt_2(\kappa,\hat{\kappa})\xrightarrow{}\cdots
  \end{multline}
  is exact, i.e., the kernel of each differential coincides
  with the image of the preceding one.

  In the coordinate-free way, the Euler-Lagrange operator can be
  computed by
  \begin{equation}\label{eq:gjsis_jets:36}
    \delta(\omega)=\ell_\omega^*(1),\qquad\omega\in\Lah^n(\pi),
  \end{equation}
  while the differentials~$\delta_1^{p,n}$ enjoy the equality
  \begin{equation}\label{eq:gjsis_jets:37}
    (\delta_1^{p,n}\Delta)(\phi_1,\dots,\phi_p)=
    \sum_{i=1}^p(-1)^{i+1}
    \ell_{\Delta,\phi_1,\dots,\hat{\phi}_i,\dots,\phi_p}(\phi_i)
    +(-1)^p\ell_{\Delta,\phi_1,\dots,\phi_{p-1}}^*(\phi_p),
  \end{equation}
  where~$\Delta\in\CDiffskalt_{p-1}(\kappa(\pi),\hat{\kappa}(\pi))$
  and~$\phi_1,\dots,\phi_p\in\kappa$. Here and below we use the
  notation~$\ell_{\Delta,\phi_1,\dots,\phi_k}(\phi)=\Ev_\phi(\Delta)(\phi_1,\dots,\phi_k)$
  In particular, if~$\psi\in\hat{\kappa}(\pi)$ then
  \begin{equation}\label{eq:gjsis_jets:38}
    \delta_1^{1,n}(\psi)=\ell_\psi-\ell_\psi^*.
  \end{equation}
\end{proposition}

\begin{remark}
  Complex~\eqref{eq:gjsis_jets:35} is exact starting from the
  term~$\Lah^n(\pi)$ independently of cohomological properties of the
  manifold~$M$. This complex is called the \emph{global variational complex}
  of the bundle~$\pi$, see~\cite{Vitolo:VS}.
\end{remark}

As a consequence of Proposition~\ref{prop:gjsis_jets:4}, we obtain the
following result:
\begin{theorem}
  \label{thm:gjsis_jets:3}
  For a vector bundle~$\pi$ one has:
  \begin{enumerate}
  \item The action functional
    \begin{equation*}
      s\mapsto\int_Mj_\infty(s)^*(\omega),\qquad s\in\Se(\pi),\quad
      \omega\in\Lah^n(\pi),
    \end{equation*}
    is stationary on a section~$s$ if and only
    if~$j_\infty(s)^*(\delta(\omega))=0$ (i.e., as we shall see below,~$s$ is
    a solution of the Euler-Lagrange equation corresponding to~$\omega$).
  \item Variationally trivial Lagrangians are total divergences, which amounts
    to the equality~$\ker\delta=\im\rmdh^{0,n-1}$.
  \item All null total divergences are total curls, i.e.,
    $\rmdh^{0,n-1}\omega=0$ if and only if~$\omega=\rmdh^{0,n-2}\theta$ for
    some~$\theta\in\Lah^{n-2}(\pi)$.
  \item If~$\psi\in\hat{\kappa}(\pi)$ then the nonlinear differential
    operator~$\Delta_\psi\colon\kappa(\pi)\to\hat{\kappa}(\pi)$ is of the
    form~$\delta\omega$ (i.e., is an Euler-Lagrange operator) if and only
    if~$\ell_\psi=\ell_\psi^*$ (the Helmholtz condition).
  \end{enumerate}
\end{theorem}

\subsection{A parallel with finite-dimensional differential geometry. I}
\label{sec:gjsis_jets:parallel-with-finite}

We want to indicate here a very useful and productive analogy between the
geometry of jet spaces (and, more generally, of differential equations) and
classical differential geometry of finite-dimensional smooth manifolds. This
parallel was exposed by A.M.~Vinogradov first within his philosophy of
\emph{Secondary Calculus} (cf.~\cite{Vinogradov:CAnPDEqSC} and references
therein) and we elaborate it further.

Two points of view on jet spaces as geometrical objects may exist. The first
one is formal, traditional and straightforward. It was described on the
previous pages and treats~$\J^\infty(\pi)$ as a particular case of general
infinite-dimensional manifolds. Such an approach, being by all means necessary
for a rigorous exposition of the theory, actually ignores essential intrinsic
structures of jet spaces actually.

Another viewpoint is completely informal but incorporates these structures
\emph{ab ovo} and allows to reveal new and non-trivial relations and results
just ``translating'' from the language of the classical differential
geometry. To facilitate such translations, we shall compile a sort of a
dictionary.

So, we consider the space~$\J^\infty(\pi)$ endowed with the Cartan
distribution~$\C$ and take for the points of the new ``manifold'' maximal
integral submanifolds of~$\C$. As it was indicated above, they are graphs
of~$j_\infty(s)$, $s\in\Se(\pi)$, and thus new ``points'' are sections
of~$\pi$ (i.e., fields).

Let~$\omega\in\Lah^n(\pi)$ be a horizontal $n$-form on~$\J^\infty(\pi)$ (or a
Lagrangian density). Then to any ``point''~$s\in\Se(\pi)$ we can put into
correspondence the number
\begin{equation*}
  \omega(s)=\int_M j_\infty(s)^*(\omega).
\end{equation*}
Thus, Lagrangians are understood as functions. Due to the identity
\begin{equation*}
  j_\infty(s)^*(\rmdh\omega)=\rmd(j_\infty(s)^*(\omega))
\end{equation*}
and the Stokes formula
\begin{equation*}
  \int_M\rmd\theta=\int_{\partial M}\theta,\qquad\theta\in\La^{n-1}(M),
\end{equation*}
``functions'' of the form~$\rmdh\omega$ vanish at all points. So, no-trivial
functions are elements of the cohomology
group~$E_1^{0,n}(\pi)=H_h^h(\pi)$. Thus, the beginning of the dictionary is
\begin{eqnarray*}
  \text{\textbf{Manifold~$M$}}&&\text{\textbf{Jet space~$\J^\infty(\pi)$}} \\[1ex]
  \text{points}&\quad\longleftrightarrow\quad
  &\text{sections of~$\pi$ (fields)} \\
  \text{functions}&\quad\longleftrightarrow\quad&\text{Lagrangians }\
  \omega=L\rmd x^1\wedge\dots\wedge\rmd x^n\\
  \text{value at a point, }\ f(x)&\quad\longleftrightarrow\quad&
  \text{integral }\
  \omega(s)=\int_Mj_\infty(s)^*\omega \\
  &&\text{(the cohomology class of~$\omega$)}
\end{eqnarray*}

The next step is to define vector fields. They should be infinitesimal
transformations of~$\J^\infty(\pi)$ that preserve the Cartan distribution (or,
equivalently, move ``points'' to ``points''). These are exactly the fields
lying in~$\Xc(\pi)$. But vector fields~$X\in\C\X(\pi)$ (i.e., lying in the
Cartan distribution) are tangent to maximal integral manifolds of the latter
and thus are trivial in the space of fields. Consequently, non-trivial vector
fields are identified with elements of~$\sym\pi$, i.e., with evolutionary
vector fields on~$\J^\infty(\pi)$. On the other hand, as it was indicated
above, they are integrable sections\footnote{By integrable sections we mean
  those ones whose graphs are maximal integral manifolds of the Cartan
  distribution.} of~$\JJ^\infty(\kappa)$. Hence, this bundle can be considered
as the tangent bundle to~$\J^\infty(\pi)$. So, the dictionary can be continued
as follows:
\begin{eqnarray*}
  \text{\textbf{Manifold~$M$}}&&\text{\textbf{Jet space~$\J^\infty(\pi)$}} \\[1ex]
  \text{vector fields,
    $\X(M)$}&\quad\longleftrightarrow\quad&
  \text{evolutionary vector fields, $\kappa(\pi)$} \\
  \text{the tangent bundle}&\quad\longleftrightarrow\quad&
  \text{the bundle of horizontal
    jets $\JJ^\infty(\kappa(\pi))$}
\end{eqnarray*}

Differential forms on a smooth manifold~$M$ may be understood as multi-linear
functions on the space of vector fields (or fibre-wise multi-linear functions
on~$T(M)$): we insert a vector field into a $p$-form and obtain a
$(p-1)$-form. In our context, such objects are exactly elements
of~$\CDiffskalt_{p-1}(\kappa,\hat{\kappa})=E_1^{p,n}(\pi)$. We call them
\emph{variational forms} of degree~$p$ and have the following parallel:
\begin{eqnarray*}
  \text{\textbf{Manifold~$M$}}&&\text{\textbf{Jet space~$\J^\infty(\pi)$}} \\[1ex]
  \text{differential forms,
    $\La^p(M)$}&\quad\longleftrightarrow\quad&\text{variational forms, $\CDiffskalt_{p-1}(\kappa,\hat{\kappa})$} \\
  \text{the de~Rham complex}&\quad\longleftrightarrow\quad&\text{the variational complex}
\end{eqnarray*}

\begin{remark}
  \label{rem:gjsis_jets:4}
  It can also be shown that smooth maps~$\J^\infty(\pi)\to\J^\infty(\pi')$
  preserving the Cartan distributions are completely determined by non-linear
  differential operators from~$\Se(\pi)$ to~$\Se(\pi')$ while the
  differentials of these maps~$\JJ^\infty(\kappa)\to\JJ^\infty(\kappa')$
  correspond to linearizations. Unfortunately, a detailed exposition of this
  parallel is out of scope of our review.
\end{remark}

We shall continue to compile our dictionary in
Subsection~\ref{sec:gjsis_jets:parallel-with-finite-1}.

\subsection{Hamiltonian formalism}
\label{sec:gjsis_jets:hamilt-form}

The objects dual to~$\CDiffskalt_{p-1}(\kappa(\pi),\hat{\kappa}(\pi))$ are the
modules of \emph{variational
  multivectors}~$D_p(\pi)=\CDiffskalt_{p-1}(\hat{\kappa}(\pi),\kappa(\pi))$. In
particular, $D_1(\pi)=\kappa(\pi)$. We also set~$D_0(\pi)=\Lah^n(\pi)/\im
\rmdh^{0,n-1}=E_1^{0,n}(\pi)$. To describe Hamiltonian formalism
on~$\J^\infty(\pi)$, we first introduce the \emph{variational Schouten
  bracket}~\cite{KerstenKrasilshchikVerbovetsky:HOpC}
\begin{equation*}
  \lshad\cdot,\cdot\rshad\colon D_p(\pi)\times D_q(\pi)\to D_{p+q-1}(\pi)
\end{equation*}
in the following way
(cf.~\cite{Krasilshchik:SBCAl,IgoninVerbovetskyVitolo:VMBGJS}, see
also~\cite{Dorfman:DSInNEvEq,GelfandDorfman:SBHOp}). If~$B=[\omega]\in
D_0(\pi)$ is a coset of a horizontal form~$\omega\in\Lah^n(\pi)$ and~$A\in
D_p(\pi)$, $p>0$, then we set
\begin{equation*}
  \lshad A,B\rshad=(-1)^p\lshad B,A\rshad=A(\delta\omega),
\end{equation*}
while for any~$B\in D_q(\pi)$, $q>0$,
and~$\psi=\delta\omega\in\hat{\kappa}(\pi)$, $\omega\in\Lah^n(\pi)$,
\begin{equation*}
  \lshad A,B\rshad(\psi)=\lshad A,B(\psi)\rshad+(-1)^{q-1}\lshad A(\psi),B\rshad,
\end{equation*}
and these two equalities define the bracket completely. In particular,
\begin{equation*}
  \lshad \phi,[\omega]\rshad=[\Ev_\phi(\omega)],\qquad\phi\in
  D_1(\pi)=\kappa(\pi),\quad
  \omega\in\Lah^n(\pi),
\end{equation*}
and
\begin{equation*}
  \lshad \phi,\phi'\rshad=\Ev_\phi(\phi')-\Ev_{\phi'}(\phi)=\{\phi,\phi'\},\qquad
  \phi,\phi'\in D_1(\pi).
\end{equation*}

\begin{proposition}
  \label{prop:gjsis_jets:5}
  The variational Schouten bracket determines a super Lie algebra structure in
  the space~$D(\pi)=\sum_{p\ge 0}D_p(\pi)$ in the following sense:
  \begin{align}
    \label{eq:gjsis_jets:39}
    &\lshad A,B\rshad=-(-1)^{(p-1)(q-1)}\lshad B,A\rshad,\\\label{eq:gjsis_jets:40}
    &(-1)^{(p+1)(r+1)}\lshad\lshad A,B\rshad,C\rshad+(-1)^{(q+1)(p+1)}\lshad\lshad
    B,C\rshad,A\rshad\\
    &\qquad+(-1)^{(r+1)(q+1)}\lshad\lshad C,A\rshad,B\rshad=0\nonumber
  \end{align}
  for all~$A\in D_p(\pi)$, $B\in D_q(\pi)$, $C\in D_r(\pi)$.
\end{proposition}

To compute the Schouten bracket explicitly, for any natural~$n$ consider the
set~$S_n^i$ of all~$(n-i)$-\emph{un-shuffles} consisting of all
permutations~$\sigma$ of the set~$\{1,\dots,n\}$ such that
\begin{equation*}
  \sigma(1)<\dots<\sigma(i),\qquad\sigma(i+1)<\dots<\sigma(n).
\end{equation*}
We formally set~$S_n^i=\varnothing$ for~$i<0$ and~$i>n$. We also use a short
notation~$\psi_{\sigma(k_1,k_2)}$
for~$\psi_{\sigma(k_1)},\dots,\psi_{\sigma(k_2)}$.

Let now~$A\in D_p(\pi)$ and~$B\in D_q(\pi)$. Then for
any~$\psi_1,\dots,\psi_{p+q-1}\in\hat{\kappa}(\pi)$ we have
\begin{multline}\label{eq:gjsis_jets:41}
  \lshad A,B\rshad(\psi_1,\dots,\psi_{p+q-2}) =\sum_{\sigma\in
  S_{p+q-2}^{q-1}}(-1)^\sigma
  \ell_{B,\psi_{\sigma(1,q-1)}}(A(\psi_{\sigma(q,p+q-2)})) \\
  -(-1)^{(p-1)(q-1)}\sum_{\sigma\in S_{p+q-2}^p}(-1)^\sigma
  B(\ell^*_{A,\psi_{\sigma(1,p-1)}}(\psi_{\sigma(p)}),
  \psi_{\sigma(p+1,p+q-2)}) \\
  -(-1)^{(p-1)(q-1)}\sum_{\sigma\in S_{p+q-2}^{p-1}}(-1)^\sigma
  \ell_{A,\psi_{\sigma(1,p-1)}}(B(\psi_{\sigma(p,p+q-2)})) \\
  +\sum_{\sigma\in S_{p+q-2}^q}(-1)^\sigma
  A(\ell^*_{B,\psi_{\sigma(1,q-1)}}
  (\psi_{\sigma(q)}),\psi_{\sigma(q+1,p+q-2)}),
\end{multline}
where~$(-1)^\sigma$ stands for the parity of the permutation~$\sigma$ and, as
before,~$\ell_{\Delta,\psi_1,\dots,\psi_k}(\phi)=\Ev_\phi(\Delta)(\psi_1,\dots,\psi_k)$.

We say that a bivector~$A\in
D_2(\pi)=\CDiffskalt(\hat{\kappa}(\pi),\kappa(\pi))$ is a \emph{Hamiltonian
  structure} on~$\J^\infty(\pi)$ if
\begin{equation}\label{eq:gjsis_jets:42}
  \lshad A,A\rshad=0.
\end{equation}

\begin{remark}
  A more appropriate name is a \emph{Poisson} structure but we follow here the
  tradition accepted in the theory of integrable systems.
\end{remark}

Given a Hamiltonian structure, one can define a \emph{Poisson bracket} (with
respect to the Hamiltonian structure~$A$) on the set of Lagrangians:
\begin{equation}\label{eq:gjsis_jets:43}
  \{\omega,\omega'\}_A=\langle A(\delta(\omega)),\delta(\omega')\rangle,\qquad
  \omega,\omega'\in D_0(\pi).
\end{equation}
Two elements are \emph{in involution} (with respect to the structure~$A$) if
\begin{equation*}
  \{\omega,\omega'\}_A=0.
\end{equation*}

\begin{proposition}
  \label{prop:gjsis_jets:6}
  For any~$A\in\CDiffskalt(\hat{\kappa}(\pi),\kappa(\pi))$ one has
  \begin{equation}\label{eq:gjsis_jets:44}
    \{\omega,\omega'\}_A=-\{\omega',\omega\}_A
  \end{equation}
  If in addition~$A$ satisfies~\eqref{eq:gjsis_jets:42} then
  \begin{equation}\label{eq:gjsis_jets:45}
    \{\omega,\{\omega',\omega''\}_A\}_A+
    \{\omega',\{\omega'',\omega\}_A\}_A+
    \{\omega'',\{\omega,\omega'\}_A\}_A=0.
  \end{equation}
  One also has
  \begin{equation}\label{eq:gjsis_jets:46}
    A_{\{\omega,\omega'\}_A}=\{A_\omega,A_{\omega'}\},
  \end{equation}
  where~$A_\omega=A(\delta(\omega))\in\kappa(\pi)$ and the curlies in the
  right-hand side denote the Jacobi bracket.
\end{proposition}

Chose a Hamiltonian structure~$A$ and consider the sequence of operators
\begin{equation}\label{eq:gjsis_jets:47}
  \fl\qquad
  0\xrightarrow{}D_0(\pi)\xrightarrow{\partial_A}D_1(\pi)\xrightarrow{}
  \cdots\xrightarrow{}D_p(\pi)\xrightarrow{\partial_A}D_{p+1}(\pi)
  \xrightarrow{}\cdots,
\end{equation}
where~$\partial_A(B)=\lshad A,B\rshad$.

\begin{proposition}
  \label{prop:gjsis_jets:7}
  Sequence~\eqref{eq:gjsis_jets:47} is a complex,
  i.e.,~$\partial_A\circ\partial_A=0$.
\end{proposition}

We say that~$\Ev_\phi$ is a \emph{Hamiltonian vector field}
if~$\phi\in\ker\partial_A$. This is equivalent to
\begin{equation}\label{eq:gjsis_jets:48}
  \Ev_\phi(\{\omega,\omega'\}_A)=\{\Ev_\phi(\omega),\omega'\}_A
  +\{\omega,\Ev_\phi(\omega')\}_A,
\end{equation}
i.e., $\Ev_\phi$ preserves the Poisson bracket. Due to
Proposition~\ref{prop:gjsis_jets:7}, a particular case of Hamiltonian fields
are fields of the form~$\Ev_{A(\delta(\omega))}$. In this case,~$\omega\in
D_0(\pi)$ is called the \emph{Hamiltonian} of the field under consideration.

We say that~$\omega\in D_0(\pi)$ is a \emph{first integral} of a
Hamiltonian field~$\Ev_\phi$ if~$\Ev_\phi(\omega)=0$. A Hamiltonian
field~$\Ev_{\phi'}$ is a \emph{symmetry} for the field~$\Ev_\phi$
if~$[\Ev_{\phi'},\Ev_\phi]=0$, or
\begin{equation*}
  (\Ev_\phi-\ell_\phi)(\phi')=0.
\end{equation*}

\begin{proposition}
  \label{prop:gjsis_jets:8}
  If~$\Ev_\phi$ is a Hamiltonian vector field with respect to a Hamiltonian
  structure~$A$ then the operator~$A\circ\delta=\partial_A\colon D_0(\pi)\to
  D_1(\pi)$ takes first integrals of~$\Ev_\phi$ to its symmetries.
\end{proposition}

A Hamiltonian structure~$B\in D_2(\pi)$ is said to be \emph{compatible} with
the structure~$A$ if~$B\in\ker\partial_A$, or
\begin{equation}\label{eq:gjsis_jets:49}
  \lshad A,B\rshad=0.
\end{equation}
This is equivalent to the fact that all the bivectors
\begin{equation}\label{eq:gjsis_jets:50}
  \lambda A+\mu B,\qquad \lambda,\mu\in\R,
\end{equation}
are Hamiltonian structures on~$\J^\infty(\pi)$. The
family~\eqref{eq:gjsis_jets:50} is called a \emph{Poisson pencil}. When two
Hamiltonian structures are given, one also says that they form a
\emph{bi-Hamiltonian structure}.

\begin{coordinates}
  Let us indicate how to verify conditions~\eqref{eq:gjsis_jets:42}
  and~\eqref{eq:gjsis_jets:49} in coordinates (the explanation will be given
  below in Subsection~\ref{sec:gjsis_jets:parallel-with-finite-1}, see
  Remark~\ref{rem:gjsis_jets:5}). Take bivectors~$A$, $B\in
  D_2(\pi)=\CDiffskalt_1(\hat{\kappa}(\pi),\kappa(\pi))$. Then~$A$ and~$B$, in
  adapted coordinates in~$\J^\infty(\pi)$, are represented as matrix
  $\C$-differential operators
  \begin{equation*}
    A=\left\Vert\sum_\sigma a_\sigma^{ij}D_\sigma\right\Vert,\qquad
    B=\left\Vert\sum_\sigma b_\sigma^{ij}D_\sigma\right\Vert,
  \end{equation*}
  where~$i$, $j=1,\dots, m=\dim\pi$. Let us put into correspondence to these
  operators the functions
  \begin{equation}\label{eq:gjsis_jets:58}
    W_A=\sum_{\sigma,i,j} a_\sigma^{ij}p_\sigma^ip^j,\qquad
    W_B=\sum_{\sigma,i,j} b_\sigma^{ij}p_\sigma^ip^j,
  \end{equation}
  where~$p_\sigma^i$ are \emph{odd} variables. Then~$A$ is a Hamiltonian
  structure if and only if
  \begin{equation}\label{eq:gjsis_jets:51}
    \delta\left(\sum_i\frac{\delta W_A}{\delta u^i}\frac{\delta W_A}{\delta p^i}\right)=0,
  \end{equation}
  while two Hamiltonian structures are compatible if and only if
  \begin{equation}\label{eq:gjsis_jets:52}
    \delta\left(\sum_i\left(\frac{\delta W_A}{\delta u^i}\frac{\delta W_B}{\delta p^i}+\frac{\delta W_B}{\delta u^i}\frac{\delta W_A}{\delta p^i}\right)\right)=0.
  \end{equation}
\end{coordinates}

\begin{theorem}[the Magri Scheme, see~\cite{Magri:SMInHEq,Krasilshchik:SBCAl}]
  \label{thm:gjsis_jets:4}
  Let~$(A,B)$ be a bi-Hamiltonian structure on~$\J^\infty(\pi)$ and assume
  that the complex~\eqref{eq:gjsis_jets:47} is acyclic in the term~$D_1(\pi)$,
  i.e., every Hamiltonian vector field with respect to~$A$ possesses a
  Hamiltonian. Assume also that two densities $\omega_1$, $\omega_2\in
  D_0(\pi)$ are given, such that
  $\partial_{A}(\omega_1)=\partial_{B}(\omega_2)$. Then:
  \begin{enumerate}
  \item There exist elements $\omega_3,\dots,\omega_s,\dots\in D_0(\pi)$
    satisfying
    \begin{equation}\label{eq:gjsis_jets:60}
      \partial_{A}(\omega_s)=\partial_{B}(\omega_{s+1}),\qquad
      s=2,3,\dots
    \end{equation}
  \item All elements $\omega_1,\dots,\omega_s,\dots$ are in involution with
    respect to both Hamiltonian structures, i.e.,
    \begin{equation*}
      \{\omega_\alpha,\omega_\beta\}_A=\{\omega_\alpha,\omega_\beta\}_B=0
    \end{equation*}
    for all $\alpha$, $\beta\ge1$.
  \end{enumerate}
\end{theorem}

\begin{example}[the KdV hierarchy]
  \label{exmp:gjsis_jets:1}
  Consider~$\J^\infty(\pi)$ for the trivial one-dimensional
  bundle~$\pi\colon\R\times\R\to\R$. Let~$x$ be the independent variable
  and~$u$ be the fibre coordinate (the unknown function). Then the adapted
  coordinates~$u=u_0$, $u_1,\dots,u_k,\dots$ in~$\J^\infty(\pi)$ arise,
  where~$u_k$ corresponds to~$\partial^ku/\partial x^k$. The operators
  \begin{equation*}
    A=D_x=\frac{\partial}{\partial x}
    +\sum_{k=0}^\infty u_{k+1}\frac{\partial}{\partial u_k}
  \end{equation*}
and
\begin{equation*}
  B=D_x^3+4uD_x+2u_1,
\end{equation*}
as it can be easily checked using~\eqref{eq:gjsis_jets:51}
and~\eqref{eq:gjsis_jets:52}, constitute a bi-Hamiltonian structure
on~$\J^\infty(\pi)$. Then obviously for the horizontal forms
\begin{equation*}
  \omega_1=\frac{1}{2}u^2\rmd x,\qquad\omega_2=\frac{1}{2}u\rmd x
\end{equation*}
one has~$\partial_A\omega_1=A(u)=u_1$ and~$\partial_B\omega_2=B(1)=u_1$, i.e.,
\begin{equation*}
  \partial_A\omega_1=\partial_B\omega_2.
\end{equation*}
The first cohomology group of~$\partial_A$ is trivial
(see~\cite{Getzler:DTHOpFCV}), and consequently we obtain an infinite series
of first integrals and the corresponding symmetries. The second,
after~$u_1$, symmetry is~$6uu_1+u_3$:
\begin{equation*}
  6uu_1+u_3=\partial_A((u^3-\frac{1}{2}u_1^2)\rmd
  x)=\partial_B(\frac{1}{2}u^2\rmd x).
\end{equation*}
The corresponding flow on~$\J^\infty(\pi)$ is governed by the evolution
equation
\begin{equation*}
  u_t=6uu_x+u_{xxx};
\end{equation*}
thus, we obtain the Korteweg-de Vries equation and the corresponding hierarchy
of commuting flows (the higher KdV equations). The entire family of commuting
flows can be obtained by applying the \emph{Lenard recursion operator}
(see~\cite{GardnerGreeneKruskalMiura:KEqGVMExS})
\begin{equation}\label{eq:gjsis_jets:53}
  R=B\circ A^{-1}=D_x^2+4u+2u_1D_x^{-1}
\end{equation}
to the right-hand side of the first flow~$u_t=u_x$ sufficiently many times.
\end{example}

\begin{example}[the Boussinesq hierarchy]
  \label{exmp:gjsis_jets:2}
  Consider the adapted coordinates~$x$, $u$, $v,\dots,u_k$, $v_k,\dots$ in the
  space~$\J^\infty(\pi)$, where~$\pi\colon\R^2\times\R\to\R$ is the trivial
  two-dimensional bundle over~$\R$. Then the operators
  \begin{equation*}
    A=
    \begin{pmatrix}
      0&D_x\\
      D_x&0
    \end{pmatrix},\qquad
    B=
    \begin{pmatrix}
      \sigma D_x^3+uD_x+\frac{1}{2}u_1&\frac{1}{2}vD_x\\
      \frac{1}{2}vD_x+\frac{1}{2}v_1&D_x
    \end{pmatrix},
  \end{equation*}
  where~$\sigma$ is a real constant, form a bi-Hamiltonian structure. For the
  differential forms~$\omega_1=2(u+v)\rmd x$ and~$\omega_2=uv\rmd x$ one
  obviously has~$\partial_A\omega_1=\partial_B\omega_2$. The arising hierarchy
  of commuting flows corresponds to the evolution equation
  \begin{equation}\label{eq:gjsis_jets:54}
    \begin{array}{l}
      u_t  =u_xv+uv_x+\sigma v_{xxx},\\
      v_t  =u_x+vv_x
    \end{array}
  \end{equation}
  which is the \emph{two-component Boussinesq system} which can be obtained
  from the Kaup equation, see~\cite{Kaup:HOrWWEqMSIt}. Note that there exists
  another Hamiltonian operator 
  \begin{equation*}
    C=
    \begin{pmatrix}
      C^{uu}&C^{uv}\\
      C^{vu}&C^{vv}
    \end{pmatrix},
  \end{equation*}
  where
  \begin{align*}
    C^{uu}&=\sigma D_x^3+\frac{3}{2}\sigma v_1D_x^2+(\sigma
    v_2+uv)D_x+\frac{1}{2}(\sigma v_3+uv_1+u_1v),\\
    C^{uv}&=\sigma D_x^3+(u+\frac{1}{4}v^2)D_x+\frac{1}{2}u_1,\\
    C^{vu}&=\sigma D_x^3+(u+\frac{1}{4}v^2)D_x+\frac{1}{2}(u_1+vv_1),\\
    C^{vv}&=vD_x+\frac{1}{4}v_1,
  \end{align*}
which is compatible both with~$A$ and~$B$. In this sense,
system~\eqref{eq:gjsis_jets:54} is \emph{tri-Hamiltonian}.
\end{example}

\begin{example}[the KdV hierarchy, II]
  \label{exmp:gjsis_jets:3}
  Let~$\pi\colon\R^3\times\R^1\to\R^1$ be the trivial three-dimensional bundle
  with the coordinates~$t$ in the base and~$u$, $v$, $w$ in the
  fibre. Introduce the adapted coordinates~$u_k$, $v_k$, $w_k$, where
  $k=0,1,2,\dots$, in~$\J^\infty(\pi)$ and consider the operators
  \begin{equation*}
    A=
    \begin{pmatrix}
      0&-1&0\\
      1&0&-6u\\
      0&6u&D_t
    \end{pmatrix},\quad
    B=
    \begin{pmatrix}
      0&-2u&-D_t-2v\\
      2u&D_t&-12u^2-2w\\
      -D_t+2v&12u^2+2w&8uD_t+4u_1
    \end{pmatrix},
  \end{equation*}
  which form a bi-Hamiltonian structure on~$\J^\infty(\pi)$. It can be easily
  seen that~$\partial_A\omega_1=\partial_B\omega_2$, where
  \begin{equation*}
    \omega_1=(uw-\frac{1}{2}v^2+2u^3)\rmd t,\qquad
    \omega_2=-(\frac{3}{2}u^2+\frac{1}{2}w)\rmd t.
  \end{equation*}
  Thus, we obtain a hierarchy of commuting flows whose first term is
  \begin{equation*}
    u_x=v,\qquad v_x=w,\qquad w_x=u_t-6uv,
  \end{equation*}
  which is obviously the KdV equation rewritten in a different way (cf.~the
  paper~\cite{Tsarev:HPSIEqCMMMP}).

  Note that the Lenard recursion operator~\eqref{eq:gjsis_jets:53} in the new
  representation of the KdV hierarchy acquires the form
  \begin{equation*}
    R=
    \begin{pmatrix}
      0&-2u&-D_t-2v\\
      2u&D_t&-12u^2-2w\\
      -D_t+2v&12u^2-2w&8uD_t+4u_1
    \end{pmatrix}
    \circ
    \begin{pmatrix}
      -36uD_t^{-1}\circ u&1&-6uD_t^{-1}\\
      -1&0&0\\
      6D_t^{-1}\circ u&0&D_t^{-1}
    \end{pmatrix}.
  \end{equation*}
\end{example}

\begin{remark}
  In a recently published paper~\cite{BarakatSoleKac:PVAlTHEq}, the authors
  formulate a much weaker than triviality of the first Poisson cohomology
  group criterion for feasibility of the Magri scheme. The criterion is given
  in the framework of Dirac structures~\cite{Dorfman:DSInNEvEq,Courant:DM}
  (though it also admits a self-contained formulation) that unfortunately lie
  beyond the scope of our review. Note that Dirac structures that merge the
  notions of symplectic and Hamiltonian operators constitute an interesting
  object for geometrical research in PDEs.
\end{remark}

\begin{remark}
  Normal forms for Hamiltonian operators of order~$\le 5$ and a
  ``variational'' analog of the Darboux Lemma were presented
  in~\cite{Astashov:NFHOpFT,AstashovVinogradov:SHOpFT}. In~\cite{SoleKacWakimoto:CPVAl},
  one can find normal forms for operators of order~$\le 11$ and some
  classification results for operators of higher order.
\end{remark}

\subsection{A parallel with finite-dimensional differential geometry. II}
\label{sec:gjsis_jets:parallel-with-finite-1}

Let us now continue to construct the dictionary started in
Section~\ref{sec:gjsis_jets:parallel-with-finite}.

Of course, variational multivectors introduced in
Section~\ref{sec:gjsis_jets:hamilt-form} are exact counterparts of classical
multivector fields in differential geometry. These fields are naturally
understood as smooth functions on~$T^*M$ skew-symmetric and multi-linear
with respect to fibre variables. Exactly the same interpretation is valid for
variational vectors if one considers the
bundle~$\tau^*\colon\JJ^\infty(\hat{\kappa})\to\J^\infty(\pi)$\footnote{Note
  that in such a way we independently arrived to Kupershmidt's notion of the
  cotangent bundle to a vector bundle, see~\cite{Kupershmidt:GJBSLHF}.}. Thus,
we have the following translations:
\begin{eqnarray*}
  \text{\textbf{Manifold~$M$}}&&\text{\textbf{Jet space~$\J^\infty(\pi)$}} \\[1ex]
  \text{multivector fields}&\quad\longleftrightarrow\quad&\text{variational multivectors, $\CDiffskalt_{p-1}(\hat{\kappa},\kappa)$} \\
  \text{Schouten bracket}&\quad\longleftrightarrow\quad&\text{variational Schouten bracket} \\
  \text{Poisson structure}&\quad\longleftrightarrow\quad&\text{Hamiltonian operator} \\
  \text{cotangent bundle,}&\quad\longleftrightarrow\quad&\text{variational cotangent bundle,}\\
  \qquad\qquad T^*M\to M&\quad\longleftrightarrow\quad&%\qquad\qquad
  \JJ^\infty(\hat{\kappa})\to\J^\infty(\pi)\end{eqnarray*}

As it is well known, the tangent space~$T^*M$ is endowed with the natural
symplectic structure~$\Omega=\rmd p\wedge\rmd q\in\La^2(M)$ and, in
addition,~$\Omega=\rmd\rho$, where the form~$\rho=p\rmd q$ is defined
invariantly as well. Similar constructions exist
on~$\JJ^\infty(\hat{\kappa})$. Let us show this.

To this end, recall (see Remark~\ref{rem:gjsis_jets:1}) that the
manifold~$\JJ^\infty(\hat{\kappa})$ is diffeomorphic
to~$\J^\infty(\pi\times_M\hat{\pi})$,
where~$\hat{\pi}=\Hom(\pi,\bigwedge^nT^*M)$. Hence, the module of variational
$1$-forms on~$\JJ^\infty(\hat{\pi})$ is isomorphic to
\begin{equation}\label{eq:gjsis_jets:55}
  \hat{\kappa}(\pi\times_M\hat{\pi})=
  \hat{\kappa}(\pi)\times_{\J^\infty(\pi)}\kappa(\pi).
\end{equation}
Then the $1$-form~$\rho_\pi$ (the analog of~$p\rmd q$) is uniquely defined
by the condition
\begin{equation}\label{eq:gjsis_jets:56}
  j_h^\infty(\psi)^*(\rho_\pi)=(\psi,0),
\end{equation}
where~$\psi\in\hat{\kappa}(\pi)$ is an arbitrary variational $1$-form
on~$\J^\infty(\pi)$.

Now, by the same reasons and dually to~\eqref{eq:gjsis_jets:55}, the module of
vector fields on~$\JJ^\infty(\hat{\kappa})$ is
\begin{equation*}
  \kappa(\pi\times_M\hat{\pi})=\kappa(\pi)\times_{\J^\infty(\pi)}\kappa(\hat{\pi}).
\end{equation*}
Thus we see that the symplectic structure~$\Omega_\pi$ must be an element of
the
module~$\CDiffskalt(\kappa(\pi\times_M\hat{\pi}),\hat{\kappa}(\pi\times_M\hat{\pi}))$. For
any element~$(\phi,\psi)\in\kappa(\pi\times_M\hat{\pi})$ we set
\begin{equation}\label{eq:gjsis_jets:57}
  \Omega_\pi(\phi,\psi)=(-\psi,\phi);
\end{equation}
this is a skew-adjoint $\C$-differential operator of order~$0$.

\begin{remark}
  The form~$\rho_\pi$ can be defined in a different way. Namely,
  let~$X=(\phi,\psi)$ be a vector field on~$\JJ^\infty(\hat{\kappa})$. Then
  its 1st component is a vector field vertical with respect to the
  projection~$\JJ^\infty(\hat{\kappa})\to\J^\infty(\pi)$ and may be understood
  as a function on~$\JJ^\infty(\hat{\kappa})$. Then we
  set~$\rho_\pi(X)=\phi$. This definition is equivalent
  to~\eqref{eq:gjsis_jets:56}.
\end{remark}

Of course, the operator~$\Omega_\pi$ is invertible and the inverse
one~$S_\pi=\Omega_\pi^{-1}$ is a bivector on the cotangent
manifold~$\JJ^\infty(\hat{\kappa})$. There exist two points of view at this
manifold (as well as at~$T^*M$). The first one treats it as a classical
(``even'') manifold. The second approach considers~$\JJ^\infty(\hat{\kappa})$
as a super-manifold with odd coordinates along the
projection~~$\JJ^\infty(\hat{\kappa})\to\J^\infty(\pi)$ and even ones in the
base. If one takes the first approach the bivector~$S_\pi$ will define the
Poisson bracket for functions on~$\JJ^\infty(\hat{\kappa})$. The second
approach leads to functions multi-linear and skew-symmetric with respect to
fibre variables. As it was stated above, these functions are identified with
variational multivectors on~$\J^\infty(\pi)$. Then the super-bracket defined
by~$S_\pi$ coincides with the Schouten bracket. The bracket is given by the
formula
\begin{equation}
  \label{eq:gjsis_jets:59}
  S_\pi(\delta\omega_1)(\omega_2)
  =\langle S_\pi(\delta\omega_1),\delta\omega_2\rangle=\left\{
  \begin{array}{ll}
    \{\omega_1,\omega_2\}, &\text{the even case,}\\
    \lshad\omega_1,\omega_2\rshad, &\text{the odd case,}
  \end{array}\right.
\end{equation}
in both cases.

\begin{remark}
  \label{rem:gjsis_jets:5}
  The above said clarifies the meaning of
  formulas~\eqref{eq:gjsis_jets:58}--\eqref{eq:gjsis_jets:52}. Namely, the
  correspondence~$A\mapsto W_A$ given by~\eqref{eq:gjsis_jets:58} describes
  how to construct the function on~$\JJ^\infty(\hat{\kappa})$ when a
  bivector~$A$ is given (to be more precise, this function is the horizontal
  cohomology class of the form~$W_A\rmd x^1\wedge\dots\wedge\rmd x^n$
  in~$H_h^n(\pi)$). The argument of~$\delta$ in~\eqref{eq:gjsis_jets:53} is
  the coordinate expression of the bracket~\eqref{eq:gjsis_jets:59} in the odd
  case while~\eqref{eq:gjsis_jets:59} itself checks triviality of its
  horizontal cohomology class.
\end{remark}

\section{Differential Equations}
\label{sec:gjsis_eqs:diff-equat}

With the concept of the jet bundle at our disposal we give a geometric
definition of differential equations.

\subsection{Definition of differential equations}
\label{sec:gjsis_eqs:defin-diff-equat}

Suppose we have a system
\begin{equation}
  \label{eq:gjsis_eqs:1}
  F_s(x^i,u^j_I)=0,\qquad s=1,\dots,l,
\end{equation}
of partial differential equations in $n$ independent variables~$x^i$
and $m$ dependent variables~$u^j$.  Equations~\eqref{eq:gjsis_eqs:1}
determine a locus in the jet space~$\J^{\infty}(\pi)$ of a vector
bundle~$\pi\colon E\to M$, such that $\dim E=m+n$, $\dim M=n$.

The subset of~$\J^{\infty}(\pi)$ defined in this way is not an
adequate geometric construction corresponding to the system at hand,
because it does not take into account differential consequences
of~\eqref{eq:gjsis_eqs:1}.  So, we extend~\eqref{eq:gjsis_eqs:1} to a
larger system
\begin{equation}
  \label{eq:gjsis_eqs:2}
  D_I(F_s)=0\qquad\text{for all multi-indices  $I$ and $s=1,\dots,l$,}
\end{equation}
and consider the locus $\mathcal{E}\subset\J^{\infty}(\pi)$ defined
by~\eqref{eq:gjsis_eqs:2}.

Thus, we get a correspondence
\begin{equation*}
  F_s(x^i,u^j_I)=0\quad\text{~\eqref{eq:gjsis_eqs:1}}\qquad
  \mapsto\qquad\mathcal{E}\subset\J^{\infty}(\pi).
\end{equation*}
This correspondence behaves nice with respect to solutions
of~\eqref{eq:gjsis_eqs:1}: they are those sections of~$\pi$ whose
infinite jets lie in~$\mathcal{E}$.  To put this another way, the
solutions of~\eqref{eq:gjsis_eqs:1} are the maximal integral
submanifolds of the Cartan distribution restricted to~$\mathcal{E}$.

This shows that $\mathcal{E}$ endowed with the Cartan distribution can
be taken as the geometric object corresponding to
system~\eqref{eq:gjsis_eqs:1}, we call such a manifold~$\mathcal{E}$ an
\emph{equation}.

An equation is generally of infinite dimension.

Without loss of generality we assume that~\eqref{eq:gjsis_eqs:1}
does not contain equations of zero order, in geometric language this
means that the projection
$\eval{\pi_{\infty,0}}_{\mathcal{E}}\colon\mathcal{E}\to\J^0(\pi)$ is
surjective.  It is obvious that at every point $\theta\in\mathcal{E}$
the Cartan plane is tangent to the equation, $\C_{\theta}\subset
T_{\theta}(\mathcal{E})$, so that the dimension of the Cartan
distribution on an equation is equal to~$n$, the same as on the jet
space.

System~\eqref{eq:gjsis_eqs:1} is a coordinate description of an
equation~$\mathcal{E}$.  Every equation has many different coordinate
descriptions.

\begin{example}
  \label{exmp:gjsis_eqs:1}
  Take the KdV equation (cf.~Example~\ref{exmp:gjsis_jets:1} from
  Section~\ref{sec:gjsis_jets:jet-spaces})
  \begin{equation}
    \label{eq:gjsis_eqs:6}
    u_t-6uu_x-u_{xxx}=0.
  \end{equation}
  The bundle~$\pi$ here is the projection $\pi\colon\R^3\to\R^2$,
  $(x,t,u)\mapsto(x,t)$.  The jet space $\J^{\infty}(\pi)$ has
  coordinates $x$,~$t$,~$u$,~$u_x$,~$u_t$,~\dots\,,~$u_I$,~\dots{} The
  equation $\mathcal{E}\subset\J^{\infty}(\pi)$ is given by the
  infinite series of equations
  \begin{align*}
    u_t&=6uu_x+u_{xxx}, \\
    u_{tx}&=6u_x^2+6uu_{xx}+u_{xxxx}, \\
    &\cdots \\
    u_{tI}&=D_I(6uu_x+u_{xxx}), \\
    &\cdots
  \end{align*}
  The Cartan distribution on~$\mathcal{E}$ is two-dimensional and generated by
  \begin{align*}
    D_x&=\frac{\partial}{\partial x}
    +\sum_s u_{s+1}\frac{\partial}{\partial u_s}, \\
    D_t&=\frac{\partial}{\partial t}
    +\sum_s D_x^s(6uu_x+u_{xxx})\frac{\partial}{\partial u_s},
  \end{align*}
  where $u_s=u_{x\dots x}$ ($s$ times).  The functions $x$,~$t$,~$u_s$ can be
  taken to be coordinates on~$\mathcal{E}$.

  The system
  \begin{equation}
    \label{eq:gjsis_eqs:5}
    u_x-v=0, \quad v_x-w=0, \quad w_x-u_t+6uv=0
  \end{equation}
  gives rise to the same equation $\mathcal{E}\subset\J^{\infty}(\pi')$ with
  $\pi'\colon\R^5\to\R^2$, $(x,t,u,v,w)\mapsto(x,t)$.  To prove
  that~\eqref{eq:gjsis_eqs:6} and~\eqref{eq:gjsis_eqs:5} define the same
  equation consider the map
  \begin{align*}
    &a\colon\J^2(\pi)\to\J^0(\pi'),\qquad
    &&a(x,t,u,u_x,u_t,u_{xx},u_{xt},u_{tt})=(x,t,u,u_x,u_{xx}), \\
    &b\colon\J^0(\pi')\to\J^0(\pi),\qquad
    &&b(x,t,u,v,w)=(x,t,u).
  \end{align*}
  Let $a^{\infty}\colon\J^{\infty}(\pi)\to\J^{\infty}(\pi')$ and
  $b^{\infty}\colon\J^{\infty}(\pi')\to\J^{\infty}(\pi)$ be the lifts of these
  maps (cf.~Remark~\ref{sec:gjsis_jets:Lie-Back}).  Then it is easy to see
  that $\eval{a^{\infty}}_{\mathcal{E}}\circ
  \eval{b^{\infty}}_{\mathcal{E}}=\eval{b^{\infty}}_{\mathcal{E}}\circ
  \eval{a^{\infty}}_{\mathcal{E}}=\eval{\id}_{\mathcal{E}}$.  Lifts preserve
  the Cartan distributions, hence the maps $\eval{a^{\infty}}_{\mathcal{E}}$
  and $\eval{b^{\infty}}_{\mathcal{E}}$ are isomorphisms of equations
  determined by~\eqref{eq:gjsis_eqs:6} and \eqref{eq:gjsis_eqs:5}.

  Thus, two different coordinate expressions~\eqref{eq:gjsis_eqs:6} and
  \eqref{eq:gjsis_eqs:5} determine the same equation~$\mathcal{E}$ included
  into two different jet spaces:
  \begin{equation*}
    \xymatrixcolsep{1pc}
    \xymatrixrowsep{0.5pc}
    \xymatrix{
      &J^{\infty}(\pi) \\
      \mathcal{E}\ar[ur]\ar[dr] \\
      &J^{\infty}(\pi').
    }
  \end{equation*}
\end{example}

As usual, the choice of functions~\eqref{eq:gjsis_eqs:1} should be
restricted by some \emph{regularity} assumptions.  Namely, we require
that there exists a subset $\Sigma\subset\{D_I(F_s)\}$ of functions
on~$\J^{\infty}(\pi)$ such that
\begin{enumerate}
\item $F_s\in\Sigma$ for all $s=1$,~\dots,~$l$;
\item the functions that belong to~$\Sigma$ define the
  equation~$\mathcal{E}$;
\item the differentials $\rmd f$ are locally
  linearly independent on~$\mathcal{E}$ for all~$f\in\Sigma$.
\end{enumerate}
We always require these conditions to be satisfied.

In this review, we shall assume that an equation at
hand~$\mathcal{E}\subset\J^{\infty}(\pi)$ is globally defined by a
relation~$F=0$, where $F$ is a section of an appropriate
$l$-dimensional vector bundle~$\xi$ over the jet
space~$\J^{\infty}(\pi)$.  This is always possible.
Equations~\eqref{eq:gjsis_eqs:1} are local coordinate expressions
for~$F=0$.  Denote by~$P$ the $\F(\pi)$-module of sections of~$\xi$,
so that $F\in P$.

The above regularity conditions imply the following very useful fact:
a function $f\in\J^{\infty}(\pi)$ vanishes on~$\mathcal{E}$,
$\eval{f}_{\mathcal{E}}=0$, if and only if $f=\Delta(F)$ for some
$\C$-differential operator~$\Delta\colon P\to\F(\pi)$.

Let $\mathcal{E}=\{F=0\}$ be the equation defined by a section $F\in
P$.  The section $F$ is called \emph{normal} if for any
$\C$-differential operator~$\Delta\colon P\to\F(\pi)$ such that
$\Delta(F)=0$ we have $\eval{\Delta}_{\mathcal{E}}=0$.  The
equation~$\mathcal{E}\subset\J^{\infty}(\pi)$ is called \emph{normal}
if it can be defined by a normal section.

\begin{example}
  A simple example of an abnormal equation is the system
  \begin{equation*}
    u_y-v_x=0,\qquad u_z-w_x=0,\qquad v_z-w_y=0.
  \end{equation*}
  Gauge equations, including Maxwell, Yang-Mills, and Einstein equations, are
  not normal as well.  Such equations are beyond the scope of our review.  On
  the other hand, the majority of equations of mathematical physics, in
  particular all evolution equations, are normal.
\end{example}

The Cartan vector fields on an equation~$\mathcal{E}$ form a Lie
algebra~$\C\X(\mathcal{E})$.  In the same manner as for jet spaces we
define the Lie algebra of \emph{symmetries} of~$\mathcal{E}$ and
spaces of \emph{the Cartan} and \emph{horizontal forms}:
\begin{align*}
  \sym\mathcal{E}&=\Xc(\mathcal{E})/\C\X(\mathcal{E}), \\
  &\text{where }\Xc(\mathcal{E})
  =\sd{X\in\X(E)}{[X,\C\X(\mathcal{E})]\subset\C\X(\mathcal{E})}, \\
  \Lac^p(\mathcal{E})
  &=\sd{\omega\in\La^p(\mathcal{E})}{i_X(\omega)
    =0\quad\forall X\in\C\X(\mathcal{E})}, \\
  \Lah^1(\mathcal{E})
  &=\La^1(\mathcal{E})/\Lac^1(\mathcal{E}), \\
  \Lah^q(\mathcal{E})&=\Lah^1(\mathcal{E})\wedge\dots\wedge\Lah^1(\mathcal{E}).
\end{align*}
We shall discuss them in more detail below.

A morphism of equations $f\colon\mathcal{E}\to\mathcal{E}'$ is a
smooth map that respects the Cartan distribution, i.e., for all points
$\theta\in\mathcal{E}$ we have
$f_*(\C_{\theta})\subset\C_{f(\theta)}$, where $\C_{\theta}$ is the
Cartan plane at a point $\theta\in\mathcal{E}$.

If a morphism $f\colon\mathcal{E}\to\mathcal{E}'$ is a fibre bundle
and the map $f_*\colon\C_{\theta}\to\C_{f(\theta)}$ is an isomorphism
of vector spaces for all points $\theta\in\mathcal{E}$ then $f$ is
called a \emph{covering}.

We shall discuss the theory of covering in
Section~\ref{sec:gjsis_nonloc:nonlocal-theory}.

\begin{example}
  \label{exmp:gjsis_eqs:3}
  For an equation~$\mathcal{E}$ we construct the quotient bundle
  \begin{equation*}
    \tau\colon\T(\mathcal{E})=T(\mathcal{E})/\C\to\mathcal{E},
  \end{equation*}
  where $T(\mathcal{E})\to\mathcal{E}$ is the tangent bundle to~$\mathcal{E}$,
  $\C\subset T(\mathcal{E})$ is the Cartan distribution thought of as a
  subbundle of~$T(\mathcal{E})$.  Given an inclusion
  $\mathcal{E}\subset\J^{\infty}(\pi)$, $\pi\colon E\to M$, as above, the
  fibre bundle $\tau\colon\T(\mathcal{E})\to\mathcal{E}$ can be identified
  with the vertical bundle with respect to the projection $\mathcal{E}\to M$.

  Every Cartan vector field $X\in\C\X(\mathcal{E})$ can be lifted to a vector
  field $\tilde{X}\in\X(\T(\mathcal{E}))$ as follows.  It suffices to define an
  action of~$\tilde{X}$ on fibre-wise linear functions on~$\T(\mathcal{E})$
  that can be naturally identified with Cartan $1$-forms
  $\omega\in\Lac^1(\mathcal{E})$.  We put
  \begin{equation*}
    \tilde{X}(\omega)=L_X(\omega),
  \end{equation*}
  where $L_X$ denotes the Lie derivative.  In coordinates, we have
  $\tilde{D_i}=D_i$, where $D_i$ in the right-hand side are the total
  derivatives on~$\T(\mathcal{E})$.
  
  Let us describe an inclusion of~$\T(\mathcal{E})$ to a jet space.  Assume
  that~$\mathcal{E}\subset\J^{\infty}(\pi)$ is given by~$F=0$.  Let
  $\chi\colon\J^{\infty}(\pi\times_M\pi)\to\J^{\infty}(\pi)$ be defined by the
  projection to the first factor.  Then
  $\T(\mathcal{E})\subset\J^{\infty}(\pi\times_M\pi)$ is defined by equations
  \begin{equation*}
    \tilde{F}=0,\qquad\tilde{\ell}_F(\bi{v})=0,
  \end{equation*}
  where the tilde over~$F$ denotes the pullback by~$\chi$ (with the tilde for
  the Cartan vector fields defined above),
  $\bi{v}\in\tilde{\kappa}(\pi)=\F(\pi\times_M\pi,\pi)$ is the projection to
  the second factor $\Gamma(\pi\times_M\pi)\to\Gamma(\pi)$.

  So, $\T(\mathcal{E})$ is an equation and the vector fields of the
  form~$\tilde{X}$ generate the Cartan distribution on it.  Thus, the bundle
  $\tau\colon\T(\mathcal{E})\to\mathcal{E}$ is a covering called the
  \emph{tangent covering} to~$\mathcal{E}$.

  In coordinates, we have $\bi{v}=(v^1,\dots,v^m)$ if the coordinates on
  $\J^{\infty}(\pi\times_M\pi)$ are $x^i$,~$u^j_I$,~$v^j_I$, with $u^j_I$ and
  $v^j_I$ corresponding to the first and second factors, respectively.  Thus,
  a coordinate description of~$\T(\mathcal{E})$ has the form
  \begin{equation*}
    F_s(x^i,u^j_I)=0,\qquad
    \sum_{\alpha,I}\frac{\partial F_j}{\partial u_I^{\alpha}}v_I^{\alpha}=0.
  \end{equation*}
  
  Note that if $\mathcal{E}$ is a normal equation then~$\T(\mathcal{E})$ is
  normal as well.
\end{example}

The reader will find more examples and a detailed discussion of coverings in
Section~\ref{sec:gjsis_nonloc:nonlocal-theory} below.

\subsection{Linearization}
\label{sec:gjsis_eqs:line-diff-eqaut}

Let $\mathcal{E}\subset\J^{\infty}(\pi)$ be an equation defined by a section
$F\in P$.  Denote by~$\kappa$ the restriction of the module
$\kappa(\pi)=\F(\pi,\pi)$ to~$\mathcal{E}$.  The \emph{linearization} of the
equation~$\mathcal{E}$ is the restriction to~$\mathcal{E}$ of the
linearization of~$F$:
\begin{equation*}
  \bar{\ell}_F=\eval{\ell_F}_{\mathcal{E}}\colon\kappa\to P.
\end{equation*}
We denote by bar the restriction of a $\C$-differential operator
to~$\mathcal{E}$ and preserve the notation of modules for their restrictions.

\begin{remark}
  The operator~$\ell_F$ is well-defined globally only if the module~$P$ has
  the form $P=\F(\pi,\pi')$.  For an arbitrary module~$P$ the operator
  $\ell_F$ is defined only locally.  But its restriction $\bar{\ell}_F$ is
  well-defined globally on the whole~$\mathcal{E}$.
\end{remark}

\begin{remark}
  \label{rem:gjsis_eqs:1}
  A $\C$-differential operator~$\Delta$ on~$\mathcal{E}$ is called
  \emph{normal} if for any $\C$-differential operator~$\square$ the condition
  $\square\circ\Delta=0$ implies $\square=0$.

  The linearization operator~$\bar{\ell}_F$ of an equation~$\mathcal{E}$ is
  normal if and only if the section~$F$ is normal. If an
  equation~$\mathcal{E}$ admits a normal representation, i.e., an
  embedding~$\mathcal{E}=\{F=0\}\subset\J^\infty(\pi)$ such that the
  corresponding operator~$\bar{\ell}_F$ is normal then there always exists
  another representation~$\mathcal{E}=\{F'=0\}\subset\J^\infty(\pi')$ for
  which~$\mathcal{E}$ acquires an evolutionary form.
\end{remark}

Recall that in coordinates the linearization~$\bar{\ell}_F$ has the form
(cf.~\eqref{eq:gjsis_jets:20}):
\begin{equation*}
  \bar{\ell}_F=
  \left\Vert\sum_I\frac{\partial F_j}{\partial u_I^\alpha}D_I
  \right\Vert_{\alpha=1,\dots,m.}^{j=1,\dots,l}
\end{equation*}

If an equation~$\mathcal{E}$ is defined by two different sections $F_1\in P_1$
an $F_2\in P_2$ in different, generally speaking, jet spaces
  \begin{equation*}
   \xymatrixcolsep{1pc}
    \xymatrixrowsep{0.5pc}
    \xymatrix{
      &J^{\infty}(\pi_1) \\
      \mathcal{E}\ar[ur]\ar[dr] \\
      &J^{\infty}(\pi_2)
    }
  \end{equation*}
  then the corresponding linearizations
  $\bar{\ell}_{F_1}\colon\kappa_1\to P_1$ and
  $\bar{\ell}_{F_2}\colon\kappa_2\to P_2$ are \emph{equivalent}
  \cite{DudnikovSamborski:LOvSPDE,IgoninKerstenKrasilshchikVerbovetskyVitolo:VBGPDE}
  in the sense that there exist $\C$-differential operators $\alpha$,
  $\beta$, $\alpha'$, $\beta'$, $s_1$, and $s_2$ on~$\mathcal{E}$
  \begin{equation}
    \label{eq:gjsis_eqs:3}
  \xymatrixcolsep{5pc}
  \xymatrix{
    \kappa_1\ar[r]_{\bar{\ell}_{F_1}}\ar@<.5ex>[d]^{\alpha}
    &P_1\ar@<.5ex>[d]^{\alpha'}\ar@/_1pc/@<-1ex>[l]_{s_1} \\
    \kappa_2\ar@<.5ex>[u]^{\beta}\ar[r]^{\bar{\ell}_{F_2}}
    &P_2\ar@<.5ex>[u]^{\beta'}\ar@/^1pc/@<1ex>[l]^{s_2}
  }
\end{equation}
such that
\begin{equation*}
  \bar{\ell}_{F_1}\beta=\beta'\bar{\ell}_{F_2},\quad
  \bar{\ell}_{F_2}\alpha=\alpha'\bar{\ell}_{F_1},\quad
  \beta\alpha=\id+s_1\bar{\ell}_{F_1},\quad
  \alpha\beta=\id+s_2\bar{\ell}_{F_2}.
\end{equation*}

\begin{example}
  Consider two presentations~\eqref{eq:gjsis_eqs:6} and~\eqref{eq:gjsis_eqs:5}
  of the KdV equation from Example~\ref{exmp:gjsis_eqs:1}.  The operators of
  diagram~\eqref{eq:gjsis_eqs:3} are:
    \begin{align*}
    \bar{\ell}_{F_1}&=D_t-D_x^3-6uD_x-6u_x, \\
    \bar{\ell}_{F_2}&=\begin{pmatrix}
      D_x & -1 & 0 \\
      0 & D_x & -1 \\
      -D_t+6v & 6u & D_x
    \end{pmatrix}
    \\
    \intertext{and} %\\
    \alpha&=\begin{pmatrix}1 \\ D_x \\ D_{xx}\end{pmatrix}, \\
    \beta&=\begin{pmatrix}1 & 0 & 0\end{pmatrix}, \\
    \alpha'&=\begin{pmatrix}\hphantom{-}0 \cr \hphantom{-}0 \cr -1
    \end{pmatrix}, \\
    \beta'&=\begin{pmatrix}-D_{xx}-6u & -D_x & -1\end{pmatrix}, \\
    s_1&=0, \\
    s_2&=\begin{pmatrix}
      0 & 0 & 0 \\
      1 & 0 & 0 \\
      D_x & 1 & 0
    \end{pmatrix}.
  \end{align*}

  The form of operators $\alpha$ and $\beta$ is obvious from the form of
  operators $a$ and $b$ in Example~\ref{exmp:gjsis_eqs:1}.  The operators
  $\alpha'$ and $\beta'$ show how equations~\eqref{eq:gjsis_eqs:6} and
  \eqref{eq:gjsis_eqs:5} are obtained one from the other.
\end{example}

\begin{example}
  Let the number of dependent variables $m$ be equal to~$1$.  Consider the
  lift $L\colon\J^{\infty}(\pi)\to\J^{\infty}(\pi)$ of the Legendre
  transformation
  \begin{equation*}
    \textstyle
    \eval{L}_{\J^1(\pi)}(x^i,u,u_{x^i})=(u_{x^i},
    \sum_{\alpha}x^{\alpha}u_{x^{\alpha}}-u,x^i).
  \end{equation*}
  It is defined wherever $\det\left\Vert u_{x^ix^j}\right\Vert\ne0$ and
  preserves the Cartan distribution.  Of course, $L$ is not a map of fibre
  bundles.

  Consider an equation~$\mathcal{E}$ defined by $F=0$, then $L^*(F)=0$ defines
  the same equation~$\mathcal{E}$. The operators in
  diagram~\eqref{eq:gjsis_eqs:3} are as follows:
  \begin{align*}
    F_1&=F, \\
    F_2&=L^*(F) \\
    \intertext{and} 
    \alpha&=-1,\qquad \alpha'=1,\qquad s_1=0, \\
    \beta&=-1,\qquad \beta'=1,\qquad s_2=0.
  \end{align*}
  Let us explain how to compute the maps~$\alpha$ and~$\beta$.  To this end,
  we have to take a symmetry of the Cartan distribution
  $X\in\Xc(\pi)/\C\X(\pi)$ and see how the generating section of~$X$
  transforms under the Legendre map.  The generating section can be computed
  by the formula $\phi=\omega(X)$, where $\omega=\rmd u-\sum_i u_{x^i}\rmd
  x^i$ is a Cartan form.  So,
  $\alpha(\phi)=\omega(L(X))=L^*(\omega)(X)=-\omega(X)=-\phi$. The same holds
  for~$\beta$.
\end{example}

\subsection{Symmetries and recursions}
We have defined symmetries of a differential equation~$\mathcal{E}$ as
elements of quotient space
\begin{equation*}
    \sym\mathcal{E}=\Xc(\mathcal{E})/\C\X(\mathcal{E}),
\end{equation*}
where $\Xc(\mathcal{E}) =\sd{X\in\X(\mathcal{E})}{[X,\C\X(\mathcal{E})]
  \subset\C\X(\mathcal{E})}$, that is, symmetries of equations are symmetries
of the Cartan distribution on it modulo trivial symmetries (the ones belonging
to the Cartan distribution).

Obviously, symmetries of an equation form a Lie algebra with respect to the
commutator.

Given an inclusion $\mathcal{E}\subset\J^{\infty}(\pi)$, each
symmetry~$X\in\sym\mathcal{E}$ contains exactly one vector
filed~$X^v\in\Xc(\mathcal{E})$ vertical with respect to the projection
$\pi_{\infty}\colon\mathcal{E}\to M$, where $M$ is the base manifold of the
bundle~$\pi$.  Next, for every such a field~$X^v$ there always exists an
evolutionary vector field~$\Ev_{\phi'}$ such that
\begin{equation*}
  X^v=\eval{\Ev_{\phi'}}_{\mathcal{E}}.
\end{equation*}
For any such a field, the
restriction~$\phi=\eval{\phi'}_{\mathcal{E}}\in\kappa=\kappa\eval{\kappa(\pi)}_{\mathcal{E}}$
depends on~$X^v$ only. Hence, we can denote the field $X^v$ by~$\Ev_\phi$.

The element~$\phi\in\kappa$ such that~$X^v=\Ev_{\phi}$ is called the
\emph{generating section} (or a \emph{characteristic}) of this symmetry.

If the equation at hand~$\mathcal{E}$ is defined by an equality $F=0$, then
the existence of a symmetry~$\Ev_\phi$ boils down to the condition
\begin{equation*}
  \eval{\Ev_{\phi'}(F)}_{\mathcal{E}}=0,
\end{equation*}
where, as above, $\phi'\in\kappa(\pi)$ is an arbitrary extension of~$\phi$,
which is equivalent to the condition
\begin{equation}
  \label{eq:gjsis_eqs:4}
  \bar{\ell}_F(\phi)=0\quad\text{on~$\mathcal{E}$.}
\end{equation}

This is the determining equation for the symmetries of the
equation~$\mathcal{E}=\{F=0\}$.

The Jacobi bracket~\eqref{eq:gjsis_jets:15} yields a Lie algebra structure
on~$\sym\mathcal{E}$ in terms of generating sections:
\begin{equation}
  \label{eq:gjsis_eqs:8}
  \{\phi,\psi\}=\bar{\ell}_\psi(\phi)-\bar{\ell}_\phi(\psi),\qquad
  \phi,\psi\in\ker\bar{\ell}_F\subset\kappa.
\end{equation}

Searching symmetries, that is solving~\eqref{eq:gjsis_eqs:4}, one usually
begins with choosing \emph{internal coordinates} on~$\mathcal{E}$.
\begin{example}
  We have seen in Example~\ref{exmp:gjsis_eqs:1} that on the KdV
  equation
    \begin{equation*}
    u_t-6uu_x-u_{xxx}=0
  \end{equation*}
  the functions $x$,~$t$,~$u_s=u_{x\dots x}$ ($s$ times) can be taken
  to be coordinates on~$\mathcal{E}$.  These are internal coordinates
  for the KdV equation.  
\end{example}

Of course, the choice of internal coordinates is not unique.

Next, one can start with finding symmetries whose generating sections
depend on derivatives of order less than some number~$k$,
$\phi=\phi(x,t,u,u_1,\dots,u_k)$.

\begin{example}
  \label{exmp:gjsis_eqs:2}
  Let us find symmetries of the KdV equation such that
  $\phi=\phi(x,t,u,u_1)$.

  The determining equation for symmetries has the form:
  \begin{equation}
    \label{eq:gjsis_eqs:7}
    (D_t-D_x^3-6uD_x-6u_x)\phi(x,t,u,u_1)=0.
  \end{equation}
  The left-hand side of this equation is a polynomial in $u_3$
  and~$u_2$.  The coefficient of the product $u_2u_3$ is equal
  to~$-3\phi_{u_1u_1}$, so $\phi_{u_1u_1}=0$, and we have
  \begin{equation*}
    \phi=\phi^0(x,t,u)+\phi^1(x,t,u)u_1.
  \end{equation*}
  Now, the left-hand side of~\eqref{eq:gjsis_eqs:7} is a polynomial in
  $u_3$,~$u_2$ and~$u_1$, with coefficient of~$u_3$ equal to
  $-3D_x(\phi^1)$.  Hence,
  \begin{equation*}
    \phi=\phi^0(x,t,u)+\phi^1(t)u_1.
  \end{equation*}
  With such a~$\phi$, the left-hand side of~\eqref{eq:gjsis_eqs:7} does not
  depend on~$u_3$ any more.  The coefficient of the product~$u_1u_2$ is
  $-3\phi^0_{uu}$, hence
  \begin{equation*}
    \phi=\phi^{00}(x,t)+u\phi^{01}(t,x)+\phi^1(t)u_1.
  \end{equation*}
  The coefficient of~$u_2$ is $-3\phi^{01}_x$, so that
  \begin{equation*}
    \phi=\phi^{00}(x,t)+u\phi^{01}(t)+\phi^1(t)u_1
  \end{equation*}
  and the left-hand side of~\eqref{eq:gjsis_eqs:7} is a polynomial in
  $u_1$ and~$u$.  The coefficient of the product~$uu_1$ shows that
  $\phi^{01}=0$.  The remaining coefficients of~$u_1$ and~$u$, and the
  free term reveal that
  \begin{equation*}
    \phi^{00}=c_0,\qquad \phi^1=6\phi^{00}t+c_1,
  \end{equation*}
  where $c_0$ and $c_1$ are arbitrary constants.  Therefore we have
  found two independent symmetries of the KdV equation
  \begin{equation}
    \label{eq:gjsis_eqs:9}
    \phi_1=u_x\quad\text{and}\quad\phi_2=6tu_x+1.
  \end{equation}
\end{example}

Let $X\in\sym\mathcal{E}=\Xc(\mathcal{E})/\C\X(\mathcal{E})$ be a
symmetry.  Vector fields $Y\in\Xc(\mathcal{E})$ belonging the
equivalence class~$X$ we shall call \emph{representatives} of~$X$.  As
we explained above, any symmetry has a unique representative of the
form~$\Ev_\phi$.  But this representative can be not the simplest one.

\begin{example}
  Symmetries~\eqref{eq:gjsis_eqs:9} of the KdV equation can be
  represented, respectively, by the fields $Y_1=-\partial/\partial x$
  and $Y_2=\partial/\partial u-6t\partial/\partial x$, because
  \begin{equation*}
    \Ev_{u_x}-Y_1=D_x,\qquad
    \Ev_{6tu_x+1}-Y_2=6tD_x.
  \end{equation*}
  The fields $Y_1$ and $Y_2$ are the lifts from the zero order jet
  space, hence they correspond to one-parameter groups of
  transformations:
  \begin{align*}
    Y_1\colon & x'=x-\epsilon &\text{(translation along~$x$)}, \\
    Y_2\colon & x'=x-6\epsilon t,\quad u'=u+\epsilon
    &\text{(Galilean symmetry)}.
  \end{align*}
\end{example}

A symmetry~$X$ is called \emph{classical} if it can be represented by
a field~$\eval{Y}_{\mathcal{E}}$, with $Y\in\X(\pi)$ being the lift
from a finite order jet space.  Classical symmetries form a subalgebra
of the Lie algebra $\sym\mathcal{E}$.  By the Lie-B\"{a}cklund theorem
(see~~\cite{KrasilshchikVinogradov:SCLDEqMP}), $Y$ is a lift from the
zero order jet space (if the number of dependent variables $m>1$) or
from the first order jet space (if the number of dependent variables
$m=1$).  Thus, in coordinate language, the generating section of a
classical symmetry has the form
\begin{equation*}
  \label{eq:gjsis_eqs:10}
  \phi=\begin{cases}(\phi_1,\dots,\phi_m) &\text{for $m>1$}, \\
    \phi(x^i,u,u_i) &\text{for $m=1$},
  \end{cases}
\end{equation*}
where $\phi_\alpha=b_{\alpha}(x^i,u^j)+\sum_{k=1}^na_k(x^i,u^j)u^{\alpha}_k$,
while $\phi$ is an arbitrary smooth function.  The vector field $Y$
on~$\J^0(\pi)$ or~$\J^1(\pi)$ that represents the symmetry with generating
section~$\phi$ has the form
\begin{equation*}
  Y=\begin{cases}\textstyle\sum_{j=1}^mb_j\frac{\partial}{\partial u^j}
    -\sum_{k=1}^na_k\frac{\partial}{\partial x^k} &\text{for $m>1$}, \\
    -\textstyle\sum_{i=1}^n\frac{\partial\phi}{\partial
      u_i}\frac{\partial}{\partial x^i}+
    (\phi-\sum_{i=1}^nu_i\frac{\partial\phi}{\partial u_i})
    \frac{\partial}{\partial u}
    +\sum_{i=1}^n(\frac{\partial\phi}{\partial x^i}+
    u_i\frac{\partial\phi}{\partial u})\frac{\partial}{\partial u_i}
    &\text{for $m=1$}.
  \end{cases}
\end{equation*}

\begin{example}
  In Example~\ref{exmp:gjsis_eqs:2} above we computed symmetries of
  the KdV equation with generating sections depending on
  $x$,~$t$,~$u$, and~$u_x$.  To find all classical symmetries, we
  allow also dependence on~$u_t$ and compute symmetries in the same
  manner as in Example~\ref{exmp:gjsis_eqs:2}.  We get two additional
  classical symmetries:
  \begin{equation*}
    \phi_3=u_t\quad\text{and}\quad\phi_4=xu_x+3tu_t+2u,
  \end{equation*}
  with the corresponding one-parameter groups of transformations being
  \begin{align*}
    \phi_3\colon & t'=t-\epsilon &&\text{(translation along~$t$)}, \\
    \phi_4\colon & x'=\rme^{-\epsilon}x,\quad t'=\rme^{-3\epsilon}t, \quad
    u'=\rme^{2\epsilon}u &&\text{(scale symmetry).}
  \end{align*}
\end{example}

Solving~\eqref{eq:gjsis_eqs:4} for sections~$\phi$ that depend on at most
$k$th order derivatives will not give us a complete description of all
symmetries.  We can find all classical symmetries or a slightly larger
subspace of symmetries, which can be considered to be a \emph{lower estimate}
of the full symmetry algebra.  Letting the maximal order~$k$ be arbitrary and
by solving~\eqref{eq:gjsis_eqs:4} describe the dependence of~$\phi$ on
derivatives of order $k$,~$k-1$, etc.  This will be an \emph{upper estimate}
of the symmetry algebra.  Sometimes these estimates allow to find more
symmetries and, in some cases, obtain a complete description of the symmetry
algebra.  We refer the reader to \cite{KrasilshchikVinogradov:SCLDEqMP} for
examples of such calculations.

A different approach is to look for a \emph{recursion operator} that is a
$\C$-differential operator $R\colon\kappa\to\kappa$ such that there exists
another $\C$-differential operator~$R'$ satisfying the condition
\begin{equation}
  \label{eq:gjsis_eqs:11}
  \bar{\ell}_FR=R'\bar{\ell}_F.
\end{equation}
Operators of the form $R=\square\bar{\ell}_F$, with $\square$ being
arbitrary, enjoy~\eqref{eq:gjsis_eqs:11} for all equations, so we
consider recursion operators modulo such trivial ones.  Obviously,
$R(\sym\mathcal{E})\subset\sym\mathcal{E}$, so that having
a recursion operator we can produce an infinite number of symmetries
from a given one.

\begin{example}\label{sec:gjsis_eqs:ex_Heat-eq}
  The heat equation $u_t=u_{xx}$ has two recursion operators of the
  first order
  \begin{equation*}
    R_1=D_x\quad\text{and}\quad R_2=2tD_x+x
  \end{equation*}
  (and, of course, the identity operator, which is of no interest).
\end{example}

For nonlinear equations we often are able to find a nontrivial
recursion operator only if we allow it to contain the ``integration''
operator~$D_x^{-1}$.  In Section~\ref{sec:gjsis_nonloc:zero-curv-repr}
we explain how to define such recursion operators in a rigorous way.

\begin{example}
  As it was already mentioned, the KdV equation has the Lenard
  recursion operator
  \begin{equation}
    \label{eq:gjsis_eqs:12}
    R=D_x^2+4u+2u_xD_x^{-1}.
  \end{equation}
  It is obvious that $R(u_x)=u_t$ and
  $R^2(u_x)=R(u_t)=u_{xxxxx}+10uu_{xxx}+20u_xu_{xx}+30u^2u_x$.  Thus, we have
  a new \emph{higher} (that is non-classical) symmetry of KdV.  We can proceed
  in the same way and compute $R^3(u_x)$,~$R^4(u_x)$ and so on.  As it was
  already mentioned (see Section~\ref{sec:gjsis_jets:hamilt-form}), all
  $R^k(u_x)$ exist, that is the computation of $D_x^{-1}$ will be possible for
  all~$k$.  So, we have constructed an infinite set of symmetries of the KdV
  equation.

  If we try to apply the Lenard recursion operator to the other two
  classical symmetries, $\phi_2$ (Galilean) and $\phi_4$ (scale), we
  get $R(\phi_2)=2\phi_4$, but $R(\phi_4)$ does not exist.  In fact,
  $R(\phi_4)$ is a \emph{nonlocal} symmetry, as explained below in
  Section~\ref{sec:gjsis_nonloc:nonlocal-theory}.
\end{example}

Now, let us explain how to compute recursion operators.  To this end,
note that $\C$-differential operators $\kappa\to\kappa$ modulo
operators of the form $\square\bar{\ell}_F$ can be naturally
identified with elements of
$\Lac^1(\mathcal{E})\otimes_{\F(\mathcal{E})}\kappa$, that is with
$\kappa$-valued Cartan $1$-forms.  This identification takes an
operator $R\colon\kappa\to\kappa$ to the form
$\omega_R\in\Lac^1(\mathcal{E})\otimes_{\F(\mathcal{E})}\kappa$, such
that $\omega_R(\Ev_{\phi})=R(\phi)$.

Next, recall that the Cartan $1$-forms $\omega\in\Lac^1(\mathcal{E})$
are functions on the tangent covering $\T(\mathcal{E})$ (see
Example~\ref{exmp:gjsis_eqs:3}) that are linear along the fibres of
the projection $\tau\colon\T(\mathcal{E})\to\mathcal{E}$.
Hence the forms
$\omega\in\Lac^1(\mathcal{E})\otimes_{\F(\mathcal{E})}\kappa$ are
elements of the pullback~$\tilde{\kappa}$ of the module~$\kappa$
by~$\tau$.  Thus, to a recursion operator~$R\colon\kappa\to\kappa$ we
assign a fibre-wise linear element $\omega_R\in\tilde{\kappa}$.

In coordinates, to the operator
$R=\left\Vert\sum_If_I^{\alpha\beta}D_I\right\Vert$ there corresponds
the element
\begin{equation*}
  \omega_R=\Bigl(\sum_{\beta,I}f_I^{1\beta}v_I^{\beta},\dots,
  \sum_{\beta,I}f_I^{m\beta}v_I^{\beta}\Bigr)
\end{equation*}
with $v_I^\beta$ being coordinates along the fibres of~$\tau$.

Condition~\eqref{eq:gjsis_eqs:11} on $R$ is equivalent to the
following condition on~$\omega_R$:
\begin{equation}
  \label{eq:gjsis_eqs:13}
  \eval{\tilde{\ell}_F(\omega_R)}_{\T(\mathcal{E})}=0.
\end{equation}

\begin{example}
  For the heat equation $u_t=u_{xx}$ the tangent covering is given by
  \begin{equation*}
    u_t=u_{xx},\qquad v_t=v_{xx}.
  \end{equation*}
  Let $x$, $t$, $u_k=u_{x\dots x}$, and $v_k=v_{x\dots x}$ be internal
  coordinates on it.  To compute recursion operators of order~$1$
  in~$D_x$ we are to solve the equation
  \begin{equation*}
    (D_t-D_x^2)(f(x,t,u_k)v_x+g(x,t,u_k)v)=0.
  \end{equation*}
  It is easy to show that the solutions are: $v$, $v_x$, and $2tv_x+xv$, with
  the corresponding operators $\id$, $D_x$, and $2tD_x+x$, as it was indicated
  in Example~\ref{sec:gjsis_eqs:ex_Heat-eq}.
\end{example}

\begin{example}
  The tangent covering over the KdV equation has the form
  \begin{equation*}
    u_t-u_{xxx}-6uu_x=0,\qquad v_t-v_{xxx}-6uv_x-6u_xv=0.
  \end{equation*}
  Computation of recursion operators amounts to solving equations like the
  following:
  \begin{equation}
    \label{eq:gjsis_eqs:16}
    (D_t-D_x^3-6uD_x-6u_x)(f_2v_{xx}+f_1v_x+f_0v+f_{-1}v_{-1}),
  \end{equation}
  where $f_i=f_i(x,t,u_k)$ and $v_{-1}$ is a new variable such that
  \begin{align}
    D_x(v_{-1})&=v\label{eq:gjsis_eqs:14}, \\
    D_t(v_{-1})&=v_{xx}+6uv.\label{eq:gjsis_eqs:15}
  \end{align}
  As we saw above, the variables~$v_I$ correspond to the total
  derivatives~$D_I$, so that the variable~$v_{-1}$ is defined
  in~\eqref{eq:gjsis_eqs:14} to correspond to~$D_x^{-1}$,
  while~\eqref{eq:gjsis_eqs:15} provides the equality
  $D_t(D_x(v_{-1}))=D_x(D_t(v_{-1}))$.  We shall defer a rigorous explanation
  of this computation until Section~\ref{sec:gjsis_nonloc:nonlocal-theory}.

  One can check that $v_{xx}+4uv+2u_xv_{-1}$ is a solution
  of~\eqref{eq:gjsis_eqs:16}, it yields the Lenard recursion
  operator~\eqref{eq:gjsis_eqs:12}.
\end{example}

The Lie algebra structure on $\ker\bar{\ell}_F\subset\kappa$ (the Jacobi
bracket on symmetries) has a natural extension to a Lie superalgebra on
$\ker\eval{\tilde{\ell}_F}_{\Lac^*(\mathcal{E})
  \otimes_{\F(\mathcal{E})}\kappa}$, called the \emph{Nijenhuis} (also
\emph{Fr\"olicher-Nijenhuis\textup{)} bracket} and denoted by
$\lshad\cdot\,,\cdot\rshad$.  For a detailed definition we refer the reader
to~\cite{KrasilshchikKersten:SROpCSDE}.  Since we identified recursion
operators with elements of
$\ker\eval{\tilde{\ell}_F}_{\Lac^1(\mathcal{E})\otimes_{\F(\mathcal{E})}\kappa}$,
we can compute the Nijenhuis bracket of them.  Two particular cases are of
importance for us here: the bracket of a recursion operator with itself (the
Nijenhuis torsion):
\begin{align*}
  \frac12\lshad R,R,\rshad(\phi_1,\phi_2)%\\
  &=\{R(\phi_1),R(\phi_2)\}-R\{R(\phi_1),\phi_2\}-R\{\phi_1,R(\phi_2)\}
  +R^2\{\phi_1,\phi_2\} \\
  &=\ell_{R,\phi_2}(R(\phi_1))-\ell_{R,\phi_1}(R(\phi_2))
  +R\ell_{R,\phi_1}(\phi_2)-R\ell_{R,\phi_2}(\phi_2),
\end{align*}
where $\phi_1$,~$\phi_2\in\kappa$, and the bracket of a recursion
operator and a symmetry with generating section
$\phi\in\ker\bar{\ell}_F\subset\kappa$ (the Lie derivative):
\begin{equation*}
  L_{\phi}(R)(\phi')=\lshad \phi,R\rshad(\phi')
  =\{\phi,R(\phi')\}-R\{\phi,\phi'\},\qquad\phi'\in\kappa,
\end{equation*}
or
\begin{equation*}
  L_{\phi}(R)=\Ev_{\phi}(R)-[\ell_{\phi},R].
\end{equation*}

A recursion operator~$R$ is called \emph{Nijenhuis} (or
\emph{hereditary}) if its Nijenhuis torsion vanishes, i.e., $\lshad
R,R\rshad=0$.  Almost all known recursion operators are Nijenhuis,
including the operators encountered above.  The main property of
Nijenhuis operators is the following: for every two symmetries with
generating sections $\phi_1$,~$\phi_2\in\kappa$ such that
$\{\phi_1,\phi_2\}=0$, $L_{\phi_1}(R)=0$, $L_{\phi_2}(R)=0$ and for
arbitrary $k_1$ and $k_2$ we have
$\{R^{k_1}(\phi_1),R^{k_2}(\phi_2)\}=0$.

\begin{example}
  For the KdV equation take $\phi_1=\phi_2=u_x$.  Obviously, we have
  $L_{u_x}(R)=0$, where $R$ is the Lenard recursion
  operator~\eqref{eq:gjsis_eqs:12}, so that all symmetries $R^k(u_x)$ commute.
\end{example}

\subsection{Conservation laws}

A \emph{conserved current}~$\omega$ on an equation~$\mathcal{E}$ with $n$
independent variables is a closed horizontal $(n-1)$-form on~$\mathcal{E}$,
i.e., a form~$\omega\in\Lah^{n-1}(\mathcal{E})$ such that $\rmdh\omega=0$.

\begin{example}
  \label{exmp:gjsis_eqs:4}
  Consider the equation of continuity in fluid dynamics
  \begin{equation*}
    \rho_t+(\rho v^1)_{x_1}+(\rho v^2)_{x_2}+(\rho v^3)_{x_3}=0.
  \end{equation*}
  The form $\omega=\rho\rmd x_1\wedge\rmd x_2\wedge\rmd x_3 -\rho
  v^1\rmd t\wedge\rmd x_2\wedge\rmd x_3 +\rho v^2\rmd t\wedge\rmd
  x_1\wedge\rmd x_3 -\rho v^3\rmd t\wedge\rmd x_1\wedge\rmd x_2$ is a
  conserved current.
\end{example}

\begin{example}
  \label{exmp:gjsis_eqs:5}
  For the KdV equation $u_t-u_{xxx}-6uu_x=0$ the forms
  \begin{align*}
    &u\rmd x+(u_{xx}+3u^2)\rmd t, \\
    &u^2\rmd x+(2uu_{xx}-u_x^2+4u^3)\rmd t, \\
    &(u_x^2/2-u^3)\rmd x+(u_xu_{xxx}-u_{xx}^2/2-3u^2u_{xx}
    +6uu_x^2-9u^4/2)\rmd t
  \end{align*}
  are conserved currents.
\end{example}

For any $\eta\in\Lah^{n-2}(\mathcal{E})$ the form $\omega=\rmdh\eta$
is always a conserved current.  Such currents are called
\emph{trivial} because they are not related to the equation properties
and hence are of no interest.  Consider the quotient space
\begin{equation*}
  H_h^{n-1}(\mathcal{E})=\sd{\omega\in\Lah^{n-1}(\mathcal{E})}{\rmdh\omega=0}
  \big/\sd{\omega\in\Lah^{n-1}(\mathcal{E})}{\omega=\rmdh\eta,
    \quad\eta\in\Lah^{n-2}(\mathcal{E})}
\end{equation*}
of \emph{horizontal cohomology} of~$\mathcal{E}$.  We also want to quotient
out the \emph{topological} conserved currents that lie in the image of the map
$\zeta\colon H^{n-1}(\mathcal{E})\to H_h^{n-1}(\mathcal{E})$ induced by the
natural projection $\La^{n-1}(\mathcal{E})\to\Lah^{n-1}(\mathcal{E})$; here
$H^{n-1}(\mathcal{E})$ is the $(n-1)$st group of the de~Rham cohomology of the
space~$\mathcal{E}$.  Such currents are related to the topology of the
equation~$\mathcal{E}$ only.  Thus, we define
$\cl(\mathcal{E})=H_h^{n-1}(\mathcal{E})/\im\zeta$ to be the set of
\emph{conservation laws} of equation~$\mathcal{E}$.

We shall now discuss how to compute conservation laws for
\emph{normal} equations.  To this end, let us consider the following
complex:
\begin{equation}
  \label{eq:gjsis_eqs:17}
  0\xrightarrow{}\Omega^0(\mathcal{E})
  \xrightarrow{\delta}\Omega^1(\mathcal{E})
  \xrightarrow{\delta}\Omega^2(\mathcal{E})
  \xrightarrow{\delta}\dots,
\end{equation}
where
\begin{equation}
  \label{eq:gjsis_eqs:20}
  \Omega^0(\mathcal{E})=\cl(\mathcal{E}),\qquad
    \Omega^p(\mathcal{E})
    =Z^p\big/\rmdh(\Lac^p(\mathcal{E})\otimes_{\F(\mathcal{E})}
    \Lah^{n-2}(\mathcal{E})),
\end{equation}
and $Z^p=\ker\rmdh\subset\Lac^p(\mathcal{E})\otimes_{\F(\mathcal{E})}
\Lah^{n-1}(\mathcal{E})$, for $p>0$.  The differential~$\delta$ is induced by
the Cartan differential~$\rmdc$, so that $\delta^2=0$.
\begin{remark}
  Complex~\eqref{eq:gjsis_eqs:17} is a part of Vinogradov's
  $\C$-spectral sequence (see \cite{Vinogradov:AlGFLFT,Vinogradov:SSAsNDEqAlGFLFTC,Vinogradov:SSLFCLLTNT,Vinogradov:CAnPDEqSC}), namely
  $\Omega^p(\mathcal{E}) =E_1^{p,n-1}(\mathcal{E})$ for $p>0$,
  $\Omega^0(\mathcal{E})
  =E_1^{0,n-1}(\mathcal{E})\big/H^{n-1}(\mathcal{E})$, and
  $\delta=d_1^{*,n-1}$.
\end{remark}

Assume that $\mathcal{E}$ is a normal equation defined by a normal
section $F\in P$.  In this case complex~\eqref{eq:gjsis_eqs:17} can be
described in the following way:
\begin{equation}
  \label{eq:gjsis_eqs:23}
  \Omega^p(\mathcal{E})=\Theta^p\big/\Theta^p_{\ell},
\end{equation}
where $\Theta^p$ is a subset of $\CDiffsk_{p-1}(\kappa,\hat{P})$ that
consists of operators $\Delta\in\CDiffsk_{p-1}(\kappa,\hat{P})$ such that
\begin{equation*}
  \bar{\ell}_F^*\Delta(\phi_1,\dots,\phi_{p-1})
  -\sum_{\alpha=1}^{p-1}\Delta^{*_{\alpha}}
  (\phi_1,\dots,\bar{\ell}_F(\phi_{\alpha}),\dots,\phi_{p-1})=0
\end{equation*}
for all $\phi_1,\dots,\phi_{p-1}\in\kappa$, where~${^{*_{\alpha}}}$ denotes
the operation of taking the adjoint with respect to the $\alpha$th argument.
The subset $\Theta^p_{\ell}\subset\Theta^p$ consists of operators
$\Delta\in\Theta^p$ of the form
\begin{equation}
  \label{eq:gjsis_eqs:22}
  \Delta(\phi_1,\dots,\phi_{p-1})=\sum_{\alpha=1}^{p-1}(-1)^{\alpha+1}\Delta'
  (\bar{\ell}_F(\phi_{\alpha}),\phi_1,\dots,
  \hat{\phi}_{\alpha},\dots,\phi_{p-1})
\end{equation}
for some $\C$-differential operators $\Delta'\colon
P\times\kappa\times\dots\times\kappa\to\hat{P}$.

In particular, for $p=1$, we have
\begin{equation}
  \label{eq:gjsis_eqs:18}
  \Omega^1(\mathcal{E})=\ker\bar{\ell}_F^*\subset\hat{P}.
\end{equation}
If $p=2$ then
\begin{equation*}
  \Omega^2(\mathcal{E})
  =\sd{\Delta\in\CDiff(\kappa,\hat{P})}
  {\bar{\ell}_F^*\Delta=\Delta^*\bar{\ell}_F}\big/
  \sd{\Delta'\bar{\ell}_F}
  {\Delta'\in\CDiff(P,\hat{P}),\ {\Delta'}^*=\Delta'}.
\end{equation*}

The differential~$\delta\colon\Omega^0(\mathcal{E})\to\Omega^1(\mathcal{E})$
is given by the formula $\delta(\omega)=\eval{\nabla^*(1)}_{\mathcal{E}}$,
with~$\nabla$ being a $\C$-differential operator from~$P$ to~$\Lah^n(\pi)$
such that $\rmdh\omega=\nabla(F)$ on~$\J^{\infty}(\pi)$.

The
differential~$\delta\colon\Omega^p(\mathcal{E})\to\Omega^{p+1}(\mathcal{E})$,
$p\geq1$, has the form
\begin{equation*}
  \delta(\Delta)(\phi_1,\dots,\phi_p)
  =\sum_{\alpha=1}^p(-1)^{\alpha+1}
  \ell_{\Delta,\phi_1,\dots,\hat{\phi}_{\alpha},\dots,\phi_p}(\phi_{\alpha})
  +{\eval{\nabla}_{\mathcal{E}}}^{*_1}(\phi_1,\dots,\phi_p),
\end{equation*}
where~$\nabla$ is a $\C$-differential operator~$\nabla\colon
P\times\kappa(\pi)\times\dots\times\kappa(\pi)\to\hat{\kappa}(\pi)$ that
satisfies the relation
\begin{equation}
  \label{eq:gjsis_eqs:21}
  \ell_F^*\Delta(\phi_1,\dots,\phi_{p-1})
  -\sum_{\alpha=1}^{p-1}\Delta^{*_{\alpha}}
  (\phi_1,\dots,\ell_F(\phi_{\alpha}),\dots,\phi_{p-1})
  =\nabla(F,\phi_1,\dots,\phi_{p-1})
\end{equation}
on~$\J^{\infty}(\pi)$.

In the case of an evolution equation $F=u_t-f$ we have
$\Delta\in\CDiffskalt(\kappa,\hat{\kappa})$ and the operator~$\nabla$ can be
chosen in the form
\begin{equation*}
  \nabla(h,\phi_1,\dots,\phi_{p-1})
  =-\ell_{\Delta,\phi_1,\dots,\phi_{p-1}}(h),
\end{equation*}
where $h\in P$.

If $p=1$ and $\psi\in\im\delta\subset\Omega^1(\mathcal{E})$ then we can put
$\nabla=-\ell_{\psi}^*$.

The above description of elements in~$\Omega^p(\mathcal{E})$ in terms of
$\C$-differential operators makes sense only for a given inclusion of the
equation~$\mathcal{E}$ to a jet space.  If we consider two inclusions

\begin{equation*}
  \xymatrixcolsep{1pc}
  \xymatrixrowsep{0.5pc}
  \xymatrix{
    &J^{\infty}(\pi_1) \\
    \mathcal{E}\ar[ur]\ar[dr] \\
    &J^{\infty}(\pi_2),
  }
\end{equation*}
so that the corresponding linearizations are equivalent and we have
diagram~\eqref{eq:gjsis_eqs:3}, then the operators $\Delta_1$ and $\Delta_2$
that define the same element of~$\Omega^p(\mathcal{E})$ with respect to two
inclusions are related as follows:
\begin{equation*}
  \Delta_1={\alpha'}^*\Delta_2
  (\alpha(\,\cdot\,),\dots,\alpha(\,\cdot\,)),\qquad
  \Delta_2={\beta'}^*\Delta_1
  (\beta(\,\cdot\,),\dots,\beta(\,\cdot\,)).
\end{equation*}

Elements of~$\Omega^1(\mathcal{E})$ are said to be \emph{cosymmetries} of the
equation~$\mathcal{E}$.  We say that isomorphism~\eqref{eq:gjsis_eqs:18} takes
a cosymmetry to its \emph{generating section} (or the \emph{characteristic})
belonging to~$\ker\bar{\ell}_F^*\subset\hat{P}$.

The ``Two-line Theorem'' by Vinogradov (see, e.g.,
\cite{KrasilshchikVinogradov:SCLDEqMP}) implies that
complex~\eqref{eq:gjsis_eqs:17} is exact at the term~$\Omega^0(\mathcal{E})$,
hence the conservation laws of~$\mathcal{E}$ form a subset of the space of
cosymmetries, $\cl(\mathcal{E})\subset\Omega^1(\mathcal{E})$.  The generating
section of a cosymmetry that belongs to $\delta(\cl(\mathcal{E}))$ is called
the \emph{generating section of the conservation law}.

As noted above, to compute the generating section of a conservation law we
extend it arbitrarily to a form $\omega\in\Lah^{n-1}(\pi)$ on the jet space
$\J^{\infty}(\pi)$, so that $\eval{\rmdh\omega}_{\mathcal{E}}=0$. Hence there
exists a $\C$-differential operator $\Delta\colon P\to\Lah^n(\pi)$ such that
$\rmdh\omega=\Delta(F)$; the element
$\psi={\eval{\Delta^*}_{\mathcal{E}}}(1)\in\hat{P}$ is the generating section
of the conservation law under consideration.

The generating section~$\psi$ of a conservation law can always be extended to
the jet space~$\J^\infty(\pi)$ in such a way that
$\langle\psi,F\rangle=\rmdh\omega$, with $\omega$ being a conserved current
for the same conservation law.

\begin{remark}
  The above describe procedure is one of the ways to compute the
  differential~$\delta\colon\Omega^0(\mathcal{E})\to\Omega^1(\mathcal{E})$
  in~\eqref{eq:gjsis_eqs:17} for an arbitrary
  equation~$\mathcal{E}$. If~$\mathcal{E}$ is presented in an evolutionary
  form then computation of the generating section is simpler and more
  straightforward. Namely, let~$t$, $x^1,\dots,x^n$ be the independent
  variables and~$[\omega]\in\Omega^0(\mathcal{E})$, where
  \begin{equation*}
    \omega=X\rmd x^1\wedge\dots\wedge x^n+\sum_{i=1}^n
    T_i\rmd t\wedge\rmd x^1\wedge\dots\wedge\hat{\rmd
      x}^i\wedge\dots\wedge\rmd x^n.
  \end{equation*}
  Then the corresponding generating section is
  \begin{equation*}
    \psi=\left(\frac{\delta X}{\delta u^1},\dots,\frac{\delta X}{\delta u^m}\right),
  \end{equation*}
  where~$\delta/\delta u^j$ is the variational derivative with respect to~$u^j$.
\end{remark}

\begin{example}
  The generating section of the conservation law from
  Example~\ref{exmp:gjsis_eqs:4} is equal to ~$1$.

  The generating sections of the conservation laws from
  Example~\ref{exmp:gjsis_eqs:5} are~$1$, $2u$, and~$u_{xx}+3u^2$.
\end{example}

The determining equation for the conservation laws (or, to be more precise, of
cosymmetries) of the equation~$\mathcal{E}=\{F=0\}$ is
\begin{equation}
  \label{eq:gjsis_eqs:19}
  \bar{\ell}_F^*(\psi)=0.
\end{equation}
This equation is dual to equation~\eqref{eq:gjsis_eqs:4} for symmetries.  The
computations involved in solving~\eqref{eq:gjsis_eqs:19} are very similar to
those used to compute symmetries.  However, unlike the case of symmetries, not
all solutions of~\eqref{eq:gjsis_eqs:19} give conservation laws.  To check if
the generating section of a cosymmetry corresponds to a conservation law we can
use the following corollary of the ``Two-line Theorem'' by Vinogradov (see
\cite{KrasilshchikVinogradov:SCLDEqMP}): if the de Rham cohomology
$H^n(\mathcal{E})=0$, complex~\eqref{eq:gjsis_eqs:17} is exact at the
term~$\Omega^1(\mathcal{E})$.  So, in this case the generating
section~$\psi\in\ker\bar{\ell}_F^*$ is the generating section of a
conservation law if and only if $\delta(\psi)=0$, that is, there exists a
self-adjoint operator $\Delta'={\Delta'}^*\in\C(P,\hat{P})$ such that
\begin{equation*}
  \ell_{\psi}+{\eval{\nabla}_{\mathcal{E}}}^*=\Delta'\bar{\ell}_F,
\end{equation*}
where $\nabla\in\C(\kappa(\pi),\hat{P})$ satisfies the equality
\begin{equation*}
  \ell_F^*(\psi)=\nabla(F)\quad\text{on~$\J^{\infty}(\pi)$.}
\end{equation*}

Note that the generating section~$\psi$ of a conservation law can always be
extended to the jet space $\J^{\infty}(\pi)$ in such a way that the horizontal
$n$-form $\langle\psi,F\rangle$ will be exact:
$\langle\psi,F\rangle=\rmdh\omega$, with $\eval{\omega}_{\mathcal{E}}$ being a
conserved current that corresponds to~$\psi.$

\subsection{A parallel with finite-dimensional differential geometry. III}
\label{sec:gjsis_eqs:parallel-with-finite}

In this section we begin compilation of a dictionary between the geometry of
normal differential equations and finite-dimensional differential geometry,
similar to the one for jet spaces from
Sections~\ref{sec:gjsis_jets:parallel-with-finite}
and~\ref{sec:gjsis_jets:parallel-with-finite-1}.

We start just as we did for jet spaces: we consider the space~$\mathcal{E}$
endowed with the Cartan distribution~$\C$ and take for the points of that
``manifold'' the maximal integral submanifolds of~$\C$, i.e., the solutions
of~$\mathcal{E}$.  Further, the dictionary reads:

\begin{eqnarray*}
  \text{\textbf{Manifold~$M$}}&&
  \text{\textbf{Normal equation~$\mathcal{E}$}} \\[1ex]
  \text{points}&\quad\longleftrightarrow\quad&\text{solutions} \\
  \text{functions~$C^{\infty}(M)$}&\quad\longleftrightarrow\quad&
  \text{conservation laws~$\cl(\mathcal{E})$} \\
  \text{the de Rham complex}&\quad\longleftrightarrow\quad&
  \text{complex~\eqref{eq:gjsis_eqs:17}}\\[-1.5ex]
  \text{of differential forms}&\quad\hphantom{\longleftrightarrow}\quad&
  \dots\xrightarrow{}\Omega^{p-1}(\mathcal{E})
  \xrightarrow{\delta}\Omega^p(\mathcal{E})
  \xrightarrow{}\dots \\
  \text{vector fields}&\quad\longleftrightarrow\quad&\text{symmetries} \\
  \text{the tangent bundle}&\quad\longleftrightarrow\quad&
  \text{the tangent covering~$\tau\colon\T(\mathcal{E})\to\mathcal{E}$}
\end{eqnarray*}

In addition to considerations from
Section~\ref{sec:gjsis_jets:parallel-with-finite} on jet spaces, this
dictionary is justified by the following facts.

First, on a finite-dimensional manifold the differential forms are
functions on the tangent bundle with odd fibres.  Correspondingly,
definition~\eqref{eq:gjsis_eqs:20} shows that
\begin{equation*}
\Omega^*(\mathcal{E})=\cl(\T(\mathcal{E})),
\end{equation*}
with elements of~$\Omega^p(\mathcal{E})$ given by fibre-wise $p$-linear
conserved currents.  Fibres of the tangent covering
$\tau\colon\T(\mathcal{E})\to\mathcal{E}$ are assumed to be odd.

Since elements of~$\Omega^p(\mathcal{E})$ can be understood as conservation
laws on~$\T(\mathcal{E})$, we can ask what are the generating section of these
conservation laws?  For the element of~$\Omega^p(\mathcal{E})$ that
corresponds to an operator $\Delta\in\CDiffsk_{p-1}(\kappa,\hat{P})$ the
generating section is $(-\nabla^{*_1},\Delta)$, where the operator~$\nabla$ is
given by~\eqref{eq:gjsis_eqs:21}.  Here we interpret skew-symmetric
$\C$-differential operators $\kappa\times\dots\times\kappa\to Q$, modulo
operators of the form~\eqref{eq:gjsis_eqs:22}, as elements of the module~$Q$
pulled back on~$\T(\mathcal{E})$.  This describes the
isomorphism~\eqref{eq:gjsis_eqs:23}.

Second, on a finite-dimensional manifold two natural actions of vector fields
on differential forms exist: the interior product and the Lie derivative.
Correspondingly, on an equation~$\mathcal{E}$ the evolution
field~$\Ev_{\phi}$, $\phi\in\sym\mathcal{E}$, induces the interior product and
the Lie derivative on~$\Omega^*(\mathcal{E})$ defined
by~\eqref{eq:gjsis_eqs:20}:
\begin{equation*}
  i_{\phi}\colon\Omega^p(\mathcal{E})\to\Omega^{p-1}(\mathcal{E}),\qquad
  L_{\phi}\colon\Omega^p(\mathcal{E})\to\Omega^{p+1}(\mathcal{E}).
\end{equation*}
These operations are related to the differential~$\delta$ by the usual
identity
\begin{equation*}
  L_{\phi}=\delta i_{\phi}+i_{\phi}\delta.
\end{equation*}

In terms of $\C$-differential operators, the interior product
$i_{\phi}\colon\Omega^p(\mathcal{E})\to\Omega^{p-1}(\mathcal{E})$ for $p>1$ is
the contraction of the operator with~$\phi$.  For $p=1$, the interior product
$i_{\phi}(\psi)$, $\psi\in\ker\bar{\ell}_F^*$, arises from the Green formula
and is the conservation law defined by the conserved
current~$\eval{\omega}_{\mathcal{E}}\in\Lah^{n-1}(\mathcal{E})$ such that
\begin{equation*}
  \langle\ell_F(\phi),\psi\rangle-\langle\phi,\ell_F^*(\psi)\rangle
  =\rmdh\omega
  \quad\text{on~$\J^{\infty}(\pi)$}.
\end{equation*}

If $\psi\in\hat{P}$ is a generating section of a conservation law,
then a symmetry~$\phi\in\sym\mathcal{E}$ acts on it by the formula
\begin{equation*}
  L_{\phi}(\psi)=\Ev_{\phi}(\psi)+\square^*(\psi),
\end{equation*}
with some operator $\square\in\CDiff(P,P)$ that satisfies the
equality~$\ell_F(\phi)=\square(F)$ on~$\J^{\infty}(\pi)$.

Third, on a finite-dimensional manifold a symplectic form gives rise
to a Poisson bracket.  The corresponding construction for an
equation~$\mathcal{E}$ relies on the notion of a \emph{symplectic
  structure} that is a closed element of~$\Omega^2(\mathcal{E})$.  We
do not assume that the symplectic form is non-degenerate, so the Poisson
bracket will be defined on a subset of~$\cl\mathcal{E}$ (recall that
conservation laws are analogues of functions on~$\mathcal{E}$.)

In terms of $\C$-differential operators, a symplectic structure is the
equivalence class of operators $\Delta\in\CDiff(\kappa,\hat{P})$ such that
\begin{equation}
  \label{eq:gjsis_eqs:25}
  \bar{\ell}_F^*\Delta=\Delta^*\bar{\ell}_F, \qquad
  \ell_{\Delta,\phi_2}(\phi_1)-\ell_{\Delta,\phi_1}(\phi_2)
  +{\eval{\nabla}_{\mathcal{E}}}^{*_1}(\phi_1,\phi_2)=0,
\end{equation}
where $\phi_1$,~$\phi_2\in\kappa$, $\nabla\colon P\times\kappa\to\kappa$ is a
$\C$-differential operator such that
\begin{equation*}
  \ell_F^*\Delta-\Delta^*\ell_F=\nabla(F,\,\cdot\,)
  \quad\text{on~$\J^{\infty}(\pi)$},
\end{equation*}
modulo operators of the form $\Delta'\bar{\ell}_F$,
$\Delta'\in\CDiff(P,\hat{P})$, ${\Delta'}^*=\Delta'$.

\begin{example}
  For the simplest WDVV equation
  \begin{equation*}
    u_{yyy}-u_{xxy}^2+u_{xxx}u_{xyy}=0
  \end{equation*}
  the operator~$D_x$ is a symplectic structure.
\end{example}

For evolution equations conditions~\eqref{eq:gjsis_eqs:25} amount to
\begin{equation*}
  \Delta^*=-\Delta, \qquad
    \ell_{\Delta,\phi_1}(\phi_2)-\ell_{\Delta,\phi_2}(\phi_1)
  =\ell_{\Delta,\phi_1}^*(\phi_2).
\end{equation*}

The construction of the Poisson bracket is similar to the one on a
finite-dimensional manifold.  Let $\Omega\in\Omega^2(\mathcal{E})$ be a
symplectic structure, i.e., $\delta(\Omega)=0$.  A conservation law with the
generating section~$\psi$ is called \emph{admissible} if there exists a
symmetry~$\phi\in\sym\mathcal{E}$ such that
\begin{equation}
  \label{eq:gjsis_eqs:24}
  \psi=i_{\phi}(\Omega).
\end{equation}

Symmetries that correspond to admissible conservation laws in the
sense of~\eqref{eq:gjsis_eqs:24} are called \emph{Hamiltonian
  symmetries}.

By definition, \emph{the Poisson bracket of two admissible conservation laws}
$\omega$ and~$\omega'$ with the generating sections $\psi$ and $\psi'$,
respectively, has the generating section
\begin{equation*}
  \{\omega,\omega'\}_{\Omega}
  =L_{\phi}(\psi')
  =i_{\phi}i_{\phi'}\Omega,
\end{equation*}
where $\phi$ and $\phi'$ are Hamiltonian symmetries corresponding to
conservation laws $\omega$ and~$\omega'$.  To the conservation law
$\{\omega,\omega'\}_{\Omega}$ there corresponds the Jacobi bracket
$\{\phi,\phi'\}$, so that Hamiltonian symmetries form a Lie algebra:
$[L_{\phi},L_{\phi'}]=L_{\{\phi,\phi'\}}$.

As we see from~\eqref{eq:gjsis_eqs:24}, a symplectic structure takes a
symmetry to a conservation law, in terms of
operators~\eqref{eq:gjsis_eqs:25}, this map has the form
$\psi=\Delta(\phi)$.

\subsection{Cotangent covering to a normal equation}
\label{sec:gjsis_eqs:cotangent-covering}

In the previous section we discussed geometry related to the tangent
covering and functions on it (forms).  But what about the cotangent
covering?

We have defined the tangent covering for an equation~$\mathcal{E}$ without
fixing an inclusion of~$\mathcal{E}$ to a jet space.  Dualizing such an
invariant definition requires use of rather complicated homological algebra,
so we will define the cotangent covering for an equation~$\mathcal{E}$
embedded to a jet space, $\mathcal{E}\subset\J^{\infty}(\pi)$, and then check
that the construction does not depend on the choice of the embedding.

For a normal equation~$\mathcal{E}$ given by a normal section $F=0, $~$F\in
P$, with $P$ being a module over~$\J^{\infty}(\pi)$, we define the
equation~$\Ts(\mathcal{E})\subset\JJ^{\infty}(\hat{P})$ by the equalities
\begin{equation*}
  \tilde{F}=0,\qquad\tilde{\ell}_F^*(\bi{p})=0,
\end{equation*}
where the tilde denotes the pullback to~$\JJ^{\infty}(\hat{P})$ and
$\bi{p}\in\tilde{\hat{P}}$ corresponds to the identity operator
$\hat{P}\to\hat{P}$ under the identification
$\tilde{\hat{P}}=\CDiff(\hat{P},\hat{P})$.  The natural projection
$\tau^*\colon\Ts(\mathcal{E})\to\mathcal{E}$ is called the
\emph{cotangent covering} to~$\mathcal{E}$.

In coordinates, we have $\bi{p}=(p^1,\dots,p^l)$ if the coordinates on
$\JJ^{\infty}(\hat{P})$ are $x^i$,~$u^j_I$,~$p^j_I$, with $u^j_I$ and $p^j_I$
being fibre coordinates along projections $\mathcal{E}\to M$ and
$\JJ^{\infty}(\hat{P})\to\mathcal{E}$, respectively.

Below we assume that not only equation~$\mathcal{E}$ is normal but the
operator $\bar{\ell}_F^*$ is normal (see Remark~\ref{rem:gjsis_eqs:1}) as
well.  Then $\Ts(\mathcal{E})$ will also be a normal equation.

\begin{remark}
  Obviously, for every~$\mathcal{E}$ the cotangent
  equation~$\Ts(\mathcal{E})$ is an Euler-Lagrange equation with the
  Lagrangian density $L=\langle F,\bi{p}\rangle$.  In applications,
  considering $\Ts(\mathcal{E})$ instead of~$\mathcal{E}$ is
  occasionally useful to handle the equation as though it were
  Lagrangian (see, e.g., \cite[Volume~1,
  Sections~3.2,~3.3]{MorseFeshbach:MTP}).  We refer
  to~\cite[Section~4.5.1]{FilippovSavchinShorokhov:VPNOp} and
  \cite[Section~5]{PopovychKunzingerIvanova:CLPSLPEq} for more details
  and references.
\end{remark}

If we have two inclusions of equation~$\mathcal{E}$ to jet spaces
\begin{equation*}
  \xymatrixcolsep{1pc}
  \xymatrixrowsep{0.5pc}
  \xymatrix{
    &J^{\infty}(\pi_1) \\
    \mathcal{E}\ar[ur]\ar[dr] \\
    &J^{\infty}(\pi_2),
  }
\end{equation*}
then the adjoint linearizations $\bar{\ell}_{F_1}^*$ and
$\bar{\ell}_{F_2}^*$ are equivalent:
\begin{equation}
  \label{eq:gjsis_eqs:29}
  \xymatrixcolsep{5pc}
  \xymatrix{
    \hat{P}_1\ar[r]_{\bar{\ell}_{F_1}^*}\ar@<.5ex>[d]^{{\beta'}^*}
    &\hat{\kappa}_1\ar@<.5ex>[d]^{\beta^*}\ar@/_1pc/@<-1ex>[l]_{s_1^*} \\
    \hat{P}_2\ar@<.5ex>[u]^{{\alpha'}^*}\ar[r]^{\bar{\ell}_{F_2}^*}
    &\hat{\kappa}_2\ar@<.5ex>[u]^{\alpha^*}\ar@/^1pc/@<1ex>[l]^{s_2^*}
  }
\end{equation}
such that
\begin{equation*}
  \bar{\ell}_{F_1}^*{\alpha'}^*=\alpha^*\bar{\ell}_{F_2}^*,\quad
  \bar{\ell}_{F_2}^*{\beta'}^*=\beta^*\bar{\ell}_{F_1}^*,\quad
  {\alpha'}^*{\beta'}^*=\id+s_1^*\bar{\ell}_{F_1}^*,\quad
  {\beta'}^*{\alpha'}^*=\id+s_2^*\bar{\ell}_{F_2}^*,
\end{equation*}
where the operators $\alpha$,~$\beta$,~$\alpha'$,~$\beta'$,~$s_1$,
and~$s_2$ are defined in~\eqref{eq:gjsis_eqs:3}.  Therefore, the
cotangent coverings constructed using $\ell_{F_1}^*$ and
$\ell_{F_2}^*$ are isomorphic, thus the cotangent coverings do not
depend on the choice of inclusion
$\mathcal{E}\subset\J^{\infty}(\pi)$.

Now we describe an isomorphism between $\sym\mathcal{E}$ and the
subspace of $\cl\Ts(\mathcal{E})$ that consists of conservation laws
with fibre-wise linear conserved currents.  Thus we will justify the
first parallel in the prolongation of our dictionary:

\begin{eqnarray}
  \text{\textbf{Manifold~$M$}}&&
  \text{\textbf{Normal equation~$\mathcal{E}$}} \nonumber \\[1ex]
  \text{the cotangent bundle}&\quad\longleftrightarrow\quad&
  \text{the cotangent covering~$\tau^*\colon
    \Ts(\mathcal{E})\to\mathcal{E}$} \nonumber \\
  \label{eq:gjsis_eqs:26}
  \text{multivector fields}&\quad\longleftrightarrow\quad&
  \text{conservation laws $\cl(\Ts(\mathcal{E}))$}
\end{eqnarray}

Let $\phi\in\kappa$ be the generating section of a symmetry
of~$\mathcal{E}$.  Extend it to an element $\phi\in\kappa(\pi)$ and
consider the Green formula
\begin{equation}
  \label{eq:gjsis_eqs:27}
  \langle\ell_F(\phi),\psi\rangle-\langle\phi,\ell_F^*(\psi)\rangle
  =\rmdh\omega(\phi,\psi),
\end{equation}
where $\psi\in\hat{P}$, $\omega(\phi,\psi)\in\Lah^{n-1}(\pi)$.  The mapping
$\psi\mapsto\omega(\phi,\psi)$ is a $\C$-differential operator
$\hat{P}\to\Lah^{n-1}(\pi)$, so that it gives rise to a closed form
$\omega_{\phi}\in\Lah^{n-1}(\Ts(\mathcal{E}))$.  The induced map
$\sym\mathcal{E}\to\cl\Ts(\mathcal{E})$, which takes the symmetry with the
generating functions~$\phi$ to the conservation law with the current
$\omega_{\phi}$, gives the desired isomorphism.

In fact, formula~\eqref{eq:gjsis_eqs:27} gives more.  It holds not only for
generating sections of symmetries of~$\mathcal{E}$, but also for an arbitrary
$\phi\in\kappa$, so that we obtain a map
$\kappa\to\Lah^{n-1}(\Ts(\mathcal{E}))$.  Since $\tilde{\kappa}$ is a direct
summand in $\kappa(\Ts(\mathcal{E}))$, the element $\omega(\phi,\psi)$ yields
an element $\rho\in\Omega^1(\Ts(\mathcal{E}))$.  One can prove that $\rho$
does not depend on the choice of inclusion $\mathcal{E}\to\J^{\infty}(\pi)$
used in its construction.

In terms of isomorphism~\eqref{eq:gjsis_eqs:23}, we have
$\rho=(\bi{p},0)$.

This is a very important element since it plays the r\^ole of the
canonical $1$-form $p\rmd q$ on a finite-dimensional cotangent space:

\begin{eqnarray*}
  \text{\textbf{Manifold~$M$}}&&
  \text{\textbf{Normal equation~$\mathcal{E}$}} \\[1ex]
  \text{the canonical $1$-form $p\rmd q$}&\quad\longleftrightarrow\quad&
  \text{$\rho\in\Omega^1(\Ts(\mathcal{E}))$} \\
  \text{the canonical symplectic form}
  &\quad\longleftrightarrow\quad&
  \text{canonical symplectic structure} \\[-1ex]
  \text{$\rmd p\wedge\rmd q$}
  &\quad\longleftrightarrow\quad&
  \text{$\Omega=\delta(\rho)\in\Omega^2(\Ts(\mathcal{E}))$}
\end{eqnarray*}

In terms of operators~\eqref{eq:gjsis_eqs:23}, the canonical symplectic
structure on~$\Ts(\mathcal{E})$ has the form
\begin{equation*}
  \Omega=\begin{pmatrix}
    0  & 1 \\
    -1 & 0.
  \end{pmatrix}
\end{equation*}

As to entry~\eqref{eq:gjsis_eqs:26} of our dictionary, we take it for the
definition of multivectors.  Since we are interested in skew-symmetric
multivectors, we assume the fibres of the cotangent covering
$\Ts(\mathcal{E})\to\mathcal{E}$ to be \emph{odd}.

We call conservation laws of~$\Ts(\mathcal{E})$ whose currents are fibre-wise
$p$-linear \emph{variational $p$-vectors} on~$\mathcal{E}$ and denote their
set by $D_p(\mathcal{E})$.  Thus, $D_0(\mathcal{E})=\cl\mathcal{E}$ and
$D_1(\mathcal{E})=\sym\mathcal{E}$.

For $D_p(\mathcal{E})$ we have a description in terms of
$\C$-differential operators similar to~\eqref{eq:gjsis_eqs:23}.  Namely,
\begin{equation*}
  D_p(\mathcal{E})=\Xi_p\big/\Xi_p^{\ell},
\end{equation*}
where $\Xi_p$ is a subset of $\CDiffsk_{p-1}(\hat{P},\kappa)$ that
consists of operators $\Delta\in\CDiffsk_{p-1}(\hat{P},\kappa)$ such that
\begin{equation*}
  \bar{\ell}_F\Delta(\psi_1,\dots,\psi_{p-1})
  -\sum_{\alpha=1}^{p-1}\Delta^{*_{\alpha}}
  (\psi_1,\dots,\bar{\ell}_F^*(\psi_{\alpha}),\dots,\psi_{p-1})=0
\end{equation*}
for all $\psi_1,\dots,\psi_{p-1}\in\hat{P}$, where~${^{*_{\alpha}}}$ denotes,
as before, the operation of taking adjoint with respect to the $\alpha$th
argument.  The subset $\Xi_p^{\ell}\subset\Xi_p$ consists of operators
$\Delta\in\Xi_p$ of the form
\begin{equation*}
  \Delta(\psi_1,\dots,\psi_{p-1})=\sum_{\alpha=1}^{p-1}(-1)^{\alpha+1}\Delta'
  (\bar{\ell}_F^*(\psi_{\alpha}),\psi_1,\dots,
  \hat{\psi}_{\alpha},\dots,\psi_{p-1})
\end{equation*}
for some $\C$-differential operators $\Delta'\colon
\hat{\kappa}\times\hat{P}\times\dots\times\hat{P}\to\kappa$.

In particular, for $p=2$, we have
\begin{equation*}
  D_2(\mathcal{E})
  =\sd{\Delta\in\CDiff(\hat{P},\kappa)}
  {\bar{\ell}_F\Delta=\Delta^*\bar{\ell}_F^*}\big/
  \sd{\Delta'\bar{\ell}_F^*}
  {\Delta'\in\CDiff(\hat{\kappa},\kappa),\ {\Delta'}^*=\Delta'}.
\end{equation*}

For the element of~$D_p(\mathcal{E})$ that corresponds to an operator
$\Delta\in\CDiffsk_{p-1}(\hat{P},\kappa)$ the generating section is
$(-\nabla^{*_1},\Delta)$, where the operator~$\nabla\colon
P\times\hat{P}\times\dots\times\hat{P}\to P$ is given by
\begin{equation}
  \label{eq:gjsis_eqs:28}
  \ell_F\Delta(\psi_1,\dots,\psi_{p-1})
  -\sum_{\alpha=1}^{p-1}\Delta^{*_{\alpha}}
  (\psi_1,\dots,\ell_F^*(\psi_{\alpha}),\dots,\psi_{p-1})
  =\nabla(F,\psi_1,\dots,\psi_{p-1}).
\end{equation}

In the case of evolution equation $F=u_t-f$ we have
$\Delta\in\CDiffskalt(\hat{\kappa},\kappa)$ and the operator~$\nabla$
can be chosen in the form:
\begin{equation*}
  \nabla(h,\psi_1,\dots,\psi_{p-1})
  =\ell_{\Delta,\psi_1,\dots,\psi_{p-1}}(h),
\end{equation*}
where $h\in P$.

Since we assume the fibres of cotangent covering to be odd, the
bracket defined on $\cl(\Ts(\mathcal{E}))$ by the canonical symplectic
structure will be the \emph{variational Schouten bracket}
\begin{equation*}
  \lshad\,\cdot\,,\cdot\,\rshad\,\colon D_k\times D_l\to D_{k+l-1}.
\end{equation*}

In terms of $\C$-differential operators, this bracket has the form:
\begin{multline*}
  \lshad\Delta_1,\Delta_2\rshad(\psi_1,\dots,\psi_{k+l-2})=\sum_{\sigma\in
  S_{k+l-2}^{l-1}}(-1)^\sigma
  \ell_{\Delta_2,\psi_{\sigma(1,l-1)}}(\Delta_1(\psi_{\sigma(l,k+l-2)})) \\
  -(-1)^{(k-1)l}\sum_{\sigma\in S_{k+l-2}^k}(-1)^\sigma
  \Delta_2(\nabla_1^{*_1}(\psi_{\sigma(1,k)}),
  \psi_{\sigma(k+1,k+l-2)}) \\
  -(-1)^{(k-1)(l-1)}\sum_{\sigma\in S_{k+l-2}^{k-1}}(-1)^\sigma
  \ell_{\Delta_1,\psi_{\sigma(1,k-1)}}(\Delta_2(\psi_{\sigma(k,k+l-2)})) \\
  +(-1)^{l-1}\sum_{\sigma\in S_{k+l-2}^l}(-1)^\sigma
  \Delta_1(\nabla_2^{*_1}(\psi_{\sigma(1,l)}),\psi_{\sigma(l+1,k+l-2)}),
\end{multline*}
where $\Delta_1\in D_k(P)$, $\Delta_2\in D_l(P)$, $\nabla_1$ and $\nabla_2$
are defined by~\eqref{eq:gjsis_eqs:28},
$\psi_1$,~\dots,~$\psi_{k+l-2}\in\hat{P}$,
cf.~\eqref{eq:gjsis_jets:41}.

The above description of variational multivectors in terms of
$\C$-differential operators makes sense only for a given inclusion of
equation~$\mathcal{E}$ to a jet space.  If we consider two inclusions

\begin{equation*}
  \xymatrixcolsep{1pc}
  \xymatrixrowsep{0.5pc}
  \xymatrix{
    &J^{\infty}(\pi_1) \\
    \mathcal{E}\ar[ur]\ar[dr] \\
    &J^{\infty}(\pi_2),
  }
\end{equation*}
so that the corresponding linearizations and adjoint linearizations are
equivalent and we have diagrams~\eqref{eq:gjsis_eqs:3}
and~\eqref{eq:gjsis_eqs:29}, then the operators $\Delta_1$ and $\Delta_2$ that
define the same element of~$D_p(\mathcal{E})$ with respect to the two
inclusions are related as follows:
\begin{equation*}
  \Delta_2=\alpha\Delta_1
  ({\alpha'}^*(\,\cdot\,),\dots,{\alpha'}^*(\,\cdot\,)),\qquad
  \Delta_1=\beta\Delta_2
  ({\beta'}^*(\,\cdot\,),\dots,{\beta'}^*(\,\cdot\,)).
\end{equation*}

An element~$A\in D_2(\mathcal{E})$ is called a \emph{Hamiltonian structure}
on~$\mathcal{E}$ if $\lshad A,A\rshad=0$.  Two Hamiltonian operators $A_1$ and
$A_2$ are said to be \emph{compatible} if $\lshad A_1,A_2\rshad=0$
(cf.~Section~\ref{sec:gjsis_jets:hamilt-form}).

\begin{remark}
  Using this definition of a Hamiltonian structure, the Hamiltonian formalisms
  on jet spaces, including the Magri scheme, explained in
  Section~\ref{sec:gjsis_jets:hamilt-form}, can be extended straightforwardly
  to the case of equations described here.
\end{remark}

\begin{example}
  The Camassa-Holm equation~\cite{CamassaHolm:ISWEPS}
  \begin{equation*}
    u_t-u_{txx}-uu_{xxx}-2u_xu_{xx}+3uu_x=0    
  \end{equation*}
has a bi-Hamiltonian structure:
  \begin{equation*}
    A_1=D_x,\qquad A_2=-D_t-uD_x+u_x.
  \end{equation*}

  This equation is often written in the form
  \begin{align*}
    &m_t+um_x+2u_xm=0, \\
    &m-u+u_{xx}=0.
  \end{align*}
  Then the bi-Hamiltonian structure takes the form
  \begin{equation*}
    A'_1=
    \begin{pmatrix}
      D_x & 0 \\
      D_x-D_x^3 & 0
  \end{pmatrix},\qquad
    A'_2=
    \begin{pmatrix}
      0 & -1 \\
      2mD_x+m_x & 0
    \end{pmatrix}.
  \end{equation*}
\end{example}

\begin{example}
  Let $\mathcal{E}$ be a bi-Hamiltonian equation given by $F=0$ and
  $A_1$ and $A_2$ be the corresponding Hamiltonian operators.  The
  Kupershmidt
  deformation~\cite{Kupershmidt:KInS,KerstenKrasilshchikVerbovetskyVitolo:InKD}~$\tilde{\mathcal{E}}$
  of~$\mathcal{E}$ has the form
  \begin{equation*}
    F+A_1^*(w)=0,\qquad A_2^*(w)=0,
  \end{equation*}
  where $w=(w^1,\dots,w^l)$ are new dependent variables.  For
  example, the KdV6
  equation~\cite{Karasu(Kalkani)KarasuSakovichSakovichTurhan:ANIGKdV}
  \begin{equation*}
    v_t+v_{xxx}+12vv_x-w_x=0,\qquad w_{xxx}+8vw_x+4wv_x=0,
  \end{equation*}
  is a Kupershmidt deformation of the KdV equation corresponding to the
  Hamiltonian operators~$A_1=D_x$ and~$A_2=D_x^3+8vD_x+4v_x$.
  
  The following two variational bivectors define a bi-Hamiltonian
  structures on~$\tilde{\mathcal{E}}$:
  \begin{equation*}
    \tilde{A}_1=
    \begin{pmatrix}
      A_1 & -A_1 \\
      0 & \ell_{F+A_1^*(w)+A_2^*(w)}
    \end{pmatrix},
    \qquad \tilde{A}_2=
    \begin{pmatrix}
      A_2 & -A_2 \\
      -\ell_{F+A_1^*(w)+A_2^*(w)} & 0
    \end{pmatrix}.
  \end{equation*}
\end{example}

\begin{example}
  The equation
  \begin{equation*}
  z_{yy} + (1/z)_{xx} +2 =0
\end{equation*}
associated with an integrable class of Weingarten
surfaces~\cite{BaranMarvan:InWSFC} is bi-Hamiltonian with operators
$D_x^2$ and $2zD_{xy}-z_yD_x+z_xD_y$.
\end{example}

\section{Nonlocal theory}
\label{sec:gjsis_nonloc:nonlocal-theory}

Nonlocal phenomena in the theory of integrable systems are quite common. Here
by \emph{nonlocality} we mean an extension of the initial system by new
variables (fields) that are related to the old ones by differential
relations. Perhaps, the simplest way to observe how nonlocal objects originate
is to analyse the action of recursion operators on symmetries.

\begin{example}
  \label{exmp:gjsis_nonloc:1}
  Consider recursion operator~\eqref{eq:gjsis_jets:53}
  $R=D_x^2+4u+2u_1D_x^{-1}$ that generates the higher KdV equations (it can be
  shown that successive application of~$R$ to the first symmetry~$\phi_1=u_1$
  results in polynomial expressions in~$u$, $u_1,\dots,u_k,\dots$, see,
  e.g.,~\cite{Krasilshchik:SMPLSH,Sergyeyev:LSGNInTDROpNApFS}). When one
  applies the operator~$R$ to the first $(x,t)$-dependent symmetry
  \begin{equation*}
    \bar{\phi}_1=tu_1+\frac{1}{6}
  \end{equation*}
  (the Galilean boost), this will result in the scaling symmetry
  \begin{equation*}
    \bar{\phi}_3=tu_3+(6tu+\frac{1}{3}x)u_1+\frac{2}{3}u,
  \end{equation*}
  but application of the recursion operator to~$\bar{\phi}_3$ leads to an
  expression that contains the nonlocal term~$D_x^{-1}(u)$ which can not be
  expressed in the geometrical terms introduced above\footnote{Of course, the
    Lenard operator itself contains a nonlocal summand, but we can consider it
    just as a convenient reformulation of the Magri
    relation~\eqref{eq:gjsis_jets:60}.}.

  An apparent way to incorporate this nonlocal object into the initial
  geometric setting is to introduce a new variable, say~$w$, that is related
  with the old one by~$w_x=u$. This relation, due to the KdV equation, implies
  another one:~$w_t=3u^2+u_{xx}$ and thus we shall result in the system
  \begin{equation*}
    u_t=6uu_x+u_{xxx},\qquad w_x=u,\qquad w_t=3u^2+u_{xx}.
  \end{equation*}
\end{example}

A general geometric formulation of this construction was first introduced
in~\cite{VinogradovKrasilshchik:MCHSNEvEqNS,KrasilshchikVinogradov:NTGDEqSCLBT}
and below we shall give a concise exposition of the theory together with a
number of applications.

\subsection{Differential coverings}
\label{sec:gjsis_nonloc:diff-cover}

The notion of a covering was already used in
Section~\ref{sec:gjsis_eqs:diff-equat} in the context of the tangent and
cotangent coverings. Here we discuss it in more detail.

Let~$\mathcal{E}\subset\J^\infty(\pi)$, where~$\pi\colon E\to M$, $\dim M=n$,
be an equation. Consider a locally trivial
bundle~$\tau\colon\tilde{\mathcal{E}}\to\mathcal{E}$ and endow the
manifold~$\tilde{\mathcal{E}}$ with an $n$-dimensional
distribution~$\tilde{\mathcal{\C}}$ in such a way that
\begin{enumerate}
\item\label{item:gjsis_nonloc:1} $\tilde{\mathcal{C}}$ is integrable and
\item\label{item:gjsis_nonloc:2} for any
  point~$\tilde{\theta}\in\tilde{\mathcal{E}}$ the
  restriction~$\eval{\rmd\tau}_{\tilde{\mathcal{\C}}_{\tilde{\theta}}}$ is a
  one-to-one correspondence between~$\tilde{\mathcal{\C}}_{\tilde{\theta}}$
  and the Cartan plane~$\C_{\tau(\tilde{\theta})}$.
\end{enumerate}
Then we say that~$\psi$ is endowed with the structure of a \emph{differential
  covering} (or simply a \emph{covering}, to be short) over~$\mathcal{E}$.

\begin{coordinates}
  Consider a trivialization of the bundle~$\tau$ and let~$w^1,\dots,w^j,\dots$
  be fibre-wise coordinates (the so-called \emph{nonlocal variables}). The
  number~$r$ of these coordinates is called the \emph{dimension} of~$\tau$.

  Let~$D_1,\dots,D_n$ be the total derivatives on~$\mathcal{E}$. By
  Property~(\ref{item:gjsis_nonloc:2}) in the definition of the
  distribution~$\tilde{\mathcal{\C}}$ there exist $\tau$-vertical vector
  fields~$X_1,\dots,X_n$ on~$\tilde{\mathcal{E}}$ such that the fields
  \begin{equation*}
    \tilde{D}_i=D_i+X_i,\qquad i=1,\dots,n, 
  \end{equation*}
  lie in~$\tilde{\mathcal{\C}}$. Then Property~(\ref{item:gjsis_nonloc:1}) is
  equivalent to the system of equations
  \begin{equation}
    \label{eq:gjsis_nonloc:1}
    D_i(X_j)-D_j(X_i)+[X_i,X_j]=0,\qquad 1\le i<j\le n,
  \end{equation}
  where~$X_1,\dots,X_n$ are $\tau$-vertical fields and~$D_i(X_j)$ denotes the
  component-wise action. Since the vector fields~$X_i$ are $\tau$-vertical
  fields, they can be presented in the form
  \begin{equation*}
    X_i=\sum_{j=1}^rX_i^j\frac{\partial}{\partial w^j},
  \end{equation*}
  where~$X_i^j$ are smooth functions on~$\tilde{\mathcal{E}}$,
  while~$\mathcal{E}$, as a manifold with distribution, is isomorphic to the
  infinite prolongation of the system of PDEs
  \begin{equation*}
    \frac{\partial w^j}{\partial x^i}=X_i^j,
    \qquad i=1,\dots,n,\quad j=1,\dots,r,
  \end{equation*}
  which extends the initial equation~$\mathcal{E}$ and is compatible over it
  due to~\eqref{eq:gjsis_nonloc:1}. This system is called the \emph{covering
    equation}.
\end{coordinates}

\begin{example}
  \label{exmp:gjsis_nonloc:2}
  Consider the one-dimensional covering over the KdV equation determined by
  \begin{equation*}
    \tilde{D}_x=D_x+u\frac{\partial}{\partial w},\qquad
    \tilde{D}_t=D_t+(3u^2+u_2)\frac{\partial}{\partial w}.
  \end{equation*}
  The covering equation in this case is
  \begin{equation*}
    \frac{\partial w}{\partial x}=u,\qquad
    \frac{\partial w}{\partial t}=3u^2+u_{xx}
  \end{equation*}
  and is isomorphic to the potential KdV equation~$w_t=3w_x^2+w_{xxx}$.

  Note the the relation between~$w$ and~$u$ may be expressed in the
  form~$w=\int u\rmd x$, or~$w=D_x^{-1}u$ and thus this is exactly the
  nonlocality that arose in Example~\ref{exmp:gjsis_nonloc:1}.
\end{example}

\begin{example}
  \label{exmp:gjsis_nonloc:4}
  Let again the base equation be the KdV and the covering be described by the
  system
  \begin{equation}
    \label{eq:gjsis_nonloc:2}
    X=u+w^2+\lambda,\qquad
    T=u_2+2wu_1+2u^2+2(w^2-\lambda)u-4\lambda(w^2+1),
  \end{equation}
  where~$\lambda\in\R$. Actually,~\eqref{eq:gjsis_nonloc:2} determines a
  one-parameter family of covering structures in the trivial
  bundle~$\mathcal{E}\times\R\to\mathcal{E}$, but the covering equation is
  isomorphic to the modified KdV equation~$w_t=6w^2w_x+w_{xxx}$ for
  any~$\lambda$. Of course, the covering under consideration is a geometric
  realization of the Miura transformation~\cite{Miura:KVEqGIRExNT}.
\end{example}

\begin{remark}
  \label{rem:gjsis_nonloc:1}
  Note that any differential substitution~$u=\phi(x,w,\dots,w_I,\dots)$ is
  associated with a covering over the initial equation, though this covering
  may be infinite-dimensional. For example, such is the covering over the KdV
  equation~$u_t-6uu_x+u_{xxx}=0$ determined with the Hirota substitution
  \begin{equation*}
    u=-2\frac{\partial^2}{\partial x^2}\ln w
  \end{equation*}
  (see~\cite{Hirota:ExSKVEqMCS}). Nevertheless, in spite of the infinite dimension
  of this covering, the covering space is isomorphic to the fourth-order
  scalar equation
  \begin{equation*}
    ww_{xt}-w_tw_x+w_{xxxx}w-4w_{xxx}w_x+3w_{xx}^2=0
  \end{equation*}
  in one unknown function.
\end{remark}

\begin{example}
  \label{exmp:gjsis_nonloc:6}
  Let~$P$ be the module of sections for some vector bundle~$\xi$
  over~$\mathcal{E}$. Then the
  bundle~$j_\infty^h\colon\JJ^\infty(P)\to\mathcal{E}$ of horizontal jets is
  an infinite-dimensional covering over~$\mathcal{E}$. If~$v_K^l$ are adapted
  coordinates in~$\JJ^\infty(P)$ then the total derivatives lifted
  to~$\JJ^\infty(P)$ are of the form
  \begin{equation}
    \label{eq:gjsis_nonloc:6}
    \tilde{D}_i=D_i+\sum_{l,K}v_{Ki}^l\frac{\partial}{\partial v_K^l}.
  \end{equation}
  This construction is generalised in the next example.
\end{example}

\begin{example}[$\Delta$-coverings]
  \label{exmp:gjsis_nonloc:5}
  Let~$\mathcal{E}$ be an equation and consider a $\C$-differential
  operator~$\Delta\colon P\to Q$, where~$P$ and~$Q$ are modules of sections
  for some vector bundles~$\xi$ and~$\zeta$
  over~$\mathcal{E}$. Let~$\Phi_\Delta\colon\JJ^\infty(P)\to \JJ^\infty(Q)$ be
  the corresponding morphism of vector bundles (see
  Proposition~\ref{prop:gjsis_jets:2}). Then, under natural conditions of
  non-degeneracy,~$\tilde{\mathcal{E}}_\Delta=\ker\Phi_\Delta$ is a sub-bundle
  in~$\xi_\infty\colon\JJ^\infty(P)\to\mathcal{E}$ that carries a natural
  structure of a covering: the total derivatives in this covering are obtained
  by restriction of the operators~\eqref{eq:gjsis_nonloc:6}
  to~$\tilde{\mathcal{E}}_\Delta$. We call this covering the
  \emph{$\Delta$-covering} over~$\mathcal{E}$.

  If the operator~$\Delta$ is locally given in the matrix
  form~$\Delta=\Vert\sum_Kd_{\alpha\beta}^KD_K\Vert$ then the
  subspace~$\tilde{\mathcal{E}}_\Delta\subset\JJ^\infty(\xi)$ is described by
  the relations
  \begin{equation*}
    \sum_{\alpha,K} d_{\alpha\beta}^Kv_K^\alpha=0
  \end{equation*}
  and their prolongations. Obviously, the tangent and cotangent coverings are
  particular cases of this construction with~$\Delta=\ell_{\mathcal{E}}$
  and~$\Delta=\ell_{\mathcal{E}}^*$, respectively\footnote{Here and below we
    consider normal equations and use the notation~$\ell_{\mathcal{E}}$
    instead of~$\bar{\ell}_F$ which is justified by the results of
    Subsection~\ref{sec:gjsis_eqs:line-diff-eqaut}}.

  $\Delta$-coverings play the key r\^{o}le in solving the following
  factorisation problem: let~$\Delta'\colon P'\to Q'$ be another
  $\C$-differential operator; how to find all operators~$A\colon P\to P'$ such
  that
  \begin{equation}
    \label{eq:gjsis_nonloc:7}
    \Delta'\circ A=B\circ\Delta,
  \end{equation}
  i.e., such that the diagram
  \begin{equation*}
    \begin{CD}
      P@>\Delta>>Q\\
      @VAVV@VVBV\\
      P'@>>\Delta'>Q'
    \end{CD}
  \end{equation*}
  is commutative? Note that any operator~$A$ of the form~$A=B'\circ\Delta$,
  where~$B'\colon Q\to P'$ is an arbitrary $\C$-differential operator, is a
  solution to~\eqref{eq:gjsis_nonloc:7}. Such solutions will be called
  \emph{trivial}.

  To find nontrivial solutions, first note that since~$\Delta'$ is a
  $\C$-differential operator it can be lifted to the covering just by changing
  the total derivatives~$D_i$ to the lifted ones~$\tilde{D}_i$. Denote this
  lift by~$\tilde{\Delta}'$. Second, let us put into correspondence to any
  operator~$A=\Vert\sum_Ka_{\alpha\beta}^K\Vert$ the vector-function
  \begin{equation*}
    \tilde{\Phi}_A=
    \eval{\left(\sum_{\alpha,K}a_{\alpha,1}^Kv_K^\alpha,
        \dots,\sum_{\alpha,K}a_{\alpha,r'}^Kv_K^\alpha\right)}_{\tilde{\mathcal{E}}_\Delta},
    \qquad r'=\dim P',
  \end{equation*}
  Then one has:
  \begin{proposition}
    \label{prop:gjsis_nonloc:1}
    Classes of solutions of Equation~\eqref{eq:gjsis_nonloc:7} modulo trivial
    ones are in one-to-one correspondence with solutions of the equation
    \begin{equation*}
      \tilde{\Delta}'(\tilde{\Phi}_A)=0.
    \end{equation*}
  \end{proposition}
  Operators satisfying~\eqref{eq:gjsis_nonloc:7} take elements of~$\ker\Delta$
  to those of~$\ker\Delta'$.
\end{example}

Consider system~\eqref{eq:gjsis_nonloc:1} that determines a covering structure
in the space~$\tilde{\mathcal{E}}$ and assume that the coefficients~$X_i^j$ of
the vertical vector fields~$X_i$ are independent of nonlocal
variables~$w^\alpha$. In this case,~\eqref{eq:gjsis_nonloc:1} reduces to
\begin{equation}
  \label{eq:gjsis_nonloc:3}
  D_i(X_j)=D_j(X_i),\qquad 1\le i<j\le n;
\end{equation}
the corresponding covering is called \emph{Abelian}. The covering in
Example~\ref{exmp:gjsis_nonloc:2} is an Abelian one, while the covering
associated with the Miura transformation (Example~\ref{exmp:gjsis_nonloc:4})
is not.

Let~$\dim\tau=1$ and define a differential horizontal $1$-form
on~$\mathcal{E}$ by setting
\begin{equation}
  \label{eq:gjsis_nonloc:4}
  \omega_\tau=\sum_{i=1}^nX_i\rmd x^i.
\end{equation}
Then~\eqref{eq:gjsis_nonloc:3} amounts to the equation
\begin{equation}
  \label{eq:gjsis_nonloc:5}
  \rmdh\omega_\tau=0,
\end{equation}
where~$\rmdh$ is the horizontal differential on~$\mathcal{E}$. Thus,
one-dimensional Abelian coverings over~$\mathcal{E}$ are in one-to-one
correspondence with closed horizontal $(n-1)$-forms.

In Example~\ref{exmp:gjsis_nonloc:4}, we presented a one-parameter family of
coverings over the KdV equation. Are these coverings different for different
values of the parameter~$\lambda$ or not and what the word ``different'' means
in this context? The answer is the following.

Consider an equation~$\mathcal{E}$ and two
coverings~$\tau_i\colon\tilde{\mathcal{E}}_i\to\mathcal{E}$, $i=1$, $2$,
over~$\mathcal{E}$. We say that these coverings are \emph{gauge equivalent}
(or simply \emph{equivalent}) if
there exists an
isomorphism~$\phi\colon\tilde{\mathcal{E}}_1\to\tilde{\mathcal{E}}_2$ of the
equations~$\tilde{\mathcal{E}}_1$ and~$\tilde{\mathcal{E}}_2$ such
that the diagram
\begin{equation*}
  \xymatrix{\tilde{\mathcal{E}}_1\ar[rr]^\phi\ar[rd]_{\tau_1}&
    &\ar[ld]^{\tau_2}\tilde{\mathcal{E}}_2\\
    &\mathcal{E}&
  }
\end{equation*}
is commutative, i.e.,~$\tau_2\circ\phi=\tau_1$. In this sense, all
coverings~\eqref{eq:gjsis_nonloc:2} are different, i.e., pair-wise
non-equivalent for different values of~$\lambda$. A general cohomological
technique to check whether a parameter can be eliminated or not was suggested
in~\cite{Marvan:HGCNSP,Marvan:SPP}.

We say that a covering~$\tau\colon\tilde{\mathcal{E}}\to\mathcal{E}$ is
\emph{trivial} if for any point~$\theta\in\mathcal{E}$ there exists a
neighbourhood~$\mathcal{U}\ni\theta$ such that
\begin{enumerate}
\item $\eval{\tau}_{\mathcal{U}}$ is a trivial bundle;
\item there exists an adapted coordinate system in~$\mathcal{U}$ for which the
  fields~$\tilde{D}_i$ are of the form~$\tilde{D}_i=D_i$, $i=1,\dots,\dim M$.
\end{enumerate}

\begin{theorem}
  \label{thm:gjsis_nonloc:1}
  There exists a one-to-one correspondence between equivalence classes of
  one-dimensional Abelian coverings over~$\mathcal{E}$ and elements of the
  horizontal cohomology group~$H_h^1(\mathcal{E})$ given
  by~\eqref{eq:gjsis_nonloc:4}. In particular, a covering~$\tau$ is trivial if
  and only if the form~$\omega_\tau$ is a co-boundary,
  i.e.,~$\omega_\tau=\rmdh(f)$.
\end{theorem}

Since~$H_h^1(\mathcal{E})$ coincides with the group of conservation laws
when~$\dim M=2$, Theorem~\ref{thm:gjsis_nonloc:1} allows one to construct
special type of coverings by conservation laws of the equation at hand.

\begin{example}
  \label{exmp:gjsis_nonloc:3}
  The Camassa-Holm equation
  \begin{equation*}
    u_t-u_{txx}+3uu_x=2u_xu_{xx}+uu_{xxx}
  \end{equation*}
  admits the conservation law
  \begin{equation*}
    \omega=(u-u_{xx})\rmd x+\frac{1}{2}(u_x^2-3u^2+2uu_{xx})\rmd t;
  \end{equation*}
  consequently, the corresponding covering is given by
  \begin{equation*}
    w_x=u-u_{xx},\qquad
    w_t=\frac{1}{2}(u_x^2-3u^2+2uu_{xx}).
  \end{equation*}
\end{example}

\begin{remark}
  Since the group~$H_h^1(\mathcal{E})$ is trivial for normal equations in the
  case~$\dim M>2$, this implies that such equations possess no nontrivial
  one-dimensional Abelian covering. Actually, there exist very strong
  indications that these equations do not have finite-dimensional coverings at
  all (see~\cite{Marvan:SSC}).
\end{remark}

\begin{example}[the KP equation]
  Consider the dispersionless Kadomtsev-Petviashvili equation
  \begin{equation*}
    (u_t-6uu_x+u_{xxx})_x=u_{yy}.
  \end{equation*}
  It admits an obvious covering
  \begin{equation}
    \label{eq:gjsis_nonloc:8}
    w_x=u_y,\qquad w_y=u_t-6uu_x+u_{xxx},
  \end{equation}
  which at first glance seems to be one-dimensional. But this is not the case,
  because Equations~\eqref{eq:gjsis_nonloc:8} do not contain information on
  the derivative~$w_t$. To incorporate these data, we must introduce infinite
  number of nonlocal variables~$w^0$, $w^1,\dots$ such that
  \begin{equation*}
    w^0=w,\ w_t^0=w^1,\,\dots,\,w_t^r=w^{r+1},\,\dots
  \end{equation*}
  and express their $x$- and $y$-derivatives
  using~\eqref{eq:gjsis_nonloc:8}. Thus, the covering is infinite-dimensional
  actually.
\end{example}

It was shown above that one-dimensional Abelian coverings can be constructed
using conservation laws of the equation. Another type of coverings is related
to Wahlquist-Estabrook prolongation
structures~\cite{WahlquistEstabrook:PSNEvEq,WahlquistEstabrook:PSNEvEqII,DoddFordy:PSQF}
and their description is based on the following \emph{ansatz}:
Let~$u_t=f(u,u_1,\dots,u_k)$ be a system of evolution equations,
$u=(u^1,\dots,u^m)$, $f=(f^1,\dots,f^m)$ being vectors and~$u_i$ denoting the
$i$th derivative with respect to~$x$. Let us look for coverings such that the
coefficients of the fields~$X$ and~$T$ in
\begin{equation*}
  \tilde{D}_x=D_x+X,\qquad\tilde{D}_t=D_t+T
\end{equation*}
depend on~$u$, $u_1,\dots,u_{k-1}$ and nonlocal variables only. Then, locally,
the description of such coverings locally reduces to representations of a
certain free Lie algebra (the so-called \emph{Wahlquist-Estabrook algebra}) in
vector fields on the fibre~$W$ of the trivial
bundle~$\tau\colon\mathcal{E}\times W\to\mathcal{E}$.

\begin{example}
  \label{exmp:gjsis_nonloc:7}
  Consider the potential KdV equation
  \begin{equation}
    \label{eq:gjsis_nonloc:9}
    u_t=u_x^2+u_{xxx}
  \end{equation}
  and let us describe coverings~$\tilde{D}_x=D_x+X$, $\tilde{D}_t=D_t+T$
  over~$\mathcal{E}$ such that the fields~$X$ and~$T$ depend on~$u$, $u_1$
  and~$u_2$ only. Straightforward computations show that all these coverings
  are of the form
  \begin{align}\label{eq:gjsis_nonloc:10}
    X&=u^2\mathbf{a}+u\mathbf{b}+\mathbf{c},\nonumber\\
    T&=(2uu_2-u_1^2+2u^2u_1)\mathbf{a}+(u_2+2uu_1)\mathbf{b}
    +u_1[\mathbf{c},\mathbf{b}]+\frac{1}{2}u^2[\mathbf{b},\mathbf{d}]
    +u[\mathbf{c},\mathbf{d}]+\mathbf{e},
  \end{align}
  where~$\mathbf{a}$, $\mathbf{b}$, $\mathbf{c}$, $\mathbf{d}$
  and~$\mathbf{e}$ are vector fields on the fibre~$W$ of the covering (i.e.,
  such that they do not depend on the equation coordinates) which enjoy the
  commutator relations
  \begin{gather*}
    2\mathbf{a}=[\mathbf{a},\mathbf{b}],\quad
    \mathbf{b}=[\mathbf{a},\mathbf{c}],\quad
    \mathbf{d}=2\mathbf{c}+[\mathbf{c},\mathbf{b}],\\[0pt]
    [\mathbf{a},\mathbf{d}]=[\mathbf{c},\mathbf{e}]=0,\\[0pt]
    [\mathbf{b},\mathbf{d}]+\frac{1}{2}[\mathbf{b},[\mathbf{b},\mathbf{d}]]=0,
    \quad [\mathbf{b},\mathbf{e}]+[\mathbf{c},[\mathbf{c},\mathbf{d}]]=0,\\[0pt]
    [\mathbf{a},\mathbf{e}]+[\mathbf{b},[\mathbf{c},\mathbf{d}]]
    +\frac{1}{2}[\mathbf{c},[\mathbf{b},\mathbf{d}]]=0.
  \end{gather*}
  Now, to find all Wahlquist-Estabrook coverings for~\eqref{eq:gjsis_nonloc:9}
  amounts to describing representations, as vector fields on~$W$, of the free
  Lie algebra generated by the elements~$\mathbf{a}$, $\mathbf{b}$,
  $\mathbf{c}$, $\mathbf{d}$, and~$\mathbf{e}$ with the above relations.

  If~$W=\R$ then all such representations, up to an isomorphism, are
  \begin{gather*}
    \mathbf{a}\mapsto\frac{\partial}{\partial w},\quad
    \mathbf{b}\mapsto(2w+\beta)\frac{\partial}{\partial w},\quad
    \mathbf{c}\mapsto(w^2+\beta w+\gamma)\frac{\partial}{\partial w},\\
    \mathbf{d}\mapsto-\Delta\frac{\partial}{\partial w},\quad
    \mathbf{e}\mapsto\Delta(w^2+\beta w+\gamma)\frac{\partial}{\partial w},
  \end{gather*}
  where~$\beta$, $\gamma\in\R$ and~$\Delta=\beta^2-4\gamma$. The corresponding
  one-dimensional coverings, up to gauge equivalence, are of the form
  \begin{equation*}
    X=(u^2+2wu+w^2+\gamma)\frac{\partial}{\partial w}
  \end{equation*}
  (the parameter~$\beta$ can be removed by a gauge transformation) and
  with~$T$ given by~\eqref{eq:gjsis_nonloc:10}; they are pair-wise inequivalent
  for different values of~$\gamma$.
\end{example}

\begin{remark}
  The term \emph{covering} also refers to the parallel between classical
  differential geometry and geometry of PDEs. Namely, if we define dimension
  of an equation~$\mathcal{E}$ (or of a jet space~$\J^\infty(\pi)$) as that of
  the corresponding Cartan distribution (i.e., the number of independent
  variables) then fibres of a differential covering become zero-dimensional,
  and this complies with the definition of a topological covering. That was
  the initial reason to name the object
  in~\cite{VinogradovKrasilshchik:MCHSNEvEqNS}.

  But the parallel goes far beyond this trivial
  observation. In~\cite{Igonin:AnCFGCPDEq}, a new powerful invariant (the
  \emph{fundamental Lie algebra}) of differential equations was proposed whose
  r\^{o}le in the theory of differential coverings is quite similar to the one
  that the fundamental group plays in topology. In particular, the fundamental
  Lie algebra allows one to enumerate (locally) all coverings over a given
  equation in the same way as conjugacy classes of subgroups of the
  fundamental group enumerate topological coverings. So, the dictionary
  evolved in the previous sections can be continued:
  
  \begin{eqnarray*}
    \text{\textbf{Manifold~$M$}}&&\text{\textbf{Differential equation~$\mathcal{E}$}} \\[1ex]
    \text{topological
      dimension}&\quad\longleftrightarrow\quad&\text{differential dimension} \\
    \text{topological
      coverings}&\quad\longleftrightarrow\quad&\text{differential coverings} \\
    \text{fundamental group}&\quad\longleftrightarrow\quad&\text{fundamental
      Lie algebra}
  \end{eqnarray*}
\end{remark}

Note also that using the fundamental Lie algebra technique the author
of~\cite{Igonin:AnCFGCPDEq} proved \emph{nonexistence} of B\"{a}cklund
transformations for some pairs of differential equations. It seems that it is
impossible to achieve such a result by other methods.

\subsection{Nonlocal symmetries}
\label{sec:gjsis_nonloc:nonlocal-symmetries}

The concept of a symmetry discussed in Section~\ref{sec:gjsis_eqs:diff-equat}
can be generalised to the nonlocal situation. Consider an example.

\begin{example}
  \label{exmp:gjsis_nonloc:8}
  Let
  \begin{equation}
    \label{eq:gjsis_nonloc:11}
    u_t=uu_x+u_{xx}
  \end{equation}
  be the Burgers equation (its Lie algebra of symmetries was fully described
  in~\cite{VinogradovKrasilshchik:MCHSNEvEqNS}). Direct computations show
  that~\eqref{eq:gjsis_nonloc:11} does not possess symmetries of the
  form~$\phi=\phi(x,t,u)$, but if one extends the setting by a new (nonlocal)
  variable~$w$ such that
  \begin{equation}
    \label{eq:gjsis_nonloc:12}
    w_x=u,\qquad w_t=\frac{1}{2}u^2+u_x
  \end{equation}
  then the equation~$\ell_{\mathcal{E}}(\phi)=0$ will acquire a new family of
  solutions of the form
  \begin{equation}
    \label{eq:gjsis_nonloc:13}
    \phi=(au-2a_x)\rme^{-\frac{1}{2}w},
  \end{equation}
  where~$a=a(x,t)$ is an arbitrary solution of the heat equation~$a_t=a_{xx}$.
\end{example}

The question is: can functions~\eqref{eq:gjsis_nonloc:13} be considered as
symmetries of Equation~\eqref{eq:gjsis_nonloc:11} in some natural sense? To
answer this question, consider an arbitrary
equation~$\mathcal{E}\subset\J^\infty(\pi)$ and a
covering~$\tau\colon\tilde{\mathcal{E}}\to\mathcal{E}$. We say that~$\phi$ is
a \emph{nonlocal symmetry} (or $\tau$-symmetry) of~$\mathcal{E}$ if it is a
symmetry of~$\tilde{\mathcal{E}}$.

\begin{coordinates}
  Let~$\mathcal{E}\subset\J^\infty(\pi)$ be an equation
  and~$\tau\colon\tilde{\mathcal{E}}\to\mathcal{E}$ be a covering locally
  given by the total derivatives
  \begin{equation*}
    \tilde{D}_i=D_i+\sum_jX_i^j\frac{\partial}{\partial w^j},\qquad
    i=1,\dots,\dim M.
  \end{equation*}
  Then any nonlocal $\tau$-symmetry is of the form
  \begin{equation}
    \label{eq:gjsis_nonloc:14}
    \tilde{\Ev}_\phi+\sum_j\psi^j\frac{\partial}{\partial w^j}.
  \end{equation}
  Here~$\phi$ is an $m$-component vector-function on~$\tilde{\mathcal{E}}$
  that satisfies the equation
  \begin{equation}\label{eq:gjsis_nonloc:15}
    \tilde{\ell}_{\mathcal{E}}(\phi)=0,
  \end{equation}
  $\psi^j$ are functions on~$\tilde{\mathcal{E}}$ such that
  \begin{equation}
    \label{eq:gjsis_nonloc:16}
    \tilde{D}_i(\psi^j)=\tilde{\ell}_{X_i^j}(\phi)
    +\sum_\alpha\frac{\partial X_i^j}{\partial w^\alpha}\psi^\alpha
  \end{equation}
  and
  \begin{equation}
    \label{eq:gjsis_nonloc:17}
    \tilde{\Ev}_\phi=
    \sum_{\substack{\text{over
          internal} \\ \text{coordinates}}}
    \tilde{D}_I(\phi^j)\frac{\partial}{\partial u_I^j}
  \end{equation}
  (recall that the ``tilde'' over a $\C$-differential operator denotes its
  natural lifting to the covering).
\end{coordinates}

\begin{example}
  \label{exmp:gjsis_nonloc:9}
  Let us consider Example~\ref{exmp:gjsis_nonloc:8} again. In the case of
  covering~\eqref{eq:gjsis_nonloc:12} Equations~\eqref{eq:gjsis_nonloc:16} take
  the form
  \begin{equation}
    \label{eq:gjsis_nonloc:18}
    \tilde{D}_x(\psi)=\phi,\qquad
    \tilde{D}_t(\psi)=u\phi+\tilde{D}_x(\phi).
  \end{equation}
  Consequently, for~$\phi$ of the form~\eqref{eq:gjsis_nonloc:13} we see that
  \begin{equation*}
    \psi=-2a\rme^{-\frac{1}{2}w}
  \end{equation*}
  satisfies~\eqref{eq:gjsis_nonloc:18}. Thus, the pair of functions~$\phi$
  and~$\psi$ determine a nonlocal symmetry of the Burgers equation in the
  sense of the above definition.
\end{example}

However, the situation of the previous example is not generic.

\begin{example}
  \label{exmp:gjsis_nonloc:10}
  Consider the covering
  \begin{equation*}
    w_x=u,\qquad w_t=3u^2+u_{xx}
  \end{equation*}
  over the KdV equation~$u_t=6uu_x+u_{xxx}$ and let us try to find nonlocal
  symmetries in this covering. In the case under consideration,
  Equations~\eqref{eq:gjsis_nonloc:16} acquire
  the form
  \begin{equation}
    \label{eq:gjsis_nonloc:19}
    \tilde{D}_x(\psi)=\phi,\qquad\tilde{D}_t(\psi)=6u\phi+\tilde{D}_x^2(\phi),
  \end{equation}
  while~\eqref{eq:gjsis_nonloc:15} is
  \begin{equation*}
    \tilde{D}_t(\phi)=6u_1\phi+6u\tilde{D}_x(\phi)+\tilde{D}_x^3(\phi).
  \end{equation*}
  The simplest solution of the last equation that depends on~$w$ is
  \begin{equation*}
    \phi=tu_5+\left(10tu+\frac{1}{3}x\right)u_3+
    4\left(5tu_1+\frac{1}{3}\right)u_2+
    2\left(15tu^2+xu+\frac{1}{3}w\right)u_1+\frac{8}{3}u^2.
  \end{equation*}
  But solving~\eqref{eq:gjsis_nonloc:19} with~$\phi$ of the above form leads
  to a contradiction: no function~$\psi$ exists on~$\mathcal{E}$ such
  that~\eqref{eq:gjsis_nonloc:19} is valid for our~$\phi$.

  Nevertheless, if we introduce another nonlocal variable~$w'$ satisfying
  \begin{equation*}
    w_x'=u^2,\qquad w_t'=4u^3-u_x^2+2uu_{xx}
  \end{equation*}
  then~\eqref{eq:gjsis_nonloc:19} will be resolved in the new setting.

  But a similar problem arises at the next step: now we need to reconstruct
  the coefficient~$\psi'$ at~$\partial/\partial w'$.
\end{example}

The procedure we encountered in Example~\ref{exmp:gjsis_nonloc:10} is typical
and we shall describe it in general terms
now. Let~$\tau\colon\tilde{\mathcal{E}}\to\mathcal{E}$ be a covering and
denote by~$\mathcal{F}$ and~$\tilde{\mathcal{F}}$ the function algebras
on~$\mathcal{E}$ and~$\tilde{\mathcal{E}}$, respectively. An $\R$-linear
map~$X\colon\mathcal{F}\to\tilde{\mathcal{F}}$ is called a
\emph{$\tau$-shadow} if
\begin{enumerate}
\item $X$ is a derivation, i.e.,
  \begin{equation*}
    X(fg)=fX(g)+gX(f)
  \end{equation*}
  for all~$f$, $g\in\mathcal{F}$;
\item the action of~$X$ preserves the Cartan distribution,
  i.e.,~$\mathrm{L}_X(\omega)\in\Lac(\tilde{\mathcal{E}})$ as soon
  as~$\omega\in\Lac(\mathcal{E})$ (or, equivalently, for any Cartan
  field~$\tilde{Y}$ on~$\tilde{\mathcal{E}}$ and its restriction~$Y$
  to~$\mathcal{F}$ the commutator~$[X,\tilde{Y}]=XY-\tilde{Y}X$ is a Cartan
  field again).
\end{enumerate}
In particular, any symmetry of~$\tilde{\mathcal{E}}$ can be considered as a
shadow in an arbitrary covering~$\tau$.

\begin{coordinates}
  Let~$\mathcal{U}$ be the set of internal coordinates on~$\mathcal{E}$. Then
  any $\tau$-shadow is given by the formula
  \begin{equation*}
    \tilde{\Ev}_\phi=
    \sum_{u_I^j\in\mathcal{U}}\tilde{D}_I(\phi^j)\frac{\partial}{\partial
    u_I^j},
  \end{equation*}
  where~$\phi^1,\dots,\phi^m$ are functions on~$\tilde{\mathcal{E}}$
  and~$\tilde{D}_1,\dots,\tilde{D}_n$ are total derivatives
  on~$\tilde{\mathcal{E}}$ (cf.~\eqref{eq:gjsis_nonloc:17}). The vector
  function~$\phi=(\phi^1,\dots,\phi^m)$ must satisfy the equation
  \begin{equation*}
    \tilde{\ell}_{\mathcal{E}}(\phi)=0.
  \end{equation*}
\end{coordinates}

We say that a $\tau$-shadow~$X$ is \emph{reconstructed} in~$\tau$ if there
exists a nonlocal $\tau$-symmetry~$\tilde{X}$ such
that~$\eval{\tilde{X}}_{\mathcal{F}}=X$. As Examples~\ref{exmp:gjsis_nonloc:9}
and~\ref{exmp:gjsis_nonloc:10} show, not all shadows can be reconstructed in a
straightforward way. A general result that describes the reconstruction
procedure was proved in~\cite{Khorkova:CLNS} (see
also~\cite{KrasilshchikVinogradov:NTGDEqSCLBT}):

\begin{proposition}
  \label{prop:gjsis_nonloc:2}
  Let~$\tau\colon\tilde{\mathcal{E}}\to\mathcal{E}$ be a covering and~$X$ be a
  $\tau$-shadow. Then there exists another
  covering~$\bar{\tau}\colon\bar{\mathcal{E}}\to\tilde{\mathcal{E}}$ and a
  $\bar{\tau}$-shadow~$\bar{X}$ such that~$\eval{\bar{X}}_{\mathcal{F}}=X$.
\end{proposition}

Thus, putting~$\tau=\tau_0$ and~$\tau_{i+1}=\tilde{\tau}_i$ and applying
Proposition~\ref{prop:gjsis_nonloc:2} sufficiently (maybe, infinitely) many
times we shall arrive to a covering in which the given shadow is
reconstructed.

\begin{coordinates}
  Actually, the results of~\cite{Khorkova:CLNS} not just state the existence
  of the needed covering but provide a canonical way to construct the one. The
  construction is in a sense tautological and mimics
  relations~\eqref{eq:gjsis_nonloc:16}.

  \begin{proposition}
    \label{prop:gjsis_nonloc:3}
    Let~$\tau\colon\tilde{\mathcal{E}}\to\mathcal{E}$ be a covering over an
    equation~$\mathcal{E}$ with nonlocal coordinates~$w^1,\dots,w^j,\dots$
    and~$\tilde{D}_1,\dots,\tilde{D}_n$ be total derivatives in this
    covering. Let also~$\phi$ be a $\tau$-shadow. Then:
    \begin{enumerate}
    \item the relations
      \begin{equation}
        \label{eq:gjsis_nonloc:20}
        \frac{\partial\tilde{w}^j}{\partial x^i}=\tilde{\ell}_{X_i^j}(\phi)
        +\sum_\alpha\frac{\partial X_i^j}{\partial
          w^\alpha}\tilde{w}^\alpha,\qquad
        i=1,\dots,n,\quad j=1,\dots,\dim\tau,
      \end{equation}
      define a covering over~$\tilde{\mathcal{E}}$ whose dimension equals that
      of~$\tau$;
    \item equations~\eqref{eq:gjsis_nonloc:16} are solvable in this covering.
    \end{enumerate}
  \end{proposition}
\end{coordinates}

\begin{remark}
  It follows from~\eqref{eq:gjsis_nonloc:20} that for an Abelian
  covering~$\tau$ the covering~$\tilde{\tau}$ is Abelian as well. Hence, at
  every step of reconstruction obstructions to
  solving~\eqref{eq:gjsis_nonloc:16} lie in the horizontal cohomology group of
  the corresponding equation. Consequently,
  if~$\tau\colon\tilde{\mathcal{E}}\to\mathcal{E}$ is an Abelian covering
  and~$H_h^1(\tilde{\mathcal{E}})=0$ then any $\tau$-shadow can be
  reconstructed to a nonlocal $\tau$-symmetry. In particular, any local
  symmetry of~$\mathcal{E}$ can be lifted to~$\tilde{\mathcal{E}}$.

  Let~$\mathcal{E}$ be an equation and~$\{[\omega^\alpha]\}$,
  $\omega^\alpha\in\Lah^1(\mathcal{E})$, be an $\R$-basis of the
  group~$H_h^1(\mathcal{E})$. Assume that
  \begin{equation*}
    \omega^\alpha=X_1^\alpha\rmd x^1+\dots+X_n^\alpha\rmd x^n
  \end{equation*}
  and consider the covering~$\tau_1\colon\mathcal{E}_1\to\mathcal{E}$
  determined by
  \begin{equation*}
    \frac{\partial w^\alpha}{\partial x^i}=X_i^\alpha
  \end{equation*}
  for all~$\alpha$ and~$i=1,\dots,n$. For~$\mathcal{E}_1$ let us construct the
  covering~$\tau_2\colon\mathcal{E}_2\to\mathcal{E}_1$ in a similar way,
  etc. The covering~$\tau_*\colon\mathcal{E}_*\to\mathcal{E}$ obtained as the
  inverse limit of the sequence
  \begin{equation*}
    \dots\xrightarrow{\tau_{i+1}}\mathcal{E}_i\xrightarrow{\tau_i}
    \mathcal{E}_{i-1}
    \xrightarrow{\tau_{i-1}}\dots\xrightarrow{\tau_2}\mathcal{E}_1
    \xrightarrow{\tau_1}\mathcal{E}
  \end{equation*}
  is called the \emph{universal Abelian covering} over~$\mathcal{E}$. By
  construction,~$H_h^1(\mathcal{E}_*)=0$.

  \begin{proposition}
    For an arbitrary Abelian
    covering~$\tau\colon\tilde{\mathcal{E}}\to\mathcal{E}$, there exists a
    uniquely (up to a gauge equivalence) defined morphism
    \begin{equation*}
      \xymatrix{
        \mathcal{E}_*\ar[rr]^{f}\ar[dr]_{\tau_*}&&
        \tilde{\mathcal{E}}\ar[dl]^{\tau}\\
        &\mathcal{E}&
      }
    \end{equation*}
    and any $\tau$-shadow can be reconstructed in~$\tau_*$. In particular, any
    symmetry of~$\mathcal{E}$ can be lifted to~$\mathcal{E}_*$.
  \end{proposition}
\end{remark}

\begin{remark}
  Though the covering~$\tilde{\tau}$ whose existence is stated in
  Proposition~\ref{prop:gjsis_nonloc:2} is determined canonically by the
  shadow~$X$, the new shadow~$\tilde{X}$ is not unique, but is defined up to
  an infinitesimal gauge symmetries of~$\tau$, i.e., up
  to~$Y\in\sym(\tilde{\mathcal{E}})$ such that~$\eval{Y}_{\mathcal{F}}=0$. Due
  to~\eqref{eq:gjsis_nonloc:16}, these symmetries are given by the equations
  \begin{equation}
    \label{eq:gjsis_nonloc:21}
    \tilde{D}_i(\psi^j)=
    \sum_\alpha\frac{\partial X_i^j}{\partial w^\alpha}\psi^\alpha
  \end{equation}
  and are of the form~$Y=\sum_\alpha\psi^\alpha\partial/\partial w^\alpha$.

  \begin{example}
    For Example~\ref{exmp:gjsis_nonloc:2}, Equations~\eqref{eq:gjsis_nonloc:21}
    take the form
    \begin{equation*}
      \tilde{D}_x(\psi)=0,\qquad
      \tilde{D}_t(\psi)=0
    \end{equation*}
    and consequently infinitesimal gauge symmetries
    are~$\gamma\partial/\partial w$ in this case, $\gamma\in\R$.

    On the other hand, if we consider Example~\ref{exmp:gjsis_nonloc:4}
    then Equations~\eqref{eq:gjsis_nonloc:21} are written as
    \begin{equation*}
      \tilde{D}_x(\psi)=2w\psi,\qquad
      \tilde{D}_t(\psi)=2(u_1+2uw-4\lambda w)\psi.
    \end{equation*}
    The only solution of this system is~$\psi=0$ and thus there is no
    ambiguity in shadow reconstruction in this case.
  \end{example}

  Non-uniqueness of the solution to the problem of reconstruction leads, in
  turn, to the problem of commutation for shadows: no well defined way to
  compute the Lie bracket of shadows is known. This problem was first
  indicated in~\cite{OlverSandersWang:GS}. A way to solve it was suggested
  in~\cite{VerbovetskyGolovkoKrasilshchik:LBNS}, but a practical realization
  of the approach is somewhat cumbersome.
\end{remark}

To conclude this subsection, let us make a remark also related to the problem
of reconstruction. Let~$\tau\colon\tilde{\mathcal{E}}\to\mathcal{E}$ be a
finite-dimensional covering and~$X$ be a symmetry of~$\mathcal{E}$. One can
(at least, locally) lift~$X$ to~$\tilde{\mathcal{E}}$ in an arbitrary
way. Then, if~$X$ is an integrable vector field (i.e., if it possesses the
corresponding one-parameter group of transformations) then the lifted
field~$\tilde{X}$ is integrable as well. If~$\tilde{X}$ is a symmetry
of~$\tilde{\mathcal{E}}$ then this means that we managed to reconstruct~$X$ up
to a nonlocal symmetry in the covering~$\tau$.

Conversely, consider the one-parameter group of
transformations~$\{\tilde{A}_\lambda\}$ corresponding to~$\tilde{X}$ and for
any~$\lambda\in\R$ define an $n$-dimensional distribution~$\tilde{\C}^\lambda$
on~$\tilde{\mathcal{E}}$ by
\begin{equation}
  \label{eq:gjsis_nonloc:22}
  \tilde{\C}^\lambda\colon\theta\mapsto\tilde{\C}_\theta^\lambda=
  \tilde{A}_{\lambda,*}\left(\tilde{\C}_{\tilde{A}_\lambda^{-1}(\theta)}\right)
\end{equation}
where~$\theta\in\tilde{\mathcal{E}}$, $\tilde{\C}_\theta$ is the Cartan plane
at the point~$\theta$ and~$F_*$ denotes the differential of the map~$F$.

\begin{proposition}
  Correspondence~\eqref{eq:gjsis_nonloc:22} determines a one-parameter
  family~$\tau_\lambda$ of pair-wise inequivalent coverings
  over~$\tilde{\mathcal{E}}$ such that~$\tau_0=\tau$.
\end{proposition}

\begin{example}
  Take the covering
  \begin{equation*}
    X=u+w^2,\qquad T=u_2+2wu_1+2u^2+2w^2u 
  \end{equation*}
  over the KdV equation and apply Proposition~\ref{prop:gjsis_nonloc:3} using
  the Galilean boost $tu_1+1/6$ for the symmetry~$X$. This will result in the
  Miura covering described in Example~\ref{exmp:gjsis_nonloc:4}.
\end{example}

Not all one-parameter families of coverings can be obtained by this procedure
(a counter-example can be found in~\cite{Cieslinski:NSWAlIsInG,Cieslinski:GInSPCNNSS}). But a weaker
result was proved in~\cite{IgoninKrasilshchik:OnPFBT}:

\begin{theorem}
  Let~$\tau_\lambda\colon\tilde{\mathcal{E}}\to\mathcal{E}$ be a one-parameter
  family of coverings regarded as a deformation of the
  covering~$\tau=\tau_0$. Then the corresponding infinitesimal deformation is
  a $\tau$-shadow.
\end{theorem}

\subsection{B\"{a}cklund transformations and zero-curvature representations}
\label{sec:gjsis_nonloc:zero-curv-repr}

In conclusion, let us briefly discuss how the constructions of B\"{a}cklund
transformations~\cite{RogersShadwick:BTTAp,RogersSchief:BDT} and
zero-curvature representations~\cite{AblowitzKaupNewellSegur:InSTFAnNP} are
translated to the geometrical language.

A \emph{B\"{a}cklund transformation} between two equations~$\mathcal{E}'$
and~$\mathcal{E}''$ with the unknown functions~$u'$ and~$u''$, respectively,
is another equation~$\mathcal{E}$ in unknown functions both~$u'$ and~$u''$
such that for any solution~$u'$ of~$\mathcal{E}'$ a solution~$u''$
of~$\mathcal{E}$ is a solution of~$\mathcal{E}''$ as well and vice
versa. If~$\mathcal{E}'$ coincides with~$\mathcal{E}''$ then one speaks about
\emph{auto-B\"{a}cklund transformation}.

\begin{example}
  \label{exmp:gjsis_nonloc:11}
  Consider the sine-Gordon equation
  \begin{equation}
    \label{eq:gjsis_nonloc:23}
    u_{xy}=\sin u.
  \end{equation}
  Then the system
  \begin{equation}
    \label{eq:gjsis_nonloc:24}
    v_y-u_y=2\lambda\sin\frac{v+u}{2},\qquad
    v_x+u_x=\frac{2}{\lambda}\sin\frac{v-u}{2},
  \end{equation}
  where~$\lambda\neq 0$ is a real parameter, determines the classical
  auto-B\"{a}cklund transformation (a one-parameter family, actually)
  for~\eqref{eq:gjsis_nonloc:23}, see~\cite{DoddBullough:BTSGEq}.
\end{example}

\begin{example}
  \label{exmp:gjsis_nonloc:12}
  The second example (which now can also be considered as a classical one) was
  found in~\cite{WahlquistEstabrook:PSNEvEq}. It is of the form
  \begin{align}
    \label{eq:gjsis_nonloc:25}
    &\left(\frac{v+w}{2}\right)_x +\left(\frac{v-w}{2}\right)^2+\lambda^2=0,
    \quad\lambda\in\R,\nonumber\\
    &\left(\frac{v-w}{2}\right)_t+
    6\left(\frac{v+w}{2}\right)_x\left(\frac{v-w}{2}\right)_x+
    \left(\frac{v-w}{2}\right)_{xxx}=0
  \end{align}
  and relates solutions of the KdV equation to each other (or, to be more
  precise, system~\eqref{eq:gjsis_nonloc:25} is a B\"{a}cklund transformation
  for the potential KdV equation, while solutions of the KdV itself are
  obtained by~$u=v_x$).
\end{example}

Analysis of these two examples (as well as other ones) shows that a
B\"{a}cklund transformation between equations~$\mathcal{E}'$
and~$\mathcal{E}''$ is a diagram
\begin{equation*}
  \xymatrix{
    &\mathcal{E}\ar[dl]_{\tau'}\ar[dr]^{\tau''}&\\
    \mathcal{E}'&&\mathcal{E}''\rlap{,}
  }
\end{equation*}
where~$\tau'$ and~$\tau''$ are coverings. The correspondence between solutions
of~$\mathcal{E}'$ and~$\mathcal{E}''$ is achieved in the following
way. Let~$u'=u'(x')$ be a solution of~$\mathcal{E}'$ and assume that~$\tau'$
is a finite-dimensional covering. Then the Cartan distribution of the
equation~$\mathcal{E}$ induces on the finite-dimensional
manifold~$\mathcal{E}_{u'}=(\tau')^{-1}(u')\subset\mathcal{E}$ an
$n$-dimensional integrable distribution. In the vicinity of a generic point
the latter possesses a $(\dim\tau')$-parameter family of maximal integral
manifolds that are projected to~$u'$ by~$\tau'$ and to the corresponding
family of solutions of~$\mathcal{E}''$ by~$\tau''$. Generically, such a
correspondence is non-trivial provided~$\tau'$ and~$\tau''$ are not gauge
equivalent.

\begin{example}
  \label{exmp:gjsis_nonloc:13}
  A common way to construct non-trivial B\"{a}cklund transformations is the
  following. Let~$\tau\colon\tilde{\mathcal{E}}\to\mathcal{E}$ be a covering
  and~$f\colon\tilde{\mathcal{E}}\to\tilde{\mathcal{E}}$ be a finite symmetry
  of~$\tilde{\mathcal{E}}$, i.e., a diffeomorphism preserving the Cartan
  distribution. Then the composition~$\tau'=\tau\circ f$ is a covering as well
  and the pair~$(\tau,\tau')$ is an auto-B\"{a}cklund transformation
  for~$\mathcal{E}$. If~$f$ is not a gauge equivalence then this
  transformation is non-trivial.

  Consider covering~\eqref{eq:gjsis_nonloc:2} from
  Example~\ref{exmp:gjsis_nonloc:4} and note that the change of the nonlocal
  variable~$w\leftrightarrow-w$ is a symmetry of the covering equation (the
  mKdV one), but is not a gauge symmetry of the covering itself. Thus, for any
  value of the parameter~$\lambda$ we get an auto-B\"{a}cklund transformation
  of the KdV equation, i.e., a one-parameter family of B\"{a}cklund
  transformations.

  Note that the Wahlquist-Estabrook construction
  (Example~\ref{exmp:gjsis_nonloc:12}) is a consequence of the latter one.
\end{example}

\begin{remark}
  Families of B\"{a}cklund transformations, like the ones from
  Examples~\ref{exmp:gjsis_nonloc:11} and~\ref{exmp:gjsis_nonloc:12}, give one
  an opportunity to construct special exact solutions of integrable equations
  (such as multi-kink solutions for the sine-Gordon equation, multi-soliton
  solutions for the KdV, etc.). The construction uses the \emph{nonlinear
    superposition principle} which, in turn, is based on the following
  informal statement:
  
  \begin{theorem}[the Bianchi Permutability Theorem]
    Assume that an equation~$\mathcal{E}$ possesses a one-parameter family of
    auto-B\"{a}cklund transformations~$\mathcal{B}_\lambda$ and
    let~$\lambda\in\R$ be the parameter. For any solution~$u=u(x)$
    of~$\mathcal{E}$ denote by~$\mathcal{B}_\lambda(u)$ the set of solutions
    obtained from~$u$ by means of~$\mathcal{B}_\lambda$. Then for
    any~$\lambda_1\neq\lambda_2$ there exists a
    solution~$u_{\lambda_1,\lambda_2}\in
    \mathcal{B}_{\lambda_1}(\mathcal{B}_{\lambda_2}(u))\cap
    \mathcal{B}_{\lambda_2}(\mathcal{B}_{\lambda_1}(u))$ that is expressed as
    a bi-differential operator applied to some
    solutions~$u_1\in\mathcal{B}_{\lambda_1}(u)$
    and~$u_2\in\mathcal{B}_{\lambda_2}(u)$.
  \end{theorem}

  This ``theorem'' was first observed by Bianchi in~\cite{Bianchi:LGD} (see
  also~\cite{RogersShadwick:BTTAp}) in application to the sine-Gordon equation
  (Example~\ref{exmp:gjsis_nonloc:11}) and since then dozens of examples were
  computed, but nevertheless a general formulation of this statement (and,
  consequently, its general proof) is unknown to us. Some hints to a rigorous
  approach to the problem can be found in~\cite{Marvan:SLPBT}.
\end{remark}

Geometrical theory of B\"{a}cklund transformations is also related to an
unorthodox approach to recursion operators~\cite{Marvan:AnLROp}. Consider an
equation~$\mathcal{E}$ and its tangent
covering~$\tau\colon\T(\mathcal{E})\to\mathcal{E}$ (see
Section~\ref{sec:gjsis_eqs:diff-equat}). Recall that symmetries
of~$\mathcal{E}$ are identified with sections of~$\tau$ that take the Cartan
distribution on~$\mathcal{E}$ to that on~$\T(\mathcal{E})$. Hence, if we
consider a diagram of the form
\begin{equation*}
  \xymatrix{
    &\tilde{\mathcal{E}}\ar[dl]_{\tau'}\ar[dr]^{\tau''}&\\
    \T(\mathcal{E})\ar[dr]_{\tau}&&\T(\mathcal{E})\ar[dl]^{\tau}\\
    &\mathcal{E}\rlap{,}&
  }
\end{equation*}
where~$\tau'$ and~$\tau''$ are coverings, then this B\"{a}cklund
transformation will relate symmetries of~$\mathcal{E}$ to each other. Thus,
this B\"{a}cklund transformation plays the r\^{o}le of a recursion operator
for symmetries of~$\mathcal{E}$.

\begin{example}
  \label{exmp:gjsis_nonloc:14}
  Consider the KdV equation~$u_t=6uu_x+u_{xxx}$ and two copies of its tangent
  covering with the new dependent variable~$v$ that enjoys the the additional
  equation
  \begin{equation*}
    v_t=6u_xv+6uv_x+v_{xxx}.
  \end{equation*}
  Introduce a nonlocal variable~$\tilde{v}$ by setting
  \begin{equation*}
    \tilde{v}_x=v,\qquad\tilde{v}_t=6uv+v_{xx}.
  \end{equation*}
  Thus, internal coordinates in~$\tilde{\mathcal{E}}$ are
  \begin{equation*}
    x,\  t,\ u=u_0,\ u_x=u_1,\,\dots,\,v=v_0,\ v_x=v_1,\,\dots,\,\tilde{v}.
  \end{equation*}
  Define the covering~$\tau'$ by
  \begin{equation*}
    \tau'\colon(x,t,u_k,v_k,\tilde{v})\mapsto(x,t,u_k,v_k)
  \end{equation*}
  and the covering~$\tau''$ by
  \begin{equation*}
    \tau'\colon(x,t,u_k,v_k,\tilde{v})\mapsto
    (x,t,u_k,D_x^k(v_2+4uv+2u_1\tilde{v})).
  \end{equation*}
  The B\"{a}cklund transformation obtained in such a way is the geometrical
  realization of the Lenard recursion operator~$R=D_x^2+4u+2u_1D_x^{-1}$ given
  in Equation~\eqref{eq:gjsis_jets:53}.
\end{example}

Another, less trivial example will be considered later (see
Example~\ref{exmp:gjsis_nonloc:16} below), after discussing the concept of
\emph{zero-curvature representations} (ZCR). 

Let~$\tau\colon\tilde{\mathcal{E}}\to\mathcal{E}$ be a covering. We say that
it is \emph{linear} if
\begin{enumerate}
\item $\tau$ is a vector bundle;
\item the action of vector fields~$\tilde{D}_1,\dots,\tilde{D}_n$
  on~$\mathcal{F}(\tilde{\mathcal{E}})$ preserves the subspace of fibre-wise
  linear functions.
\end{enumerate}

\begin{coordinates}
  Let~$v^1,\dots,v^r,\dots$ be local coordinates along the fibre of~$\tau$ and
  the covering be given by the total derivatives
  \begin{equation}
    \label{eq:gjsis_nonloc:26}
    \tilde{D}_i=D_i+\sum_r X_i^r\frac{\partial}{\partial v^r},\qquad
    i=1,\dots,n.
  \end{equation}
  Then the covering is linear if and only if the coefficients~$X_i^r$
  in~\eqref{eq:gjsis_nonloc:26} are of the form
  \begin{equation*}
    X_i^r=\sum_\alpha X_{i\alpha}^rv^\alpha,
  \end{equation*}
  where~$X_{i\alpha}^r$ are smooth functions on~$\mathcal{E}$. If we now
  identify the vertical terms~$X_i=\sum_r X_i^r\partial/\partial v^r$
  in~\eqref{eq:gjsis_nonloc:26} with the function-valued matrices
  \begin{equation*}
    X_i=
    \begin{pmatrix}
      X_{i1}^1&\dots&X_{i1}^n\\
      \dots&\dots&\dots\\
      X_{in}^1&\dots&X_{in}^n
    \end{pmatrix}
  \end{equation*}
  then~\eqref{eq:gjsis_nonloc:26} will be rewritten as
  \begin{equation*}
    \tilde{D}_i=D_i+X_i,\qquad i=1,\dots,n,
  \end{equation*}
  while the conditions~$[\tilde{D}_i,\tilde{D}_j]=0$ will acquire the form
  \begin{equation*}
    D_i(X_j)-D_j(X_i)+[X_i,X_j]=0,\qquad 1\le i<j\le n.
  \end{equation*}
  In other words, we arrive to the classical definition of a ZCR
  (cf.~\cite{AblowitzKaupNewellSegur:InSTFAnNP}).
\end{coordinates}

\begin{example}
  \label{exmp:gjsis_nonloc:15}
  The well known two-dimensional ZCR for the KdV equation
  (see~\cite{AblowitzKaupNewellSegur:InSTFAnNP}) is given by
  \begin{equation*}
    \tilde{D}_x=D_x+A,\qquad\tilde{D}_t=D_t+B,
  \end{equation*}
  where
  \begin{equation*}
    A=
    \begin{pmatrix}
      0&1\\
      -u+\lambda&0
    \end{pmatrix},\qquad
    B=
    \begin{pmatrix}
      -u_x&2u-4\lambda\\
      -u_{xx}-2u^2+2\lambda u+4\lambda^2&u_x
    \end{pmatrix}.
  \end{equation*}
\end{example}

\begin{example}[vacuum Einstein equations]
  \label{exmp:gjsis_nonloc:16}
  Consider the Lewis metric~$\rmd s^2=2f(x,y)\rmd x\rmd y+\sum_{j\le
    j}g_{ij}\rmd z^i\rmd z^j$ in~$\R^4$ with coordinates~$x$, $y$, $z^1$,
  and~$z^2$ (see~\cite{Lewis:SSSEqAxSGF}). The the vacuum Einstein equations read
  \begin{equation}
    \label{eq:gjsis_nonloc:27}
    (\sqrt{\det g}g_xg^{-1})_y+(\sqrt{\det g}g_yg^{-1})_x=0.
  \end{equation}
  After a re-parameterisation,~\eqref{eq:gjsis_nonloc:27} acquires the form
  \begin{equation}
    \label{eq:gjsis_nonloc:28}
    u_{xy}=\frac{u_xu_y-v_xv_y}{u}-\frac{1}{2}\frac{u_x+u_y}{x+y},\qquad
    v_{xy}=\frac{v_xu+y+u_xv_y}{u}-\frac{1}{2}\frac{v_x+v_y}{x+y}.
  \end{equation}
  B\"{a}cklund transformations and ZCR for~\eqref{eq:gjsis_nonloc:28} were
  constructed in many papers (see,
  e.g.,~\cite{BelinskiiZakharov:InEinEqMInSPTCExSS,Harrison:BTErEqGR,Maison:ArSAxSEinEqCIn,Maison:CInSAxSEinEq}). The latter is of the
  form
  \begin{equation*}
    \tilde{D}_x=D_x+A,\qquad\tilde{D}_y=D_y+B,
  \end{equation*}
  where
  \begin{equation}
    \label{eq:gjsis_nonloc:29}
    A=\frac{1}{2}
    \begin{pmatrix}
      -\frac{(\theta+1)u_x}{u}&\frac{(\theta+1)v_x}{u^2}\\
      (\theta-1)v_x,&\frac{(\theta+1)u_x}{u}
    \end{pmatrix},\qquad
    B=\frac{1}{2\theta}
    \begin{pmatrix}
      -\frac{(\theta+1)u_y}{u}&\frac{(\theta+1)v_y}{u^2}\\
      (1-\theta)v_y&\frac{(\theta+1)u_y}{u}
    \end{pmatrix};
  \end{equation}
  here~$\theta=\sqrt{(\lambda+y)(\lambda-x)}$ and~$\lambda$ is the spectral
  parameter.
  
  Using ZCR~\eqref{eq:gjsis_nonloc:29}, a three-dimensional covering
  over~$\T(\mathcal{E})$ can be constructed
  (see~\cite{Marvan:ROpVEinEqS}). Let~$U$ and~$V$ be the variables
  in~$\T(\mathcal{E})$ corresponding to~$u$ and~$v$, respectively, and~$w^1$,
  $w^2$, $w^3$ be the nonlocal variables. Then the covering is given by the
  relations
  \begin{align*}
    w_x^1&=\frac{1-\theta}{2}v_xw^2+\frac{1+\theta}{2u^2}v_xw^3-
    \frac{1+\theta}{2u}U_x+\frac{1+\theta}{2u^2}u_xU,\\
    w_x^2&=-\frac{1+\theta}{u^2}v_xw^1-\frac{1+\theta}{u}u_xw^2-
    \frac{1+\theta}{u^3}v_xU+\frac{1+\theta}{2u^2}V_x,\\
    w_x^3&=(\theta-1)v_xw^1+\frac{1+\theta}{u}u_xw^3+\frac{\theta-1}{2}V_x
  \end{align*}
  and
  \begin{align*}
    w_y^1&=\frac{\theta-1}{2\theta}v_yw^2+\frac{1+\theta}{2\theta u^2}v_yw^3+
    \frac{1+\theta}{2\theta u^2}u_yU-\frac{1+\theta}{2\theta u}U_y,\\
    w_y^2&=-\frac{1+\theta}{\theta u^2}v_yw^1-\frac{1+\theta}{\theta u}u_yw^2-
    \frac{1+\theta}{\theta u^3}v_yU+\frac{1+\theta}{2\theta u^2}V_y,\\
    w_y^3&=\frac{1-\theta}{\theta}v_yw^1+\frac{1+\theta}{\theta u}u_yw^3+
    \frac{1-\theta}{2\theta}V_y.
  \end{align*}
  This covering gives rise to a B\"{a}cklund transformation of the form
  \begin{equation*}
    \theta U'=2uw^1+U,\qquad\theta V'=-u^2w^2-w^3,
  \end{equation*}
  i.e., to a recursion operator for symmetries.
\end{example}

\begin{remark}
  Note that with an arbitrary
  covering~$\tau\colon\tilde{\mathcal{E}}\to\mathcal{E}$ one can naturally
  associate a linear
  covering~$\tau^v\colon\T^v\tilde{\mathcal{E}}\to\tilde{\mathcal{E}}$. The
  space~$\T^v\tilde{\mathcal{E}}$ is a submanifold in~$\T\tilde{\mathcal{E}}$
  and consists of tangent vectors that vanish under the action of the
  differential~$\tau_*$.

  Note also the existence of the exact sequence of coverings
  \begin{equation*}
    \xymatrix{
      0\ar[r]&\T^v\tilde{\mathcal{E}}\ar[r]\ar[dr]_{\tau^v}&
      \T\tilde{\mathcal{E}}\ar[d]^{\tilde{\tau}}\ar[r]
      &\T^s\tilde{\mathcal{E}}\ar[r]\ar[dl]^{\tau^s}&0\\
      &&\tilde{\mathcal{E}}\rlap{,}&&
    }
  \end{equation*}
  where~$\tau^s\colon\T^s\tilde{\mathcal{E}}\to\tilde{\mathcal{E}}$ is the
  quotient. The terms of this sequence possess the following characteristic
  property: integrable sections\footnote{We say that a section is integrable
    if it preserves the Cartan distributions.} of~$\tau^v$ are infinitesimal
  gauge symmetries of~$\tau$, integrable sections of~$\tilde{\tau}$ are, as it
  was mentioned above, nonlocal $\tau$-symmetries of~$\mathcal{E}$, and
  integrable sections of~$\tau^s$ are $\tau$-shadows.
\end{remark}

\section*{Concluding remarks}
\label{sec:gjsis_concl:concluding-remarks}

We described a geometrical approach to partial differential equations which
proved to be efficient both from the theoretical point of view and in a lot of
applications. Based on this approach, in particular, Hamiltonian formalism for
arbitrary normal ($2$-line) equations is constructed. On the other hand, a
number of interesting and important problems are waiting for their
solution. We intend to continue the research along the following lines:
\begin{itemize}
\item Generalisation of the Hamiltonian formalism from normal
  equations to arbitrary $p$-line ones that, in particular, include
  gauge-invariant systems.
\item Incorporation of Dirac structures into the above described scheme and
  elaboration of their computation and use.
\item Further development of the nonlocal theory and, in particular, analysis
  of differential coverings over the systems with the more than two
  independent variables and generalisation of the theory of variational
  brackets to nonlocal structures.
\end{itemize}

\section*{Acknowledgements}
This work was supported in part by the NWO-RFBR grant
047.017.015, RFBR-Consortium E.I.N.S.T.E.IN grant 09-01-92438 and
RFBR-CNRS grant 08-07-92496.

We express gratitude to the participants of our seminar at the Independent
University of Moscow (\url{http://gdeq.org}), where the topics of our paper
were discussed.  We are grateful to our colleagues whom we collaborated with
for years and especially to Paul Kersten from Twente University (The
Netherlands), Michal Marvan from the Silesian University in Opava (Czech
Republic), Raffaele Vitolo from Salento University (Italy), and Sergey Igonin
from the Utrecht University (The Netherlands).  We are also grateful to Victor
Kac for a number of useful remarks. We express our special thanks to an
anonymous reader for his/her attention to the first version of this text and
suggestions to improve the exposition.

%\bibliographystyle{xamsplain}
%\bibliography{gjsis_mrabbrev,gjsis}
\providecommand{\MR}{\relax\ifhmode\unskip\space\fi MR }
\providecommand{\MRhref}[2]{%
  \href{http://www.ams.org/mathscinet-getitem?mr=#1}{#2}
}
\providecommand{\href}[2]{#2}
\providecommand{\urlprefix}{URL }
\providecommand*{\eprint}[2][]{%
\href{http://arXiv.org/abs/#2}{\begingroup \Url{arXiv:#2}}%
}

\end{document}